\documentclass[12pt]{amsart}
\usepackage[margin=1in]{geometry}
\usepackage{standalone}
\usepackage{amsmath,amsthm,amssymb,amsfonts,leftindex,stmaryrd}
\usepackage{harpoon,float}
\usepackage{mathtools,enumitem,multirow,outlines,extarrows}
\usepackage{graphicx, subcaption}
\usepackage{grffile}
\usepackage[makeroom]{cancel}
\usepackage[dvipsnames]{xcolor}
\usepackage{tikz, tikz-cd, tikzscale}
\usetikzlibrary{spath3, knots, shapes.geometric}
\usetikzlibrary{arrows.meta,arrows}
\tikzcdset{labels={font=\everymath\expandafter{\the\everymath\displaystyle}}}
\tikzcdset{cells={font=\everymath\expandafter{\the\everymath\displaystyle}}}
\usepackage{nicematrix,blkarray}
\NiceMatrixOptions{cell-space-limits = 2pt}
\usepackage{hyperref}
\hypersetup{
    colorlinks,
    citecolor=black,
    filecolor=black,
    linkcolor=MidnightBlue,
    urlcolor=black
}
\usepackage[nameinlink,capitalize]{cleveref}
\usepackage{thmtools}
\usepackage{todonotes}

\newcommand{\mbb}{\mathbb}
\newcommand{\mbf}{\mathbf}
\newcommand{\mcal}{\mathcal}
\newcommand{\mrm}{\mathrm}
\newcommand{\msf}{\mathsf}

\newcommand{\btri}{\blacktriangle}

\newcommand{\bsq}{\blacksquare}

\newcommand{\tit}{\textit}

\newcommand{\pic}{\includegraphics}
\newcommand{\picw}[2]{\pic[width = #1\textwidth]{#2}}

\newcommand{\br}{\overline}
\newcommand{\td}{\widetilde}
\newcommand*{\defeq}{\mathrel{\vcenter{\baselineskip0.5ex \lineskiplimit0pt\hbox{\scriptsize.}\hbox{\scriptsize.}}}%
=}

\newcommand{\DMO}{\DeclareMathOperator}
\DMO{\dom}{dom}
\DMO{\cod}{cod}
\DMO{\image}{im}
\DMO{\coker}{coker}
\DMO{\tr}{tr}
\DMO{\trdeg}{trdeg}
\DMO{\adj}{adj}
\DMO{\id}{id}
\DMO{\rk}{rank}
\DMO{\graph}{graph}
\DMO{\sgn}{sgn}
\DMO{\spn}{span}
\DMO{\supp}{supp}
\DMO{\codim}{codim}
\DMO{\Hom}{Hom}
\DMO{\End}{End}
\DMO{\Alt}{Alt}
\DMO{\cha}{char}
\DMO{\diag}{diag}
\DMO{\interior}{int}
\DMO{\Tor}{Tor}
\DMO{\Ext}{Ext}
\DMO{\Spec}{Spec}
\DMO{\Proj}{Proj}
\DMO{\rank}{rk}
\DMO{\ann}{ann}
\DMO{\Frac}{Frac}
\DMO{\Ass}{Ass}
\DMO{\zd}{zd}
\DMO{\nil}{nil}
\DMO*{\colim}{colim}
\DMO{\Bl}{Bl}
\DMO{\cd}{cd}
\DMO{\ob}{ob}
\DMO{\cone}{cone}
\DMO{\Tot}{Tot}
\DMO{\Coder}{Coder}
\DMO{\Coeff}{Coeff}
\DMO{\Mat}{Mat}
\newcommand{\bx}{\mathbin{\square}}
\DMO{\rot}{rot}

\renewcommand{\k}{\Bbbk}
\newcommand{\R}{\mbb{R}}

\newcommand{\Z}{\mbb{Z}}

\newcommand{\mcalA}{\mcal{A}}
\newcommand{\mcalB}{\mcal{B}}

\newcommand{\mcalI}{\mcal{I}}
\newcommand{\mcalF}{\mcal{F}}
\newcommand{\mcalC}{\mcal{C}}
\newcommand{\mcalM}{\mcal{M}}
\newcommand{\mcalN}{\mcal{N}}
\newcommand{\mcalQ}{\mcal{Q}}

\newcommand{\mcalV}{\mcal{V}}
\newcommand{\mcalL}{\mcal{L}}
\newcommand{\mcalW}{\mcal{W}}
\newcommand{\mcalP}{\mcal{P}}

\newcommand{\mbfT}{\mbf{T}}
\newcommand{\mbfS}{\mbf{S}}
\newcommand{\mbfR}{\mbf{R}}
\newcommand{\mbfQ}{\mbf{Q}}
\newcommand{\mbfP}{\mbf{P}}
\newcommand{\Aug}{\mcalA ug}
\newcommand{\del}{\partial}

\newcommand{\ep}{\epsilon}
\newcommand{\tep}{\tilde{\ep}}
\newcommand{\al}{\alpha}
\newcommand{\be}{\beta}
\newcommand{\de}{\delta}
\newcommand{\la}{\lambda}
\newcommand{\LA}{\Lambda}
\newcommand{\DE}{\Delta}
\newcommand{\tht}{\theta}

\newcommand{\inn}{\subset}

\newcommand{\cont}{\supset}

\newcommand{\bbra}[1]{\llbracket #1 \rrbracket}
\newcommand{\agl}[1]{\langle #1 \rangle}

\newcommand{\flr}[1]{\lfloor #1 \rfloor}

\def\mch#1#2{\ensuremath{\left(\kern-.3em\left(\genfrac{}{}{0pt}{}{#1}{#2}\right)\kern-.3em\right)}}
\newcommand{\mono}{\hookrightarrow}

\newcommand{\xrar}{\xrightarrow}

\newcommand{\iso}{\xlongrightarrow{\sim}}

\newcommand{\dar}{\downarrow}
\newcommand{\uar}{\uparrow}

\newcommand{\bimod}[1]{#1\text{-}\msf{bimod}}

\theoremstyle{definition}
\newtheorem{thm}{Theorem}[section]
\newtheorem{coro}[thm]{Corollary}
\newtheorem{defn}[thm]{Definition}
\newtheorem{prop}[thm]{Proposition}
\newtheorem{lemma}[thm]{Lemma}
\newtheorem{exmp}[thm]{Example}
\newtheorem{deflem}[thm]{Definition/Lemma}
\newtheorem{conj}{Conjecture}
\newtheorem*{ack}{Acknowledgements}
\theoremstyle{remark}
\newtheorem{rmk}[thm]{Remark}
\numberwithin{equation}{subsection}

\begin{document}
\title{\texorpdfstring{$A_\infty$}{A-infinity} Sabloff Duality via the LSFT Algebra}
\author{Zhenyi Chen}
\begin{abstract}
    We use Ng's LSFT algebra to upgrade Sabloff duality of Legendrian knots to a quasi-isomorphism of $A_\infty$ bimodules over the positive augmentation category $\Aug_+$. We also extend the Ekholm-Etnyre-Sabloff exact sequence to an exact sequence of $\Aug_+$-bimodules, using a quotient category $\mcalC$ of short Reeb chords. In addition, we define curved augmentations of the LSFT algebra and show that they can be used to construct a homotopy inverse of the $A_\infty$ Sabloff map, together with all higher homotopies. The above results suggest a conjectural recipe for an explicit weak relative Calabi-Yau structure on the quotient $A_\infty$ functor $\pi:\Aug_+\to \mcalC$.
\end{abstract}
\maketitle
\tableofcontents

\newpage

\section{Introduction}
The study of Legendrian links was revolutionized by the introduction of Legendrian contact homology, which encodes holomorphic disk counts in the symplectization of the ambient contact manifold. Eliashberg laid out the ideas of contact homology and its Legendrian version in \cite{Eli98}. Independently, Chekanov \cite{Che02} gave the first rigorous construction combinatorially for Legendrian links $\LA$ in the standard contact $\R^3$.
The resulting differential graded algebra (dga) is known as the \tit{Chekanov-Elisashberg algebra}, denoted by $\mcalA = \mcalA(\LA)$. 

As part of symplectic field theory \cite{EGH10}, the Chekanov-Eliashberg algebra is functorial with respect to exact Lagrangian cobordisms \cite{EHK16}\cite{Kar20}. In particular, every exact Lagrangian filling $L$ of $\LA$ induces an augmentation $\ep_L$ of $\mcalA$. Chekanov used augmentations to truncate $\mcalA$, producing \tit{linearized contact homology}, which is a very effective Legendrian invariant. 
This truncation procedure was later extended by \cite{CEKSW}\cite{BC14}\cite{NRSSZ} to form two $A_\infty$ categories, $\Aug_-$ and $\Aug_+$, whose objects are augmentations and whose hom-complexes are denoted by $\Hom_-$ and $\Hom_+$, respectively. These complexes relate to cohomologies of a filling $L$ via the isomorphisms \cite{DR16}\cite{NRSSZ}\[H^*\Hom_+(\ep_L,\ep_L) \simeq H^*(L)\,,\quad H^*\Hom_-(\ep_L,\ep_L) \simeq H^*(L,\LA)\,.\]

Moreover, there is an exact triangle of complexes, known as the \tit{Ekholm-Etnyre-Sabloff exact sequence} \cite{EES09}\cite{BC14}\cite{NRSSZ} \begin{equation}\label{EES sequence}\Hom_-(\ep,\ep')\to \Hom_+(\ep,\ep')\to \mcalC(\ep,\ep')\to{}\end{equation} for any augmentations $\ep,\ep'$, where $\mcalC(\ep,\ep')$ is spanned by a set of short Reeb chords corresponding to critical points of an auxiliary Morse function on $\LA$. When $\ep = \ep' = \ep_L$, this exact triangle reproduces the long exact sequence in cohomology: \[H^*(L,\LA)\to H^*(L)\to H^*(\LA)\to{}\,.\] 

This parallel goes one step further: there is a quasi-isomorphism, known as the \tit{Sabloff duality map} \cite{Sab06}\cite{EES09}\cite{BC14}\cite{NRSSZ}
\begin{equation}\label{linear Sabloff}\Hom_-(\ep,\ep')^\dagger[-2]\iso \Hom_+(\ep',\ep)\end{equation} which recovers Poincar\'{e} duality \[H_{2-*}(L,\LA)\simeq H^*(L)\] when $\ep = \ep' = \ep_L$. The superscript ${}^\dagger$ denotes the linear dual with negated grading.

The following two theorems upgrade the above relations to the level of $A_\infty$ bimodules over $\Aug_+$.
\begin{thm}[cf. {\Cref{Q+- triangle}}]
\label{intro A infty EES}
 Let $\LA$ be a Legendrian link with marked points and any Morse function. We define the \tit{negative augmentation bimodule} $\mcalQ_-$, whose underlying complex is $\Hom_-$ (\Cref{Q- def}), and a quotient $\pi:\Aug_+\to\mcalC$ called the \tit{short chord category} (\Cref{short chord category}) . There is a short exact sequence of $\Aug_+$-bimodules consisting of strict morphisms
\[
    \begin{tikzcd}
        0\rar & \mcalQ_-\rar &\Aug_+\rar{\pi} & \pi^*\mcalC\rar & 0
    \end{tikzcd}
    \] whose underlying sequence of $\k$-quivers is split, extending \eqref{EES sequence}.	
\end{thm}
\begin{thm}[cf. {\Cref{A infty Sabloff thm}}]
    \label{intro A infty Sabloff}
    Let $\LA$ be a Legendrian knot with a single marked point and any Morse function. We construct a quasi-isomorphism of $\Aug_+$-bimodules
    \[\mcalQ_-^\dagger[-1]\iso \Aug_+[1]\] called the \tit{$A_\infty$ Sabloff map}, extending \eqref{linear Sabloff}. The superscript ${}^\dagger$ denotes the linear dual bimodule (\Cref{linear dual bimodule def}).
\end{thm}
The key tool in constructing the $A_\infty$ Sabloff map is Ng's LSFT algebra, and we believe the same method can be generalized to cover Legendrian links with multiple marked points. In their upcoming work \cite{MS25}, Ma and Sabloff gives a different proof of these two theorems over $\Z/2$. 

The LSFT algebra of a Legendrian knot, introduced in \cite{Ng10}, is a filtered curved dga that counts holomorphic disks with multiple positive punctures, and its zeroth associated graded is the Chekanov-Eliashberg algebra. For multi-component Legendrian links, however, the LSFT differential does not square to something simple, as investigated in \cite{Ng23}. In this paper, we define a variant of the LSFT algebra that generalizes the definition for Legendrian knots while still being a curved dga:
\begin{thm}[cf. {\Cref{curvature formula}}]
    \label{intro LSFT}
    The \tit{composable LSFT algebra} of any pointed Legendrian link, defined in \Cref{composable LSFT algebra}, is a filtered curved dga.
\end{thm}

Sabloff duality is governed by disks with two positive punctures, which leaves one wondering: what is the role played by disks with more positive punctures? To investigate this question, we study lifts of augmentations to \tit{curved augmentations} of the LSFT algebra (see \Cref{curved augmentations}). It turns out such lifts supply an explicit homotopy inverse to the $A_\infty$ Sabloff map, as well as all higher homotopies (see \Cref{homotopy equations lemma}). When $\rot(\LA) = 0$, this result can be phrased as follows:

\begin{thm}[cf. {\Cref{higher homotopy thm}}]
\label{intro homotopies}
Let $J$ be the nerve of the groupoid $0\leftrightarrows 1$. If every augmentation of $\LA$ lifts to a $1$-curved augmentation, then we construct an explicit map of $\infty$-categories from $J$ to the dg nerve of the dg category of $\Aug_+$-bimodules, which sends the edge $1\to 0$ to the $A_\infty$ Sabloff map.
\end{thm}

   
This result demonstrates how disks with multiple positive punctures provide further homotopy coherence to the algebraic structure of Legendrian contact homology. In a similar spirit, Ng \cite{Ng23} used the LSFT operations to construct an $L_\infty$ structure on the commutatized Chekanov-Eliashberg dga. 

Conjecturally, the homotopy inverse provided by \Cref{intro homotopies} is part of an explicit weak Calabi-Yau structure on the functor $\pi:\Aug\to\mcalC$. We discuss partial results and related works in this direction in \Cref{CY structures}. 

We conclude the introduction by explaining the logical dependence among the sections. \Cref{background} contains background materials that can be skipped depending on familiarity of the reader. \Cref{composable LSFT algebra} can be read by itself. \Cref{aug+ bimods} requires \Cref{CE subsection}, \Cref{A infty background}, and \Cref{gradings}. \Cref{A infty Sabloff and higher homotopies section} requires almost everything preceding it.

\begin{ack}
    I am very grateful to my advisor Eric Zaslow for his invaluable help and encouragement. I worked in parallel with Jiajie Ma and Josh Sabloff, and I would like to thank them for sharing their insights and results. I thank Lenny Ng, Alex Karapetyan, Wenyuan Li, Ezra Getzler, and Boris Tsygan for very helpful discussions. Special thanks to Joanna Ciatti for her support. This project was partially supported by the NSF grant DMS-2104087.
\end{ack}

\section{Background}
\label{background}
\subsection{Chekanov-Eliashberg Algebra}
\label{CE subsection}
In this subsection, we briefly review the Chekanov-Eliashberg dga of Legendrian links in $\R^3$. 
For a more thorough discussion, see for example \cite{EN23}. We follow the sign convention of \cite{Ng23}, which is different but equivalent to the more common convention of \cite{ENS02} (cf. \cite[Remark 2.17]{Ng23}).

Equip $\R^3 = J^1(\R)$ with the standard contact form $\alpha = dz - ydx$. Let $\Lambda\inn\R^3$ be an oriented Legendrian link with generic Lagrangian projection $\Lambda_{xy} \defeq \pi_{xy}(\Lambda)$. Label the Reeb chords of $\LA$ by $a_1,\cdots,a_n$, which correspond to the crossings in $\LA_{xy}$. Denote the starting/ending point of $a_j$ by $a_j^-,a_j^+$, respectively. In the Lagrangian projection, $a_j^-$ is the under-crossing, and $a_j^+$ is the over-crossing. Place marked points $\star_1,\cdots,\star_m$ on $\Lambda$, such that $\pi_{xy}(\star_i)$ avoid the crossings and each component of $\Lambda$ contains at least one marked point. 

Fix a commutative ring $\k$. The \tit{Chekanov-Eliashberg algebra} $\mcalA = \mcalA(\Lambda)$ is a dga whose homotopy type is a Legendrian isotopy invariant of $\Lambda$, defined as follows. To each Reeb chord $a_j$, associate a variable $q_j$; to each marked point $\star_i$, associate two variables $t_i, t_i^{-1}$. Denote 
$\mbfQ = \{q_1,\cdots,q_n\}$, $\mbfT = \{t^{\pm 1}_1,\cdots,t^{\pm 1}_m\}$, and \[\mbfS^0 = \mbfQ\cup \mbfT = \{q_1,\cdots,q_n, t^{\pm 1}_1,\cdots,t^{\pm 1}_m\}\,.\]

Gradings of the generators take values in $\Z$ and require a Maslov potential to define. We postpone this endeavor to \Cref{gradings}, where we define gradings for the LSFT generators, which include all the Chekanov-Eliashberg generators. Define $\mcalA$ to be the tensor algebra of the graded free $\k$-module with basis $\mbfS^0$, modulo the relations $t_it^{-1}_i = t_i^{-1}t_i = 1$ for $1\leq i\leq m$. 

\begin{figure}[H]
\centering
  \begin{tikzpicture}
      \draw[ultra thick] (1,-1) -- (-1,1);
      \draw[ultra thick] (1,1) -- (.2,.2);
      \draw[ultra thick] (-1,-1) -- (-.2,-.2);
      \node (p) at (.5,0) {$+$};
      \node (p) at (-.5,0) {$+$};
      \node (q) at (0,.5) {$-$};
      \node (q) at (0,-.5) {$-$};
  \end{tikzpicture}
\caption{Reeb signs}
\label{Reeb signs}
\end{figure}
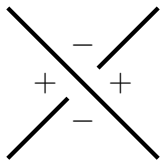

In the Lagrangian projection, the differential is given by counting boundary-punctured oriented topological disks in $\R^2$ whose boundary lies on $\LA$ and boundary punctures limit to crossings. Depending on which corner the disk covers at a crossing, each puncture is assigned a \tit{Reeb sign}, depicted in \Cref{Reeb signs}.  

The Chekanov-Eliashberg differential only counts disks with a single positive puncture. Let $\DE_1(\LA;p_j)$ be the set of disks with a unique positive puncture, which limits to the crossing $a_j$. See \Cref{single positive puncture disks rmk} for the precise definition. Given any $\DE\in\DE_1(\LA;p_j)$, let $w'(\DE)\in\mcalA$ be the word that successively record generators $t_i^{\pm 1}$'s and $q_j$'s as $\del\DE$ traverses through marked points $\star_i$'s and as it encounters crossings $a_j$'s at its punctures, respectively. See \eqref{w'} for the exact definition.

\begin{defn}
\label{CE diff def}
    The \tit{Chekanov-Eliashberg differential} $\del:\mcalA\to\mcalA$ is defined by
    \begin{align*}
		\del  q_j &= \sum_{\DE\in \DE_1(\LA;p_j)}\sgn_\del(\DE)w'(\DE)\\
		\del  t_i &= \del   t_i^{-1} = 0
	\end{align*}
where $\sgn_\del(\DE)$ is $-1$ to the number of shaded corners of $\DE$, according to \Cref{Orientation signs for CE}.
\end{defn}
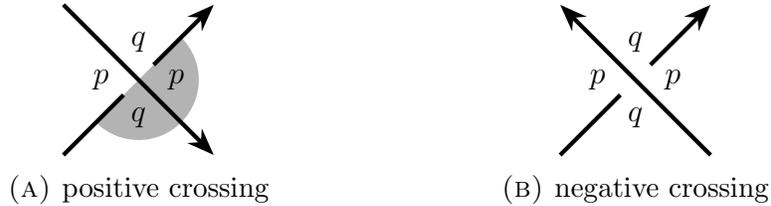
\begin{figure}[H]
\centering
\begin{subfigure}{.4\textwidth}
  \centering
  \begin{tikzpicture}
      \fill[fill=gray!60] (0,0) --  (45:.8) arc(45:-135:.8) -- cycle;
      \draw[ultra thick, Stealth-] (1,-1) -- (-1,1);
      \draw[ultra thick, Stealth-] (1,1) -- (.2,.2);
      \draw[ultra thick] (-1,-1) -- (-.2,-.2);
      \node (p) at (.5,0) {$p$};
      \node (p) at (-.5,0) {$p$};
      \node (q) at (0,.5) {$q$};
      \node (q) at (0,-.5) {$q$};
  \end{tikzpicture}
  \caption{positive crossing}
  \label{Orientation signs positive crossings p}
\end{subfigure}%
\begin{subfigure}{.4\textwidth}
  \centering
  \begin{tikzpicture}
      \draw[ultra thick, -Stealth] (1,-1) -- (-1,1);
      \draw[ultra thick, Stealth-] (1,1) -- (.2,.2);
      \draw[ultra thick] (-1,-1) -- (-.2,-.2);
      \node (p) at (.5,0) {$p$};
      \node (p) at (-.5,0) {$p$};
      \node (q) at (0,.5) {$q$};
      \node (q) at (0,-.5) {$q$};
  \end{tikzpicture}
  \caption{negative crossing}
  \label{Orientation signs negative crossings p}
\end{subfigure}
\caption{Orientation signs for the Chekanov-Eliashberg differential $\del$}
\label{Orientation signs for CE}
\end{figure}
Let $\rot(\LA)$ denote the rotation number of $\LA$ (see \Cref{rot num def}). Pick a Maslov potential so that $|t_i| \equiv 0 \mod 2\rot(\LA)$ for each $i$, which exists by \Cref{rotation number base point grading}. This choice of grading is necessary to make the following definition nonempty, since $t_i$'s are invertible.
\begin{defn}
    An \tit{augmentation} of $\LA$ is a map $\ep: \mcalA\to \k$ of differential $\Z/2\rot(\LA)$-graded algebras, where $\k$ is considered a dga concentrated in degree $0$ and equipped with the zero differential. 
\end{defn}

\subsection{\texorpdfstring{$A_\infty$}{A-infinity} Categories and Bimodules}
\label{A infty background}
In this subsection, we record the necessary definitions about $A_\infty$ categories and bimodules, following \cite{Tra08}\cite{LM08}. Readers familiar with $A_\infty$ structures can skip this subsection.

Fix a commutative ring $\k$. All definitions in this subsection are stated in terms of $\k$-modules, but they work more generally for $\k$-bimodules verbatim, which we will need later when discussing composable $m$-copy algebras.

\begin{rmk}
    In this subsection, ``graded'' means $\Z$-graded, but later when dealing with Legendrian links of non-zero rotation number, we will also allow gradings to take values in $\Z/n$ for some $n\in\Z$. Almost everything below carry over to the $\Z/n$-graded setting, except for the pretriangulated dg structure on the category of $A_\infty$-bimodules. See \Cref{triangulated rmk}.
\end{rmk}





We first define the additive objects underlying $A_\infty$ categories and bimodules, which one should think of as graded $\k$-linear categories without compositions.
\begin{defn}
    A \tit{graded $\k$-linear quiver} (or simply, \tit{$\k$-quiver}) $\mcalV$ consists of
    \begin{itemize}
    \item a set of objects, denoted by $\ob \mcalV$
    \item a graded $\k$-module $\mcalV(\ep,\ep')$, for each pair $\ep,\ep'\in\ob\mcalV$.
\end{itemize}
A \tit{degree $k$ map} of $\k$-quivers $f: \mcalV\to \mcalW$ consists of
\begin{itemize}
    \item a set map $f: \ob \mcalV\to \ob\mcalW$
    \item a degree $k$ map $f^{\ep_1,\ep_2}:\mcalV(\ep_1,\ep_2)\to \mcalW(f(\ep_1),f(\ep_2))$, for each pair $\ep_1,\ep_2\in\ob\mcalV$. 
\end{itemize}
\end{defn}

Fix a set $Z$. A \tit{$\k$-quiver on $Z$} is a $\k$-quiver with object set $Z$. The class of $\k$-quivers on $Z$ underlies a monoidal abelian category, which we now describe. A morphism of $\k$-quivers on $Z$ is a degree $0$ map of $\k$-quivers that equals to $1_Z$ on objects. Direct sums, kernels, and cokernels are taken object-wise. The tensor product of two $\k$-quivers $\mcalV,\mcalW$ on $Z$ is given by \[(\mcalV\otimes_Z \mcalW)(\ep_1,\ep_2)\defeq \bigoplus_{\ep\in Z}\mcalV(\ep,\ep_2)\otimes\mcalW(\ep_1,\ep)\,.\](When the set of objects is clear, we drop the subscript in $\otimes_Z$.) The monoidal unit $\k_Z$ is given by $\k_Z(\ep_1,\ep_2) = \k$ if $\ep_1=\ep_2$ and zero otherwise. 
\begin{defn}
    A \tit{dg $\k$-quiver} is a $\k$-quiver $\mcalV$ equipped with a degree $1$ map $d:\mcalV\to\mcalV$ fixing objects and satisfying $d^2 = 0$. The \tit{cohomology} of a dg $\k$-quiver is the $\k$-quiver given by taking cohomology object-wise. 
    
    A \tit{map of dg $\k$-quivers} is a degree $0$ map of $\k$-quivers that commutes with the differentials. Such a map is a \tit{quasi-isomorphism} if the induced map on cohomology is an isomorphism of $\k$-quivers.
\end{defn}

\begin{defn}
    A \tit{(counital) $\k$-cocategory} is a $\k$-quiver $\mcalB$ equipped with a coassociative comultiplication $\DE:\mcalB\to\mcalB\otimes_{\ob\mcalB} \mcalB$, and a counit $\mcalB\to \k_{\ob\mcalB}$. 
    
    A \tit{map of $\k$-cocategories} $f:\mcalB\to\mcalB'$ is a degree $0$ map of the underlying $\k$-quivers that satisfies $(f\otimes f)\circ \DE = \DE'\circ f$ and commutes with the counits. (Note that the set map $f: \ob\mcalB\to\ob\mcalB'$ induces a map of $\k$-quivers $\k_f:\k_{\ob\mcalB}\to\k_{\ob\mcalB'}$.)
\end{defn}
    
\begin{defn}A \tit{coaugmentation} of a $\k$-cocategory $\mcalB$ is a map of $\k$-cocategories $\k_{\ob\mcalB}\to \mcalB$ such that the composition $\k_{\ob\mcalB}\to \mcalB\to \k_{\ob\mcalB}$ with the counit is the identity. (Note that $\k_Z$ is canonically a $\k$-cocategory.) A map of coaugmented $\k$-cocategories is required to commute with the coaugmentations.
    \end{defn}
 Here is the most important example of a cocategory for this paper: the \tit{tensor (coaugmented) cocategory} of a $\k$-quiver $\mcalV$ is $T\mcalV \defeq \bigoplus_{k\geq 0} \mcalV^{\otimes k}$, where $\otimes = \otimes_{\ob \mcalV}$. The comultiplication is given by
    \begin{equation}
        \label{tensor comultiplication}
        \Delta(v_1\cdots v_n) = \sum_{i=0}^n v_1\cdots v_i\otimes v_{i+1}\cdots v_n\,.
    \end{equation} The counit $T\mcalV\to \k_{\ob \mcalV}$ is the projection to the zeroth tensor power, and the coaugmentation $\k_{\ob \mcalV}\to T\mcalV$ is the inclusion of the zeroth tensor power. 
    For convenience, we write \[\mcalV(\ep_1,\cdots,\ep_{k+1}) \defeq \mcalV(\ep_k, \ep_{k+1})\otimes\cdots\otimes \mcalV(\ep_1,\ep_2)\]and note that $\mcalV^{\otimes k}(\ep_1,\ep_{k+1}) = \bigoplus_{\ep_2,\cdots,\ep_{k}\in \ob\mcalV}\mcalV(\ep_1, \ep_2,\cdots,\ep_k,\ep_{k+1})$.
    
    Let $\mcalV[1]$ denote the (cohomologically graded) \tit{suspension} (or \tit{shift}) of $\mcalV$, given by \[\mcalV[1]^n \defeq \mcalV^{n+1}\,.\] Write $s: \mcalV\to \mcalV[1]$ for the obvious degree $-1$ isomorphism. Define the \tit{reduced degree} \[\|v\|\defeq |sv| = |v| - 1\,.\]
    The \tit{bar (coaugmented) cocategory} of $V$ is the tensor cocategory $B\mcalV\defeq  T(\mcalV[1])$ of its shift. We write elements $sv_1\otimes \cdots\otimes sv_n\in B\mcalV$ in the bar notation $[v_1|\cdots|v_n]$ of Eilenberg and MacLane.
\begin{defn}
    A \tit{coderivation} on  a $\k$-cocategory $\mcalB$ is a map $b:\mcalB\to\mcalB$ fixing objects and satisfying the \tit{co-Leibniz rule}: $\DE\circ b = (b\otimes 1+1\otimes b)\circ\DE$. A \tit{codifferential} is a degree $1$ coderivation satisfying $b^2 = 0$.
    A \tit{dg $\k$-cocategory} is a $\k$-cocategory equipped with a codifferential.
\end{defn}
\begin{defn}    
    An \tit{$A_\infty$ category (over $\k$)} is a $\k$-quiver $\mcalA$, equipped with a codifferential $b$ on $B\mcalA$, satisfying $b(1) = 0$. The coaugmented dg $\k$-cocategory $(B\mcalA,b)$ is called the \tit{bar complex} of $\mcalA$. 
    An \tit{$A_\infty$ algebra} is an $A_\infty$ category with one object.
\end{defn}    

As observed in \cite{Sta63}\cite{GJ90}, a coderivation $b$ on $B\mcalA$ satisfying $b(1) = 0$ is equivalent to the collection of degree $1$ \tit{components} 
\[
\begin{tikzcd}
    b_k: \mcalA[1]^{\otimes k}\rar{\iota} &B\mcalA\rar{b} &B\mcalA\rar{\pi} &\mcalA[1]
\end{tikzcd}
\] where $\iota$ is inclusion, $\pi$ is projection, and $k \geq 1$. (If $b_k = 0$ for $k > 2$, we say $\mcalA$ is a \tit{dg category}.) The requirement $b^2 = 0$ is then equivalent to the relations
\begin{equation}
\label{A_infty relations}
    \sum_{i+j+k = n} b_{i+1+k}(1^{\otimes i}\otimes b_j\otimes 1^{\otimes k}) = 0
\end{equation}
for all $n \geq 0$. We can also get rid of the shift and define the \tit{un-suspended components} \[\beta_k \defeq s^{-1}\circ b_k\circ s^{\otimes k}: \mcalA^{\otimes k}\to \mcalA\] of degree $2-k$. This introduces Koszul signs into the above relations (and even more signs when elements are inserted). For later reference, we have
\[\beta_k(a_k,\cdots, a_1) = (-1)^{(k-1)|a_k| + \cdots + |a_2|}s^{-1}b_k[a_k|\cdots|a_1]\,.\]

The components $b_k$ decompose further into maps of graded $\k$-modules
\[b_k^{\ep_1,\cdots,\ep_{k+1}}:\mcalA[1](\ep_1,\cdots,\ep_{k+1})\to \mcalA[1](\ep_1,\ep_{k+1})\]
for $\ep_1,\cdots,\ep_{k+1}\in\ob\mcalA$. In particular, the relations \eqref{A_infty relations} only involve finitely many objects at a time. 


\begin{defn}
\label{A_infty functor}
    Suppose $(\mcalA, b)$ and $(\mcalA',b')$ are  $A_\infty$ categories. An \tit{$A_\infty$ functor} $f: \mcalA\to \mcalA'$ is a map of coaugmented dg cocategories $f: B\mcalA\to B\mcalA'$.
\end{defn}
    Equivalently, the data consists of degree $0$ maps of $\k$-quivers
    $f_k:\mcalA[1]^{\otimes k}\to \mcalA'[1]$
    for $k\geq 1$, which agree on objects and satisfy
  \begin{equation}
  \label{A infty functor relations}
       \sum_{i+j+k = n} f_{i+1+k}(1^{\otimes i}\otimes b_j\otimes 1^{\otimes k}) = \sum_{m\geq 0, i_1+\cdots +i_m=n} b'_{m}(f_{i_1}\otimes \cdots \otimes f_{i_m})\,.
    \end{equation} for each $n\geq 1$. 

\begin{defn} An $A_\infty$ functor $f$ is \tit{strict} if $f_k = 0$ for $k > 1$.
\end{defn}
If $f$ is strict, then the relations \eqref{A infty functor relations} simplify to
\begin{equation}
    \label{strict A infty functor relation}
    f_1\circ b_n
       = b'_n\circ f_1^{\otimes n}\,.
\end{equation}


\begin{defn}
    An $A_\infty$ category $(\mcalA,b)$ is \tit{strictly unital} if there is a degree $0$ map $e:\k_{\ob\mcalA}\to\mcalA$ of $\k$-quivers on $\ob\mcalA$, satisfying
    \[\beta_1\circ e = 0\,,\ \beta_2(e\otimes 1_\mcalA) = \beta_2(1_\mcalA\otimes e) = 1_\mcalA\,,\ \beta_{i+1+j}(1_\mcalA^{\otimes i}\otimes e \otimes 1_\mcalA^{\otimes j}) = 0\]
    for any $i+j > 1$.
\end{defn}

We now move on to $A_\infty$ bimodules.
\begin{defn}
    Let $\mcalB$ be a cocategory. A \tit{(counital) $\mcalB$-bicomodule} is a $\k$-quiver $\mcalM$ on $\ob\mcalB$, equipped with counital coassociative coactions $\DE_L: \mcalM\to \mcalB\otimes \mcalM$ and $\DE_R: \mcalM\to \mcalM\otimes \mcalB$.
\end{defn}

Let $\mcalV,\mcalW$ be $\k$-quivers on $Z$. All bicomodules appearing in this paper are of the following form: define a bicomodule $T^\mcalW \mcalV \defeq T\mcalV\otimes \mcalW\otimes T\mcalV$ over the coalgebra $T\mcalV$, whose left and right actions $\Delta_L: T^\mcalW \mcalV\to T\mcalV\otimes  T^\mcalW \mcalV$ and $\Delta_R: T^\mcalW \mcalV\to T^\mcalW \mcalV\otimes T\mcalV$ are given by
\[\Delta_L(v_1\cdots v_k\,w\,v_{k+1}\cdots v_n) = \sum_{i=0}^k v_1\cdots v_i\otimes v_{i+1}\cdots v_k\,w\,v_{k+1}\cdots v_n\,,\]
\[\Delta_R(v_1\cdots v_k\,w\,v_{k+1}\cdots v_n) = \sum_{i=k}^n v_1\cdots v_k\,w\,v_{k+1}\cdots v_i\otimes v_{i+1}\cdots v_n\,.\] Denote $B^\mcalW\mcalV \defeq T^{\mcalW[1]}\mcalV[1]$, which is a bicomodule over $B\mcalV$.


\begin{defn}
    Let $\mcalB$ be a \tit{dg cocategory}. A \tit{dg $\mcalB$-bicomodule} is a $\mcalB$-bicomodule $\mcalV$ equipped with a \tit{codifferential}: a degree $1$ map $b^\mcalV:\mcalV\to\mcalV$ of $\k$-quivers on $\ob\mcalV$, satisfying $(b^\mcalV)^2 = 0$ and the left and right co-Leibniz rules.   
\end{defn}

\begin{defn}
\label{A infty bimodule def} 
    Suppose $(\mcalA,b)$ is an $A_\infty$ category. An \tit{$A_\infty$ bimodule} over $\mcalA$ (or simply, \tit{$\mcalA$-bimodule}) is a $\k$-quiver $\mcalM$ on $\ob \mcalA$, equipped with a codifferential $b^\mcalM$ on the $B\mcalA$-bicomodule $B^\mcalM\mcalA$. We call the dg $B\mcalA$-bicomodule $(B^\mcalM\mcalA, b^\mcalM)$ the \tit{bar complex} of $\mcalM$.

\end{defn}

Equivalently \cite{Tra08}, an $\mcalA$-bimodule structure on $\mcalM$ consists of degree $1$ maps \[b^\mcalM_{i|j}: \mcalA[1]^{\otimes i}\otimes \mcalM[1]\otimes \mcalA[1]^{\otimes j}\to \mcalM[1]\] of $\k$-quivers on $\ob \mcalA$, for $i,j\geq 0$, satisfying the \tit{($m,n$)-th relation}
\begin{equation}
\label{A_infty bimodule relations}
\begin{aligned}
    &\sum_{i+j = m, k+\ell = n} b_{i|\ell}^\mcalM\big(1^{\otimes i}\otimes b^\mcalM_{j|k}\otimes 1^{\otimes \ell}\big)\\
    + &\sum_{i+j+k = m} b^\mcalM_{i+1+k|n}\big(1^{\otimes i}\otimes b_j\otimes 1^{\otimes k}\otimes 1_{\mcalM[1]}\otimes 1^{\otimes n}\big)\\
    + &\sum_{i+j+k = n} b^\mcalM_{m|i+1+k}\big(1^{\otimes m}\otimes 1_{\mcalM[1]}\otimes 1^{\otimes i}\otimes b_j\otimes 1^{\otimes k}\big) = 0
\end{aligned}
\end{equation}
for every $m,n\geq 0$.
The \tit{un-suspended components} \[\beta^\mcalM_{i|j} \defeq s^{-1}\circ b^\mcalM_{i|j}\circ s^{\otimes (i+1+j)}: \mcalA^{\otimes i}\otimes \mcalM\otimes \mcalA^{\otimes j}\to \mcalM\] have degree $1-i-j$, for $i,j\geq 0$, and they satisfy a similar set of relations with Koszul signs. For later reference, we have
\begin{align*}&\quad\,\beta^\mcalM_{i|j}(a_i',\cdots, a_1', x, a_j,\cdots, a_1)\\ &= (-1)^{(i+j)|a_i| + \cdots + (j+1)|a_1| + j|x| + (j-1)|a_j'|+\cdots +|a_2|}s^{-1}b_{i|j}^\mcalM[a_i|\cdots| a_1| x| a_j'|\cdots| a_1']\,.\end{align*}

The components $b_{i|j}^\mcalM$ decompose into maps of graded $\k$-modules
 \begin{align*}(b_{i|j}^\mcalM)^{\ep_1,\cdots,\ep_{j+1}}_{\ep_1',\cdots,\ep_{i+1}'}:\mcalA[1](\ep_1',\cdots,\ep_{i+1}')\otimes \mcalM[1](\ep_{j+1},\ep_1')\otimes \mcalA[1](\ep_1,\cdots,\ep_{j+1})\to \mcalM[1](\ep_1,\ep_{i+1}')
    \end{align*}  
for $\ep_1,\cdots,\ep_{j+1},\ep_1',\cdots,\ep_{i+1}'\in \ob\mcalA$. Hence, the $A_\infty$ relations \eqref{A_infty bimodule relations} only involve finitely many objects at a time.

Here is an example of a bimodule that will be important to us in this paper.

\begin{defn}
The \tit{diagonal bimodule} of an $A_\infty$ category $\mcalA$ is the $\k$-quiver $\mcalA$ equipped with codifferential $b^\mcalA$, given by
\begin{equation}
\label{diag bimod formula}(b^\mcalA_{i|j})^{\ep_1,\cdots,\ep_{j+1}}_{\ep_1',\cdots,\ep_{i+1}'}[a_i'|\cdots |a_1'|a|a_j|\cdots| a_1] = b_{i+1+j}^{\ep_1,\cdots,\ep_{j+1},\ep_1',\cdots,\ep_{i+1}'}[a_i'|\cdots |a_1'|a|a_j|\cdots| a_1]\end{equation}
for any $i,j\geq 0$, $\ep_1,\cdots,\ep_{j+1},\ep_1',\cdots,\ep_{i+1}'\in\ob\mcalA$.
\end{defn}
We introduce the following notion of sub-bimodules: 
\begin{defn}
A \tit{$\k$-subquiver} $\mcalW$ of a $\k$-quiver $\mcalV$ consists of 
\begin{itemize}
	\item a subset $\ob\mcalW\inn\ob\mcalV$
	\item a graded $\k$-submodule $\mcalW(\ep,\ep')\inn \mcalV(\ep,\ep')$, for each $\ep,\ep'\in\ob\mcalW$.
\end{itemize} A $\k$-subquiver is \tit{wide} if $\ob\mcalW = \ob\mcalV$.
Let $(\mcalM,b^\mcalM)$ be an $\mcalA$-bimodule. An \tit{$\mcalA$-sub-bimodule} of $\mcalM$ is a wide $\k$-subquiver $\mcalN\inn \mcalM$ such that $b_{i|j}^\mcalM$ maps $\mcalA[1]^{\otimes i}\otimes \mcalN[1]\otimes \mcalA[1]^{\otimes j} $ into $\mcalN[1]$, for every $i,j\geq 0$.
\end{defn}

\begin{defn}
\label{wide ideal}
Let $\mcalA$ be an $A_\infty$ category. A \tit{wide ideal} of $\mcalA$ is an $\mcalA$-sub-bimodule of the diagonal bimodule $\mcalA$. A \tit{wide quotient} of $\mcalA$ is a strict $A_\infty$ functor $f:\mcalA\to\mcalA'$ such that $\ob\mcalA' = \ob\mcalA$, $f=1_{\ob\mcalA}$ on objects, and $f_1^{\ep_1,\ep_2}:\mcalA[1](\ep_2,\ep_1)\to \mcalA'[1](\ep_2,\ep_1)$ is surjective for each $\ep_1,\ep_2\in\ob\mcalA$.
\end{defn}

\begin{prop}
\label{wide ideal correspondence}
There is a one-to-one correspondence between the wide ideals and the wide quotients of $\mcalA$.
\end{prop}
\begin{proof}
If $\mcalI\inn\mcalA$ is a wide ideal, then the quotient $\k$-quiver $\mcalA/\mcalI$ is canonically an $A_\infty$ category, as follows. Since tensoring is right exact, it is not hard to see that the kernel of the map $\mcalA[1]^{\otimes k}\to (\mcalA/\mcalI)[1]^{\otimes k}$ is $\sum_{i+j = k-1} \mcalA[1]^{\otimes i}\otimes \mcalI[1]\otimes \mcalA[1]^{\otimes j}$, which is mapped into $\mcalI[1]$ by $b_k$. Hence, $b_k$ descends to a map $b^{\mcalA/\mcalI}_k: (\mcalA/\mcalI)[1]^{\otimes k}\to (\mcalA/\mcalI)[1]$. 
	Since $\{b_k\}$ satisfies the $A_\infty$ relations \eqref{A_infty relations}, so do $\{b^{\mcalA/\mcalI}_k\}$, making $(\mcalA/\mcalI,b^{\mcalA/\mcalI})$ into an $A_\infty$ category. By definition, the quotient map $\mcalA\to\mcalA/\mcalI$ of underlying $\k$-quivers satisfies the relation \eqref{strict A infty functor relation}, and so it is a wide quotient of $\mcalA$. Conversely, one can check that if $f:\mcalA\to \mcalA'$ is a wide quotient, then the kernel of its un-suspended component $\varphi_1:\mcalA\to \mcalA'$ is a wide ideal of $\mcalA$.
\end{proof}

\begin{defn}
    A \tit{degree $k$ pre-morphism} of $\mcalA$-bimodules $f:\mcalM\to \mcalN$ is a degree $k$ map of $B\mcalA$-bicomodules $f:B^\mcalM\mcalA\to B^\mcalN\mcalA$. Equivalently, it consists of degree $k$ maps  \[f_{i|j}: \mcalA[1]^{\otimes i}\otimes \mcalM[1]\otimes \mcalA[1]^{\otimes j}\to \mcalN[1]\]of $\k$-quivers on $\ob \mcalA$, for $i,j\geq 0$.  The \tit{un-suspended components} \[\varphi_{i|j} \defeq s^{-1}\circ f_{i|j}\circ s^{\otimes (i+1+j)}: \mcalA^{\otimes i}\otimes \mcalM\otimes \mcalA^{\otimes j}\to \mcalN\] have degree $k-i-j$. We say a pre-morphism $f$ is \tit{strict} if $f_{i|j} = 0$ for $i + j > 0$. 
\end{defn}    

The class of $\mcalA$-bimodules underlies a dg category, denoted by $\bimod{\mcalA}$. The hom-space $\bimod{\mcalA}(\mcalM,\mcalN)$ is the graded $\k$-module whose degree $k$ summand consists of degree $k$ pre-morphisms $\mcalM\to\mcalN$. The differential of a pre-morphism $f:\mcalM\to\mcalN$ is $ b^\mcalN f-(-1)^{|f|}fb^\mcalM$. Composition of pre-morphisms is just the composition of $B\mcalA$-bicomodule maps. 

\begin{defn}
    A \tit{morphism} of $\mcalA$-bimodules is a closed pre-morphism of degree $0$, or equivalently, a map $B^\mcalM\mcalA\to B^\mcalN\mcalA$ of dg $B\mcalA$-bicomodules. 
\end{defn}

The components again decompose into maps of graded $\k$-modules:
\[
    (f_{i|j})^{\ep_1,\cdots,\ep_{j+1}}_{\ep_1',\cdots,\ep_{i+1}'}:\mcalA[1](\ep_1',\cdots,\ep_{i+1}')\otimes \mcalM[1](\ep_{j+1},\ep_1')\otimes \mcalA[1](\ep_1,\cdots,\ep_{j+1})\to \mcalN[1](\ep_1,\ep_{i+1}')
    \]
    where $i,j\geq 0$ and $\ep_1,\cdots,\ep_{j+1},\ep_1',\cdots,\ep_{i+1}'\in \ob\mcalA$.

\begin{defn}
    A morphism $f:\mcalM\to \mcalN$ of $\mcalA$-bimodules is a \tit{quasi-isomorphism} if its un-suspended component $\varphi_{0|0}:(\mcalM,\beta_{0|0}^\mcalM)\to(\mcalN,\beta_{0|0}^\mcalN)$ is a quasi-isomorphism of dg $\k$-quivers.
\end{defn}

The dg category $\bimod{\mcalA}$ is (strongly) pretriangulated, in the sense of \cite{BK90}\cite{Dri04}, with shifts and cones defined as follows. For a proof, see \cite{Sei08}.

\begin{defn}
    The \tit{$n$-shift} of an $\mcalA$-bimodule $\mcalM$ is the $\k$-quiver $\mcalM[n]$, equipped with the codifferential $b^{\mcalM[n]}: B^{\mcalM[n]}A\to B^{\mcalM[n]}\mcalA$ given by
    \[b^{\mcalM[n]}_{i|j}\defeq (-1)^ns^n\circ b^\mcalM_{i|j}\circ (1^{\otimes i}\otimes s^{-n}\otimes 1^{\otimes j}): \mcalA[1]^{\otimes i}\otimes \mcalM[n+1]\otimes \mcalA[1]^{\otimes j}\to \mcalM[n+1]\,.\] 
    Similarly, the \tit{$n$-shift} $f[n]$ of a degree $k$ pre-morphism $f: \mcalM\to \mcalN$ is given by
    \[f[n]_{i|j}\defeq (-1)^{kn}s^n\circ f_{i|j}\circ (1^{\otimes i}\otimes s^{-n}\otimes 1^{\otimes j}): \mcalA[1]^{\otimes i}\otimes \mcalM[n+1]\otimes \mcalA[1]^{\otimes j}\to \mcalN[n+1]\,.\]
\end{defn}
\begin{defn}
    The \tit{mapping cone} of a morphism $f: \mcalM\to \mcalN$ is the $\mcalA$-bimodule \[\cone(f)\defeq \mcalM[1]\oplus \mcalN\] whose codifferential $b^{\cone(f)}$ on $B^{\cone(f)}\mcalA = B^{\mcalM[1]}\mcalA\oplus B^\mcalN\mcalA$ is given by
    \[b^{\cone(f)}_{i|j}\defeq \begin{bmatrix}
        b^{\mcalM[1]}_{i|j} & \\
        f_{i|j}\circ 1^{\otimes i}\otimes s^{-1} \otimes 1^{\otimes j} & b^\mcalN_{i|j}
    \end{bmatrix}\,.\]
\end{defn}
Hence, the \tit{homotopy category} $H^0(\bimod{\mcalA})$ is triangulated. Let $D(\bimod{\mcalA})$ denote the \tit{derived category} of $\mcalA$-bimodules, obtained by localizing $H^0(\bimod{\mcalA})$ at quasi-isomorphisms. 
\vspace{-2ex}
\begin{rmk}
\label{triangulated rmk}
    It is not known to the author whether the homotopy category of bimodules over a $\Z/n$-graded $A_\infty$ category is triangulated, but cf. \cite{PX97}\cite{Kel05}\cite{Liu23}. 
\end{rmk}
\begin{defn}
    Suppose $f,g: \mcalM\to \mcalN$ are two $\mcalA$-bimodule morphisms. A \tit{homotopy} between $f$ and $g$ is a degree $-1$ pre-morphism $h: \mcalM\to\mcalN$ such that
    $b^\mcalN h + hb^\mcalM = f-g$.
\end{defn}

Let $V$ be a graded $\k$-module. Denote $V^\dagger = \Hom_\k(V,\k)$, graded by $(V^\dagger)^k = \Hom_\k(V^{-k},\k)$. In particular, if $V$ has basis $v_1,\cdots,v_n$, then the dual basis $v_1^\dagger,\cdots, v_n^\dagger$ has grading $|v_i^\dagger| = -|v_i|$. More generally, if $\mcalV$ is a $\k$-quiver on $Z$, define its \tit{dual} $\mcalV^\dagger$ to be the $\k$-quiver on $Z$ with
\[\mcalV^\dagger(\ep_1,\ep_2)\defeq \mcalV(\ep_2,\ep_1)^\dagger\,.\] Note the reversal of objects. The below definition generalizes one made in \cite{Tra08}, although our signs differ from Tradler's.

\begin{defn}
\label{linear dual bimodule def}
Suppose $(\mcalM,b^\mcalM)$ is an $\mcalA$-bimodule. Its \tit{linear dual $\mcalA$-bimodule} is the $\k$-quiver $\mcalM^\dagger$ with codifferential $b^{\mcalM^\dagger}$, whose un-suspended components are given by 
\begin{equation}
\label{linear dual}
\begin{aligned}
(\beta^{\mcalM^\dagger}_{i|j})^{
\ep_1,\cdots,\ep_{j+1}}_{\ep_1',\cdots,\ep_{i+1}'}(a_i',\cdots,a_1',y,a_j,\cdots,a_1)(x) =(-1)^{\bsq}y\big(
(\beta^{\mcalM}_{j|i})_{\ep_1,\cdots,\ep_{j+1}}^{\ep_1',\cdots,\ep_{i+1}'}(a_j,\cdots,a_1,x,a_i',\cdots,a_1')\big)
\end{aligned}
\end{equation}
for any $\ep_1,\cdots,\ep_{j+1},\ep_1',\cdots,\ep_{i+1}'\in \ob\mcalA$, $a_j\otimes\cdots\otimes a_1\in\mcalA(\ep_1,\cdots,\ep_{j+1})$, $a_i'\otimes\cdots\otimes a_1'\in \mcalA(\ep_1',\cdots,\ep_{i+1}')$, $x\in \mcalM(\ep_{i+1}',\ep_1)$, and $y\in \mcalM^\dagger(\ep_{j+1},\ep_1')$, where $\bsq = (|a_1'|+\cdots +|a_i'|)(|y|+|a_1|+\cdots+|a_j|+|x|) + |y|(1+i+j)+(i+1)(j+1)$.

If $f:\mcalM\to \mcalN$ is a morphism of $\mcalA$-bimodules, then its \tit{linear dual morphism} $f^\dagger:\mcalN^\dagger\to\mcalM^\dagger$ has un-suspended components
\begin{equation}
\label{linear dual morphism}
\varphi^{\dagger}_{i|j}(a_i',\cdots,a_1',y,a_j,\cdots,a_1)(x) = (-1)^{\bsq'}y\big(\varphi_{j|i}(a_j,\cdots,a_1,x,a_i',\cdots,a_1')\big)
\end{equation}
for any $x\in \mcalM$ and $y\in \mcalN^\dagger$, where $\bsq' = (|a_1'|+\cdots +|a_i'|)(|y|+|a_1|+\cdots+|a_j|+|x|) + |y|(i+j)+(i+1)(j+1)$.
\end{defn}


We will also need to pullback bimodules along $A_\infty$ functors. To do that, we first recall the definition of cotensor products of bicomodules \cite{MM65}:

\begin{defn}
    Let $\mcalB$ be a $\k$-cocategory and $\mcalV$, $\mcalW$ be $\mcalB$-bicomodules. The \tit{cotensor product bicomodule} $\mcalV\bx_\mcalB\mcalW$ is the kernel of 
    $\DE_R^\mcalV\otimes 1 - 1\otimes \DE_L^\mcalW: \mcalV\otimes \mcalW\to \mcalV\otimes\mcalB\otimes\mcalW$, equipped with coactions $\DE_L = \DE_L^\mcalV$ and $\DE_R = \DE_R^\mcalW$.

    Let $f:\mcalB\to\mcalB'$ be a map of $\k$-cocategories and $\mcalV$ a $\mcalB'$-bicomodule. Then $\mcalB$ is naturally a $(\mcalB,\mcalB')$-bicomodule and a $(\mcalB',\mcalB)$-bicomodule. The \tit{coinduction} of $\mcalV$ along $f$ is the $\mcalB$-bicomodule $f^*\mcalV\defeq \mcalB\bx_{\mcalB'}\mcalV\bx_{\mcalB'}\mcalB$.   
\end{defn}

\begin{defn}
\label{pullback def}
    Suppose $f: \mcalA\to \mcalA'$ is an $A_\infty$ functor, and $(\mcalM,b^\mcalM)$ is an $\mcalA'$-bimodule. Then the coinduction
    \[B\mcalA\bx_{B\mcalA'}B^\mcalM\mcalA'\bx_{B\mcalA'}B\mcalA = B^\mcalM\mcalA\]
    is a dg $B\mcalA$-bicomodule, whose codifferential we denote by $b^{f^*\mcalM}$.
    Define the \tit{pullback} $f^*\mcalM$ to be the $\mcalA$-bimodule $(\mcalM,b^{f^*\mcalM})$.
\end{defn}
If $f$ is strict, then the pullback codifferential is easy to describe: its un-suspended components are given by
\[\beta^{f^*\mcalM}_{i|j}(a_i',\cdots ,a_1',x,a_j,\cdots ,a_1) = \beta^{\mcalM}_{i|j}(\varphi_{0|0}(a_i'),\cdots ,\varphi_{0|0}(a_1'),x,\varphi_{0|0}(a_j),\cdots ,\varphi_{0|0}(a_1))\,.\]

\section{Composable LSFT Algebra}
\label{composable LSFT algebra}
While the Chekanov-Eliashberg algebra counts holomorphic disks with a single positive puncture, rational symplectic field theory aims to encompass the information of all disks. However, it turns out disks with multiple positive punctures admit boundary bubbling phenomena that do not occur for disks with a single positive puncture. To account for these boundary bubblings, Cieliebak and Latschev put forth a program \cite{CL09} employing operations from string topology of Chas and Sullivan. In \cite{Ng10}, Ng implemented this program in the case of Legendrian knots in $\R^3$, producing a filtered curved dga called the LSFT algebra, whose zeroth associated graded is the Chekanov-Eliashberg algebra. Ng later dealt with multi-component Legendrian links in \cite{Ng23}, but in this case the LSFT differential no longer squares to something simple.

In this section, we define a variant called the \tit{composable LSFT algebra}, which we show is a filtered curved dga for any Legendrian links, with curvature given by a formula generalizing \cite[Proposition 3.15]{Ng10}. Here, ``composable'' means $\k^m$-linear, where $m$ is the number of base points (see \Cref{composability rmk}). Most definitions and lemmas carry over to the composable version, with minimal modifications.

\subsection{Generators and Gradings}
\label{gradings}
Let $\Lambda$ be an oriented Legendrian link in $\R^3$. Place \tit{marked points} $\star_1,\cdots,\star_m$ on $\LA$ such that there is at least one on each component. With respect to the orientation of $\Lambda$, place a \tit{base point} $\bullet_i$ immediately after each marked point $\star_i$. In \cite{Ng23}, the above data is called a \tit{pointed Legendrian link}. 

Suppose $\LA$ has generic Lagrangian projection $\LA_{xy}$. Label its Reeb chords by $a_1,\cdots,a_n$. To each Reeb chord $a_j$, associate two variables $p_j,q_j$. To each marked point $\star_i$, associate two variables $t_i,t_i^{-1}$. Denote $\mbfP = \{p_1,\cdots,p_n\}$, $\mbfR = \mbfP\cup \mbfQ = \{p_1,\cdots,p_n,q_1,\cdots, q_n\}$,
\[\mbfS = \mbfP\cup \mbfS^0 = \{p_1,\cdots,p_n,q_1,\cdots,q_n,t^{\pm 1}_1,\cdots,t^{\pm 1}_m\}\,.\]
Define an involution on $\mbfR$ by
\[\bar{s} = \begin{cases}
    q_j & \text{if } s = p_j\\
    p_j & \text{if } s = q_j
\end{cases}\,.\]
To define gradings and link gradings of the generators, we will associate an oriented broken path $\gamma_s$ lying on $\LA$, called the \tit{capping path}, to each $s\in\mbfS$. (Each $\gamma_s$ is a composable based string, as in \Cref{composable strings def}.) The starting ($-$) and ending ($+$) points of $\gamma_s$ are determined by the following: 
\begin{defn}
\label{internal link grading def}
The \tit{internal link grading} of the LSFT generators of $\Lambda$ is the pair of functions $e_-,e_+:\mbfS\to \{1,\cdots, m\}$, defined as follows. Since each component of $\Lambda$ contains at least one marked point $\star_i$, removing $\star_1,\cdots,\star_m$ leaves a disjoint union of oriented open intervals:
\[\Lambda\setminus\{\star_1,\cdots,\star_m\} = I_1\sqcup \cdots\sqcup I_m\]
where $\bullet_i\in I_i$. Then the values of $e_-,e_+$ on $p_j,q_j,t_i^{\pm 1}$ are uniquely defined by the following:
\begin{enumerate}[label = (\arabic*)]
    \item $a_j^-$ lies in $I_{e_-(p_j)} = I_{e_+(q_j)}$;
    \item $a_j^+$ lies in $I_{e_+(p_j)} = I_{e_-(q_j)}$;
    \item $e_+(t_i) = e_-(t_i^{-1}) = i$;
    \item $\bullet_{e_-(t_i)} = \bullet_{e_+(t_i^{-1})}$ immediately precedes $\bullet_i$ according to the orientation of $\LA$.
\end{enumerate}
 
\end{defn}

Now we are ready to define the capping paths. For each $1\leq i\leq m$, let $\gamma_{t_i}$ be the oriented segment of $\LA$ starting at $\bullet_{e_-(t_i)}$ and ending at $\bullet_{e_+(t_i)}$; let $\gamma_{t_i^{-1}}\defeq -\gamma_{t_i}$ be the reverse. For each $1\leq j\leq n$, let $\gamma_j^\pm$ be the oriented segment of $\LA$ starting at $\bullet_{e_\pm(p_j)} = \bullet_{e_\mp(q_j)}$ and ending at $a_j^\pm$. Then set
\begin{align*}
    \gamma_{p_j}&\defeq \gamma_j^-\cdot (-\gamma_j^+)\\
    \gamma_{q_j}&\defeq \gamma_j^+\cdot (-\gamma_j^-)
\end{align*}
where $\gamma_1\cdot \gamma_2$ denotes the concatenated path that traverses $\gamma_1$ and then $\gamma_2$.

To assign gradings to the generators, we need to choose an auxillary piece of data:

\begin{defn}[{\cite[Definition 2.3]{Ng23}}]
\label{Maslov potential def}
    A \tit{Maslov potential} on $\LA$ is a map $\mu:\{\bullet_1,\cdots,\bullet_m\}\to \R$ such that $\mu(\bullet_i)\in\R$ is a lift of the unit tangent vector ($\in S^1 = \R/\Z$) of $\LA_{xy}$ at $\bullet_i$.
\end{defn}

Given a Maslov potential $\mu$ on $\LA$, define the grading
\begin{equation}
    \label{grading def}
    |s|\defeq \flr{2(-\rot( \gamma_s) + \mu(\bullet_{e_+(s)}) - \mu(\bullet_{e_-(s)}))}\in\Z
\end{equation}
for each $s\in \mbfS$, where $\rot(\gamma)\in\R$ for any oriented immersed curve $\gamma\inn\R^2$ is the number of counterclockwise revolution in $S^1$ made by the unit tangent vector $\gamma'/|\gamma'|$, and $\rot(\gamma)\defeq \rot(\pi_{xy}\circ \gamma)$ for $\gamma\inn\R^3$. If $\gamma\inn\R^2$ has unit speed, then $\rot(\gamma) = \frac{1}{2\pi}\int\gamma''$.
\begin{defn}
\label{rot num def}
    If $\LA\inn\R^3$ is an oriented Legendrian knot, then its \tit{rotation number} $\rot(\LA)$ is defined as above. If $\LA$ is an oriented Legendrian link with components $\LA_1,\cdots,\LA_\ell$, then
    \[\rot(\LA)\defeq \begin{cases}
        0 & \rot(\LA_1)=\cdots = \rot(\LA_\ell)=0\\
        \gcd\{\rot(\LA_1),\cdots, \rot(\LA_\ell)\} & \text{otherwise}
    \end{cases}\,.\]
\end{defn}
\begin{prop}[{\cite[Proposition 2.7]{Ng23}}]
\label{grading prop}
For any Maslov potential $\mu$, we have 
    \begin{enumerate}[label = (\arabic*)]
        \item $|t_i| = -|t_i^{-1}|\in 2\Z$ for any marked point $\star_i$;
        \item If $\LA_0$ is a component of $\LA$, then $\sum_{\bullet_i\in \LA_0} |t_i| = -2\rot(\LA_0)$;
        \item $|p_j| + |q_j| = 1$ for any Reeb chord $a_j$.
    \end{enumerate}
Moreover, the set of all possible gradings is a torsor for $\Z^m/(1,\cdots,1)$, which acts by changing $\mu$. If we increase $\mu$ by $(k_1,\cdots,k_m)\in\Z^m$, then $|s|$ increases by $2(k_{e_+(s)}-k_{e_-(s)})$, for any $s\in \mbfS$. In particular, a grading is determined by its values on $t_1,\cdots, t_m$.
\end{prop}
\begin{lemma}
    \label{rotation number base point grading} For any pointed Legendrian link $\LA$, there is a Maslov potential such that $|t_i| \equiv 0\mod 2\rot(\LA)$ for each $i$.
\end{lemma}
\begin{proof}
     Since the grading of any particular $t_i$ is determined by the Maslov potential restricted to the component containing $\bullet_i$, we can assume $\LA$ is connected. Label the base points $\bullet_1,\cdots,\bullet_m$ so that $\bullet_i$ appears before $\bullet_{i+1}$ according to the orientation of $\LA$. Pick any lift $\mu(\bullet_m)\in\R$ of the unit tangent vector of $\LA$ at $\bullet_m$. For $i < m$, set
    \[\mu(\bullet_i) = \mu(\bullet_m) + \rot(\gamma_{t_1})+\cdots + \rot(\gamma_{t_i})\] which defines a valid Maslov potential. Then for each $i < m$, we have $\mu(\bullet_{e_+(t_i)}) - \mu(\bullet_{e_-(t_i)}) = \mu(\bullet_{i}) - \mu(\bullet_{i-1}) = \rot(\gamma_{t_i})$, and so $|t_i| = 0$. For $t_m$, we have \[|t_m|  
 = \flr{2(-\rot( \gamma_s) + \mu(\bullet_m) - \mu(\bullet_{m-1}))} = -2(\rot(\gamma_{t_1})+\cdots+\rot(\gamma_{t_m})) = -2\rot(\LA)\] which concludes the proof. 
\end{proof}


Let $\k^m = \k e_1\oplus \cdots \oplus \k e_m$ be the semisimple ring with idempotents $e_1,\cdots,e_m$. By a $\k^m$-algebra, we mean a monoid in the category of $\k^m$-bimodules. 
\begin{defn}
\label{km bimodule}
    Suppose $\mbf{S}$ is a finite set equipped with functions $e_-,e_+:\mbf{S}\to\{1,\cdots,m\}$ and grading $|\cdot|:\mbfS\to\Z/n$ for some $n\in\Z$. Define a $\Z/n$-graded $\k^m$-bimodule $V(\mbfS;e_-,e_+)$ as follows: the underlying graded $\k$-module is free with basis $\mbfS$; the idempotents act by \[e_ise_j = \de_{i,e_-(s)}\de_{j,e_+(s)}s\] for any $s\in\mbfS$.
\end{defn}
\begin{rmk}
\label{composability rmk}
    The tensor algebra over $\k^m$ of $V(\mbfS;e_-,e_+)$ is exactly the path algebra of the graded quiver with vertex set $\{1,\cdots,m\}$, edge set $\mbfS$, source function $e_-$, and target function $e_+$. The adjective ``composable'' refers to composable arrows/paths in this quiver.
\end{rmk}
Fix a Maslov potential of $\LA$, and equip $\mbfS$ with the grading \eqref{grading def}. Let $\mcalL_0=\mcalL_0(\LA)$ be the tensor algebra of the $\Z$-graded $\k^m$-bimodule $V(\mbfS;e_-,e_+)$, modulo the relations $t_it_i^{-1} = t_i^{-1}t_i = e_i$. Define a descending filtration \[\mcalL_0 = \mcalF^0\mcalL_0 \cont \mcalF^1\mcalL_0\cont\cdots\] by letting $\mcalF^1\mcalL_0$ be the two-sided ideal generated by $p_1,\cdots,p_n$, and $\mcalF^k\mcalL_0\defeq (\mcalF^1\mcalL_0)^k$.

\begin{defn}
    The underlying filtered graded $\k^m$-algebra of the \tit{composable LSFT algebra} of $\LA$ is the completion \[\mcalL\defeq \widehat{\mcalL_0}\] with respect to the above filtration. Define the graded $\k$-module $\mcalL^\mrm{cyc}$ as the quotient of $\mcalL$ by the $\k$-submodule generated by the commutators $[x,y] = xy - (-1)^{|x||y|}yx$, where $x,y\in \mcalL$. For any $x\in\mcalL$, write $[x]\in \mcalL^\mrm{cyc}$ for its image under the quotient map.
\end{defn}

In the remainder of this section, we will follow \cite{Ng10}\cite{Ng23} and define the following structures on $\mcalL$ and $\mcalL^\mrm{cyc}$:
\begin{itemize}
    \item the \tit{SFT bracket} $\{\cdot,\cdot\}:\mcalL^\mrm{cyc}\otimes \mcalL\to \mcalL$;
    \item the \tit{string differential} $\de:\mcalL\to\mcalL$;
    \item the \tit{Hamiltonian} $h\in\mcalL^\mrm{cyc}$;
    \item the \tit{LSFT differential} $d: \mcalL\to\mcalL$. 
\end{itemize}
To define these operations, we need the geometric interpretation of (cyclically) composable words as broken strings on $\LA$, which is the subject of the following subsection.


\subsection{Composable Words and Strings}

\begin{figure}[t]
    \centering
    \begin{subfigure}{.5\textwidth}
        \centering
        \trimbox{30pt 140pt 5pt 130pt}{
\begin{tikzpicture} 
   \begin{knot}[
     consider self intersections = true,
     end tolerance = 1pt,
     clip width = 5pt,
     draft mode = off
   ]
   \strand[thick, rounded corners] (1,.5) .. controls (-6,-7) and (-6,7) .. (1,-.5);
   \strand[thick, looseness = 4](1,.5) to [out = 45, in = -45] (1,-.5);
   \strand[thick, rounded corners, yshift = 2cm] (1,.5) .. controls (-6,-7) and (-6,7) .. (1,-.5);
   \strand[thick, looseness = 4, yshift = 2cm](1,.5) to [out = 45, in = -45] (1,-.5);
   \flipcrossings{1,3,4}
   \end{knot} 
	\draw[ultra thick, color = teal] (-2.2,3.6) .. controls (-1.2,3.3) and (-.5,2.7) ..  (0.3,2) .. controls (0,1.7) and (-.58,1.22) .. (-.67,1.16) .. controls (-2,2.14) and (-3.4,2.15) .. (-3.98,1.3) .. controls (-4.15,1.6) and (-4.15,2.4) .. (-3.92,2.8);
	\path[tips, -Stealth, thick, shift = {(-2.75,.21)}] (1,0) -- (0,0);
	\path[tips, -Stealth, thick, shift = {(-2.73,-1.78)}] (-1,0) -- (0,0);
	\path[tips, -Stealth, thick, shift = {(-2.5,1.92)}, color  = teal] (-1,0) -- (0,0);
   \node at (-2.5,3.78) {\tiny$\bigstar$};
   \node[above = 3pt] at (-2.6,3.78) {\scriptsize$\bigstar_2$};
   \node at (-4.17,2.65) {\tiny$\bigstar$};
   \node[left = 3pt] at (-4.17,2.6) {\scriptsize$\bigstar_1$};
   \node at (-1,-1.22) {\tiny$\bigstar$};
   \node[below = 3pt] at (-1,-1.25) {\scriptsize$\bigstar_3$};
   \node[] at (-2.2,3.73) {\footnotesize$\bullet$}; 
   \node[above = 4pt] at (-2,3.73) {\footnotesize$\bullet_2$};
   \node[] at (-4.05,2.9) {\footnotesize$\bullet$};
   \node[left = 3pt] at (-4.05,3.1) {\footnotesize$\bullet_1$};
   \node[] at (-.75,-1.05) {\footnotesize$\bullet$};  
   \node[right = 3pt] at (-.75,-1.2) {\small$\bullet_3$};
   \node[below = 5pt] at (0.5,2) {\footnotesize$q_3$};
   \node[right = 5pt] at (0.5,2) {\footnotesize$p_3$};
   \node[left = 9pt,teal] at (0.5,2) {\footnotesize$\tau_3$};
   \node[below = 5pt] at (-.65,1) {\footnotesize$q_2$};
   \node[right = 5pt] at (-.65,1) {\footnotesize$p_2$};
   \node[above = 8pt,teal] at (-.65,1) {\footnotesize$\tau_2$};
   \node[below = 7pt] at (-3.95,1) {\footnotesize$q_1$};
   \node[right = 4pt] at (-4,1) {\footnotesize$p_1$};
   \node[above = 16pt,teal] at (-3.85,1) {\footnotesize$\tau_1$};
   \node[below = 5pt] at (.5,0) {\footnotesize$q_4$};
   \node[right = 5pt] at (.5,0) {\footnotesize$p_4$};
\end{tikzpicture}
}
        \caption{composable based string}
        \label{composable based string}
    \end{subfigure}%
    \begin{subfigure}{.5\textwidth}
        \centering
        \trimbox{30pt 140pt 5pt 130pt}{
\begin{tikzpicture} 
   \begin{knot}[
     consider self intersections = true,
     end tolerance = 1pt,
     clip width = 5pt,
     draft mode = off
   ]
   \strand[thick, rounded corners] (1,.5) .. controls (-6,-7) and (-6,7) .. (1,-.5);
   \strand[thick, looseness = 4](1,.5) to [out = 45, in = -45] (1,-.5);
   \strand[thick, rounded corners, yshift = 2cm] (1,.5) .. controls (-6,-7) and (-6,7) .. (1,-.5);
   \strand[thick, looseness = 4, yshift = 2cm](1,.5) to [out = 45, in = -45] (1,-.5);
   \flipcrossings{1,3,4}
   \end{knot} 
	\draw[ultra thick, color = teal] (-.9,1) .. controls (-2.3,1.95) and (-3.4,1.8) .. (-3.83,1) .. controls (-3.4,.2) and (-2.3,0.05)  .. cycle;
	\path[tips, -Stealth, thick, shift = {(-2.75,.21)}] (1,0) -- (0,0);
	\path[tips, -Stealth, thick, shift = {(-2.73,-1.78)}] (-1,0) -- (0,0);
	\path[tips, -Stealth, thick, shift = {(-2.5,1.66)}, color  = teal] (-1,0) -- (0,0);
   \node at (-2.5,3.78) {\tiny$\bigstar$};
   \node[above = 3pt] at (-2.6,3.78) {\scriptsize$\bigstar_2$};
   \node at (-4.17,2.65) {\tiny$\bigstar$};
   \node[left = 3pt] at (-4.17,2.6) {\scriptsize$\bigstar_1$};
   \node at (-1,-1.22) {\tiny$\bigstar$};
   \node[below = 3pt] at (-1,-1.25) {\scriptsize$\bigstar_3$};
   \node[] at (-2.2,3.73) {\footnotesize$\bullet$}; 
   \node[] at (-2.2,3.73) {\footnotesize$\bullet$}; 
   \node[above = 4pt] at (-2,3.73) {\small$\bullet_2$};
   \node[] at (-4.05,2.9) {\footnotesize$\bullet$};
   \node[left = 3pt] at (-4.05,3.1) {\small$\bullet_1$};
   \node[] at (-.75,-1.05) {\footnotesize$\bullet$}; 
   \node[right = 3pt] at (-.75,-1.2) {\small$\bullet_3$};
   \node[below = 5pt] at (0.5,2) {\footnotesize$q_3$};
   \node[right = 5pt] at (0.5,2) {\footnotesize$p_3$};
   \node[below = 5pt] at (-.65,1) {\footnotesize$q_2$};
   \node[right = 5pt] at (-.65,1) {\footnotesize$p_2$};
   \node[left = 10.5pt,teal] at (-.65,1) {\footnotesize$\tau_0$};
   \node[below = 7pt] at (-3.95,1) {\footnotesize$q_1$};
   \node[left = 2pt] at (-4,1) {\footnotesize$p_1$};
   \node[right = 7pt,teal] at (-4,1) {\footnotesize$\tau_1$};
   \node[below = 5pt] at (.5,0) {\footnotesize$q_4$};
   \node[right = 5pt] at (.5,0) {\footnotesize$p_4$}; 
\end{tikzpicture}
}
        \caption{cyclically composable string}
        \label{cyclically composable string}
    \end{subfigure}
    \caption{Illustrated here are the projections of two generic composable strings with holomorphic corners, which are offset from the link for readability. The labelings at crossings are explained in \Cref{holomorphic corner rmk}.\\ \hspace{\textwidth}\hspace*{1em} On the left, we have a composable based string of length $3$, corresponding to the composable word $t_1^{-1}q_1q_2p_3$. On the right, we have a cyclically composable string of length $2$, corresponding to the cyclically composable word $p_2p_1$. Unlike the broken closed strings considered in \cite{Ng23}, our composable strings are allowed to jump across Reeb chords but not between base points.}
    \label{strings picture}
\end{figure}
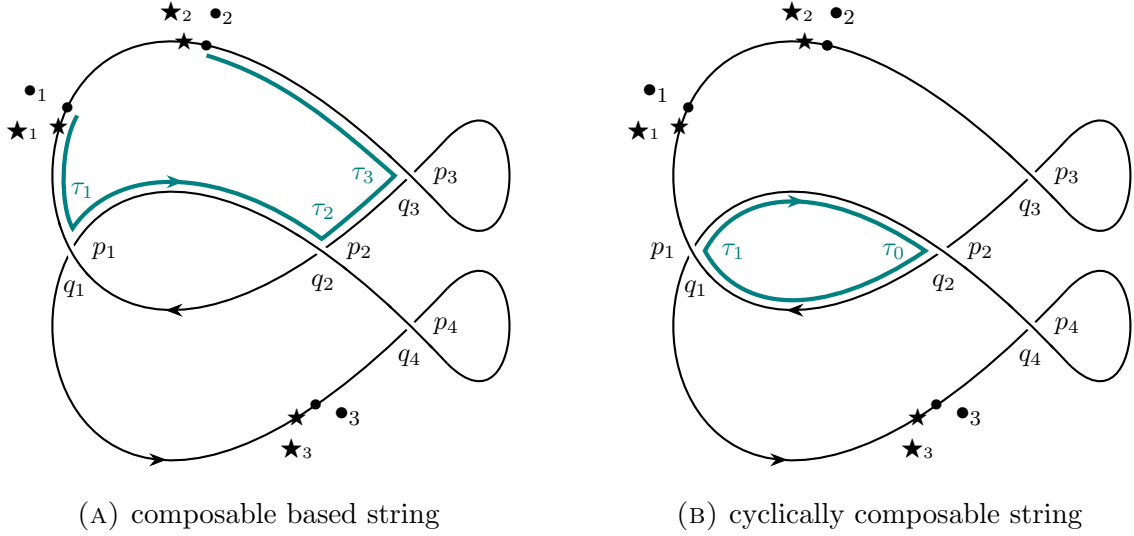

Define the set of \tit{composable words} (reduced according to the relations $t_it_i^{-1} = t_i^{-1}t_i = e_i$)
\[\mcalW = \{s_1\cdots s_k\mid s_i\in\mbfS, e_+(s_i) = e_-(s_{i+1}) \text{ for }1\leq i\leq k-1\}\cup \{e_1,\cdots,e_m\}\,.\]
Also define the set of \tit{cyclically composable words}
\[\mcalW_\circ = \{s_1\cdots s_k\in \mcalW\mid e_+(s_k) = e_-(s_1)\}\,.\]
Extend the internal link grading on generators to words by setting
\begin{enumerate}[label = (\arabic*)]
    \item $e_-(s_1\cdots s_k) \defeq e_-(s_1)$;
    \item $e_+(s_1\cdots s_k) \defeq e_+(s_k)$;
    \item $e_-(e_i) = e_+(e_i) \defeq i$,
\end{enumerate}
for any $s_1,\cdots,s_k\in\mbfS$.
\begin{rmk}
Note that $\mcalW$ is a $\k$-linear basis of $\mcalL$, but to define a $\k^m$-linear map out of $\mcalL$, one needs to make sure it preserves $e_-,e_+$ on $\mcalW$.
\end{rmk}
The following definition is similar to and generalizes \cite[Definition 3.1]{Ng10}. Also compare with \cite[Definition 2.9]{Ng23}. See \Cref{strings picture} for illustrations.

\begin{defn}
\label{composable strings def}
    For any $k\geq 0$, pick cyclically ordered points $\tau_0,\tau_1,\cdots,\tau_k$ on an oriented circle $S^1$. A \tit{composable based string} of length $k$ is a smooth map
    \[\gamma: S^1\setminus\{\tau_0,\tau_1,\cdots,\tau_k\}\to \Lambda\]
    such that
    \begin{enumerate}[label = (\arabic*)]
        \item for each $1\leq i\leq k$, there is some Reeb chord $a_j$ such that
        $e_\pm^i(\gamma)\defeq \lim_{\tau\to \tau_i^\pm} \gamma(\tau)= a_j^\pm\text{ or } a_j^\mp$ (we say $\gamma$ has a \tit{jump} at the Reeb chord $a_j$);
        \item $e_\pm(\gamma)\defeq \lim_{\tau\to \tau_0^\pm} \gamma(\tau)\in\{\bullet_1,\cdots,\bullet_m\}$.
    \end{enumerate}
    Let $\gamma_{xy}: S^1\setminus\{\tau_0\}\to \LA_{xy}$ denote the unique continuous extension of $\pi_{xy}\circ\gamma$. 
    
    Similarly, a \tit{cyclically composable string} of length $k$ is a smooth map
\[\gamma: S^1\setminus\{\tau_1,\cdots,\tau_k\}\to \Lambda\] satisfying condition (1) for each $1\leq i\leq k$.  Let $\gamma_{xy}: S^1\to \Lambda$ denote the continuous extension of $\pi_{xy}\circ\gamma$.

All strings are considered up to orientation-preserving reparametrizations of the domain.

A \tit{homotopy} of cyclically composable strings (resp. composable based strings) is a homotopy of smooth maps that preserves each $e_\pm^i$ (resp. $e_\pm^i$'s and $e_\pm$).
\end{defn}
Define the maps
\[\begin{cases}w:\{\text{composable based strings}\}\to \mcalW\\ w:\{\text{cyclically composable strings}\}\to \mcalW_\circ\end{cases}\]
as follows. Suppose $\gamma$ is a composable based string of length $k$. For each $1\leq i\leq k$, let
\begin{equation}
    \label{wi of gamma}
    w_i(\gamma)\defeq \begin{cases}
    p_j & e_\pm^i(\gamma) = a_j^\pm\\
    q_j & e_\pm^i(\gamma) = a_j^\mp
\end{cases}\,.
\end{equation}
For each $0\leq i\leq k$, let 
\begin{equation}
\label{w i i+1 of gamma}
w_{i,i+1}(\gamma)\defeq t_{i_1}^{n_1}\cdots t_{i_r}^{n_r}\end{equation}
where $n_\ell\in\Z$ is the number of times the oriented interval $\gamma(\tau_i,\tau_{i+1})$ crosses $\star_{i_\ell}$, counted according to the orientation of $\LA$. Now define
\[w(\gamma) \defeq  w_{01}(\gamma)w_1(\gamma)w_{12}(\gamma)\cdots w_k(\gamma)w_{k0}(\gamma)\,.\]
Similarly, if $\gamma$ is a cyclically composable string of length $k$, then define
\[w(\gamma) \defeq w_1(\gamma)w_{12}(\gamma)\cdots w_k(\gamma)w_{k1}(\gamma)\,.\]
Homotopies preserve $w$, since they preserve jumps at Reeb chords (and the base points). Moreover, we have

\begin{figure}[H]
    \centering
     \begin{tikzpicture}
      \draw[ultra thick] (1,-1) -- (-1,1);
      \draw[ultra thick] (1,1) -- (.2,.2);
      \draw[ultra thick] (-1,-1) -- (-.2,-.2);
      \node (p) at (.5,0) {$p_j$};
      \node (p) at (-.5,0) {$p_j$};
      \node (q) at (0,.5) {$q_j$};
      \node (q) at (0,-.5) {$q_j$};
  \end{tikzpicture}
    \caption{Labeling the quadrants at the crossing $a_j$ by the associated LSFT generators $p_j,q_j$. Note that opposite quadrants share the same label.}
    \label{corners}
\end{figure}
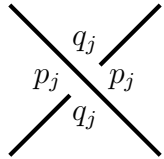

\begin{prop}
The maps
\[\begin{cases}w:\{\text{composable based strings}\}/\text{homotopy}\to \mcalW\\ w:\{\text{cyclically composable strings}\}/\text{homotopy}\to \mcalW_\circ\end{cases}
\]
are bijections.
\end{prop}
\begin{proof}
    Given any $s_1\cdots s_k\in\mcalW$, we have
    $w(\gamma_{s_1}\cdots \gamma_{s_k}) = s_1\cdots s_k$, so the first $w$ is surjective. If $s_1\cdots s_k\in\mcalW_\circ$, then $\gamma_{s_1}\cdots \gamma_{s_k}$ is cyclically composable, so the second $w$ is surjective.

    Suppose $\gamma_1,\gamma_2$ are composable based strings of length $k$ with $w(\gamma_1) = w(\gamma_2)$, then $w_i(\gamma_1) = w_i(\gamma_2)$ for each $1\leq i\leq k$ and $w_{i,i+1}(\gamma_1) = w_{i,i+1}(\gamma_2)$ for each $0\leq i\leq k$. Since $w_i(\gamma_1) = w_i(\gamma_2)$, we get $e^i_\pm(\gamma_1) = e^i_\pm(\gamma_2)$. Moreover, we have $e_-(\gamma_1) = e_-(\gamma_2)$, both being $e_-$ of the first generator in the word $w(\gamma_1) = w(\gamma_2)$, and similarly for $e_+$. Hence, after a homotopy, we can assume $\gamma_1,\gamma_2$ agree on some neighborhood of each puncture. It remains to relate $\gamma_1, \gamma_2$ via a compactly supported homotopy on each $(\tau_i,\tau_{i+1})$. Since $w_{i,i+1}(\gamma_1) = w_{i,i+1}(\gamma_2)$, we can homotope to get $(\gamma_1|_{(\tau_i,\tau_{i+1})})^{-1}\{\star_1,\cdots,\star_m\} = (\gamma_2|_{(\tau_i,\tau_{i+1})})^{-1}\{\star_1,\cdots,\star_m\}$ and make $\gamma_1,\gamma_2$ agree on a neighborhood of this preimage. Removing this preimage, we have divided $(\tau_i,\tau_{i+1})$ into subintervals, each of which is mapped by $\gamma_1$ and $\gamma_2$ into the same component of $\LA\setminus\{\star_1,\cdots,\star_m\}$. But each component of $\LA\setminus\{\star_1,\cdots,\star_m\}$ is contractible, so we are done. The argument for cyclically composable strings is similar.
\end{proof}
Ng introduced the following conditions in \cite{Ng10}, which allow us to work with the projections of strings instead of the strings themselves.
\begin{defn}
    A composable based string $\gamma: S^1\setminus\{\tau_0,\tau_1,\cdots,\tau_k\}\to\LA$  or a cyclically composable string $\gamma: S^1\setminus\{\tau_1,\cdots,\tau_k\}\to\LA$ is \tit{generic} if 
    \begin{enumerate}[label = (\arabic*)]
        \item $\gamma'(\tau_i^\pm)\defeq \lim_{\tau\to\tau_i^\pm}\gamma'(\tau)\neq 0$ for each $1\leq i\leq k$;
        \item $\gamma'(\tau)\neq 0$ whenever $\gamma(\tau)$ is a Reeb chord endpoint.
    \end{enumerate}
Suppose $\gamma$ is generic. We say $\gamma$ \tit{has holomorphic corners} if $(\gamma_i'(\tau_i^-),\gamma(\tau_i^+)-\gamma(\tau_i^-),\gamma'(\tau_i^+))$ is a positive frame in $\R^3$, for each $1\leq i\leq k$. Equivalently, $\gamma_{xy}$ makes a left turn at each $\tau_i$. 
\end{defn}
\begin{lemma}[{\cite[Section 3.2]{Ng10}}]
    Each composable based string/cyclically composable string is homotopic to one that is generic with holomorphic corners.
\end{lemma}
\begin{rmk}
\label{holomorphic corner rmk}
    The projection $\pi_{xy}$ gives a bijection between cyclically composable strings that are generic with holomorphic corners and oriented closed curves on $\LA_{xy}$ that have left-turning corners at crossings and non-zero velocities when they pass crossings (with or without turning). There is a similar bijection for composable based strings. Under these bijections, jumps at Reeb chords of a string $\gamma$ correspond to left-turning corners of $\gamma_{xy}$. In this case, the generators associated to jumps of $\gamma$ can be read off from the corners of $\gamma_{xy}$, by labeling the quadrants at a crossing as specified in \Cref{corners}.
\end{rmk}
\subsection{SFT Bracket}
\label{SFT Bracket}
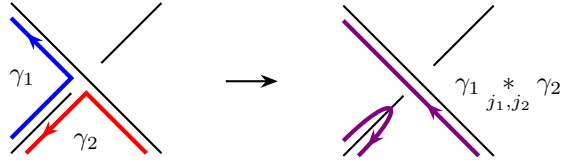
\begin{figure}[H]
\centering
  \begin{tikzpicture}
      \node[outer xsep = 20pt] (pic1) at (0,0) {\begin{tikzpicture}
    \draw[thick] (1,-1) -- (-1,1);
    \draw[thick] (1,1) -- (.2,.2);
    \draw[thick] (-1,-1) -- (-.2,-.2);
    \draw[color = red, ultra thick] (.8,-1) -- (0, -.2) -- (-.8,-1);
    \draw[color = blue, ultra thick] (-1,.8) -- (-.2,0) -- (-1,-.8);
    \node at (-.85, 0) {\footnotesize $\gamma_1$}; 
    \node at (0, -.85) {\footnotesize $\gamma_2$}; 
    \path[tips, -Stealth, color = blue, thick,shift = {(.2,-.2)}] (-.4,.2)--(-1,.8);
    \path[tips, -Stealth, color = red, thick,shift = {(.2,.2)}] (-.4, -.6)--(-.8,-1);
\end{tikzpicture}};
      \node[outer xsep = 20pt] (pic2) at (5,0) {\begin{tikzpicture}
    \draw[thick] (1,-1) -- (-1,1);
    \draw[thick] (1,1) -- (.2,.2);
    \draw[thick] (-1,-1) -- (-.2,-.2);
    \draw[color = violet, ultra thick] (-1,.8) -- (.8,-1);
    \draw[color = violet, ultra thick] (-.8,-1) .. controls (-.1,-.3) and (-.3,-.1) .. (-1,-.8);
    \path[tips, -Stealth, color = violet, thick, shift={(.1,-.1)}] (.3,-.5)--(0, -.2);
    \path[tips, -Stealth, color = violet, thick,shift = {(.1,.1)}] (-.5,-.65)--(-.78,-.98);
    \node at (1.2,-.2) {\footnotesize $\gamma_1\underset{j_1,j_2}{*}\gamma_2$}; 
\end{tikzpicture}};
      \draw[thick, -Stealth] (pic1) -- (pic2);
  \end{tikzpicture}
\caption{Local picture of the gluing operation for two generic composable strings with holomorphic corners}
\label{SFT bracket local}
\end{figure}

Suppose $\gamma_1$ is a cyclically composable string of length $k_1$, $\gamma_2$ is a composable based string of length $k_2$, and they each have a jump at the same Reeb chord but in opposite directions. Namely, there are $1\leq j_1\leq k_1$ and $1\leq j_2\leq k_2$ such that
\[e_\pm^{j_1}(\gamma_1) = e_\mp^{j_2}(\gamma_2)\,.\]
Then we can glue $\gamma_1$ and $\gamma_2$ at this Reeb chord to obtain a composable based string $\gamma_1\underset{j_1,j_2}{*}
\gamma_2$ of length $k_1+k_2-2$. See \Cref{SFT bracket local} for an illustration. 
After inserting signs and summing over all such pairs $(j_1,j_2)$, this gluing procedure defines the following operation: 
\begin{defn}
\label{SFT bracket}
    View $\mcalL^\mrm{cyc}\otimes_{\k} \mcalL$ as a $\k^m$-bimodule by letting $\k^m$ act on $\mcalL$. The \tit{SFT bracket}
    \[\{\cdot,\cdot\}:\mcalL^\mrm{cyc}\otimes_{\k} \mcalL\to \mcalL\]
    is a degree $1$ $\k^m$-linear map defined as follows. For any $s,s'\in \mbfS$, define
    \[\{s,s'\} \defeq \begin{cases}
        1 & (s,s') = (p_j,q_j)\text{ for some } j\\
        -1 & (s,s') = (q_j,p_j)\text{ for some } j\\
        0 & \text{otherwise}
    \end{cases}\,.\]
    For $s_1,\cdots,s_k,s'\in \mbfS$, and $w = s_1\cdots s_k$, define
    \[\{[w],s'\} \defeq \sum_{i=1}^{k-1}(-1)^{(|s_1|+\cdots+|s_i|)(|s_{i+1}|+\cdots+|s_k|)}s_{i+1}\cdots s_ks_1\cdots s_{i-1}\{s_i,s'\}\]
    where the sign makes the definition independent of the representative $w\in[w]$. Now extend by $\k$-linearity and the following Leibniz rule (in the second variable):
    \[\{x,yz\} = \{x,y\}z+(-1)^{(|x|+1)|y|}y\{x,z\}\]
    for $x\in \mcalL^\mrm{cyc}$ and $y,z\in\mcalL$. It is easy to see that $\{\cdot,\cdot\}$ preserves $e_\pm$ in the second slot, so it is indeed $\k^m$-linear. The SFT bracket descends to the cyclic quotient\[\{\cdot,\cdot\}:\mcalL^\mrm{cyc}\otimes \mcalL^\mrm{cyc}\to \mcalL^\mrm{cyc} \,.\]
\end{defn}

\begin{rmk}
    Note that the relations $t_it_i^{-1} = t_i^{-1}t_i = e_i$ cause no issues in the above definition, since $\{\cdot,\cdot\}$ vanishes on $t_i^{\pm 1}$.
\end{rmk}

\begin{prop}[{\cite[Proposition 3.4]{Ng10}}] For any $x,y\in \mcalL^\mrm{cyc}$ and $z\in\mcalL$, we have 
    \[\{x,\{y,z\}\} = \{\{x,y\},z\} + (-1)^{(|x|+1)(|y|+1)}\{y,\{x,z\}\}\,.\]
\end{prop}
\subsection{String Differential}
\begin{figure}[H]
    \centering
    \begin{tikzpicture}
        \node[outer xsep = 20pt] (pic1) at (0,0) {\begin{tikzpicture}
    \draw[thick] (1,-1) -- (-1,1);
    \draw[thick] (1,1) -- (.2,.2);
    \draw[thick] (-1,-1) -- (-.2,-.2);
    \draw[color = blue, ultra thick] (-1,.8) -- (.8,-1);
    \path[tips, -Stealth, color = blue, thick] (.3,-.5)--(0, -.2);
        \node at (-.8, -.1) {\footnotesize $\gamma$}; 
\end{tikzpicture}};
        \node[outer xsep = 20pt] (pic2) at (5,0) {\begin{tikzpicture}
    \draw[thick] (1,-1) -- (-1,1);
    \draw[thick] (1,1) -- (.2,.2);
    \draw[thick] (-1,-1) -- (-.2,-.2);
    \draw[color = blue, ultra thick] (-1,.8) -- (-.2,0) .. controls (-.9,-.8) and (-.8,-.9) .. (0,-.2) --(.8,-1);
    \path[tips, -Stealth, color = blue, thick,shift = {(.2,-.2)}] (.4,-.6)--(0, -.2) ;
        \node at (1.2,-.2) {\footnotesize $\delta(\gamma;\tau)$}; 
\end{tikzpicture}};
        \draw[thick, -Stealth] (pic1) -- (pic2); 
    \end{tikzpicture}
    \caption{Local picture of the string differential}
    \label{string differential local}
\end{figure}
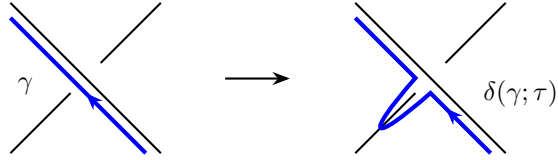

Let $\gamma$ be a composable based string of length $k$ that is generic with holomorphic corners. We say $\tau\in S^1\setminus \{\tau_0,\tau_1,\cdots,\tau_k\}$ is \tit{interior Reeb} for $\gamma$ if $\gamma(\tau)$ is an endpoint of some Reeb chord $a_j$. Given such a $\tau$, define a composable based string $\de(\gamma;\tau)$ of length $k+2$, which is obtained from $\gamma$ by replacing $\tau$ with a jump to the other endpoint of the same Reeb chord and then a jump back. We can homotope $\de(\gamma;\tau)$ so that it becomes generic with holomoprhic corners, and this pertubation is unique up to homotopy through such composable based strings \cite[Section 2.8]{Ng10}. The Lagrangian projection of such a pertubation is depicted in \Cref{string differential local}.
Explicitly, if $w(\gamma) = s_1\cdots s_k$ and $\tau$ lies between $s_i$ and $s_{i+1}$, then
\[w(\de(\gamma;\tau)) = \begin{cases}
    s_1\cdots s_i (p_j q_j) s_{i+1}\cdots s_k & \gamma(\tau) = a_j^-\\
    s_1\cdots s_i (q_j p_j) s_{i+1}\cdots s_k & \gamma(\tau) = a_j^+
\end{cases}\,.\]
Define signs $\ep_1 = \pm 1$ if $\gamma(\tau) = a_j^\mp$, and $\ep_2 = \pm 1$ depending on if $\gamma'(\tau)$ agrees or disagree with the orientation of $\LA$, and 
\[\ep(\gamma;\tau) \defeq (-1)^{|s_1|+\cdots+|s_i|}\ep_1\ep_2\,.\]
We are now ready to make the following definition.

\begin{defn}
\label{string differential}
    The \tit{string differential} \[\de: \mcalL\to \mcalL\] is a degree $-1$ $\k^m$-linear map defined as follows. Given $w\in\mcalW$, choose a generic composable based string $\gamma$ with holomorphic corners such that $w(\gamma) = w$. Then define
    \[\de(w) \defeq\sum_{\tau \text{ interior Reeb for } \gamma} \ep(\gamma;\tau)w(\de(\gamma;\tau))\] and extend by $\k$-linearity. Since $\de$ clearly preserves $e_\pm$, it is indeed $\k^m$-linear. 
\end{defn}

\begin{prop}[{\cite[Proposition 2.26]{Ng23}}] The string differential satisfies
\begin{enumerate}[label = (\arabic*)]
    \item $\de(xy) = (\de x)y + (-1)^{|x|}x(\de y)$ for any $x,y\in\mcalL$;
    \item $\de^2 = 0$.
\end{enumerate}
\end{prop}



\begin{defn}
\label{bullet operator}
    For each base point $\bullet_i$, define a linear map $\bullet_i:\mcalL^\mrm{cyc}\to\mcalL$ as follows. Given a word $w\in \mcalW$, pick a composable based string \[\gamma_w: S^1\setminus\{\tau_0,\tau_1\cdots, \tau_k\}\to \LA\] representing $w$, satisfying $\gamma_w'(\tau)\neq 0$ whenever $\gamma_w(\tau) = \bullet_j$ for some $j$. Let $n_i(\gamma_w)$ denote the number of times $\gamma_w|_{(\tau_k,\tau_1)}$ passes through $\bullet_i$, counted with signs. Now define
    \[\bullet_i([w]) \defeq (-1)^{|w|}\sum_{w'\in [w]} n_i(\gamma_{w'})w'\]
    and extend by linearity. Let $\bullet \defeq \bullet_1+\cdots +\bullet_m$.
\end{defn}

The following generalizes \cite[Proposition 3.8(5)]{Ng10}.

\begin{prop}
\label{string derivation property}
    For any $x\in \mcalL^\mrm{cyc}$ and $y\in\mcalL$, we have
    \[\de\{x,y\} = \{\de x,y\} - (-1)^{|x|}\{x,\de y\} + [\bullet(x),y]\,.\]
\end{prop}
\begin{proof}
The proof is nearly identical to Ng's. The only difference is that we need to keep track of more base points. Let \[f(x,y) = \de\{x,y\}-\{\de x,y\} + (-1)^{|x|}\{x,\de y\}\] and note that  
$[\bullet(x),y] = \mathord{\bullet}_{e_-(y)}(x)y - (-1)^{|x||y|}y\mathord{\bullet}_{e_+(y)}(x)$.

By linearity, we can assume $y\in\mcalW$ and $x = [x']$ for some $x'\in\mcalW_\circ$. Let $\gamma_x,\gamma_y$ be generic strings with holomorphic corners representing $x',y$, respectively. 

Clearly, terms in $f(x,y)$ where $\de$ and $\{\cdot, \cdot\}$ don't interact cancel pairwise. If they do interact, then $\gamma_x$ and $\gamma_y$ must share a segment. By linearity, we can assume they share exactly one segment. Note that no jumps at Reeb chords can occur on the shared segment, since $\gamma_x$ and $\gamma_y$ have holomorphic corners.
  
If the segment does not pass through $\bullet_{e_-(y)}$ or $\bullet_{e_+(y)}$, then there are $9$ possible configurations, each of which gives two canceling terms in $f(x,y)$ (cf. \cite[Figures 3.7 and 3.8]{Ng10}). Since $\bullet_{e_\pm(y)}(x) = 0$, the desired equation holds. 

Otherwise, we have two cases: $\bullet_{e_-(y)} = \bullet_{e_+(y)}$ or $\bullet_{e_-(y)} \neq \bullet_{e_+(y)}$. The first case is no different from having only one base point, which is covered by Ng's proof. For the second case, we can assume the segment passes through one of $\bullet_{e_\pm(y)}$ exactly once but not the other. There are $6$ possible configurations for each of $\bullet_{e_\pm(y)}$, illustrated in \Cref{e- not e+}. Below we analyze one configuration (with a specified orientation) in detail, and we leave the rest to the reader. 
  \setlength{\arrayrulewidth}{1pt}
\begin{figure}[t]
    \centering
        \begin{tabular}{|c|c|c|}
            \hline
             \trimbox{1cm -0.2cm 0cm 1cm}{
\begin{tikzpicture}[scale = 0.8]
\vspace*{-5ex}
    \draw[thick] (-0.3,1.4) -- (-0.3,0) -- (0.3,0);
    \draw[thick] (-1,0) -- (-0.5,0) -- (-0.5,2);
    \draw[dashed,thick,looseness = 3] (-.5,2) to [out = 90, in = 180] (-1,0);
    \draw[fill] (1.4,0) circle (2.4pt) node[below=2pt] {\scriptsize $e_\mp(y)$};
    \draw[fill] (-.3,1.4) circle (2.4pt) node[right=2pt] {\scriptsize $e_\pm(y)$};
    \draw[dashed,thick] (0.3,0) to (1.4,0);
    \node at (-1.2,1) {$x$};  
    \node[anchor = south east] at (-.5, 0) {\footnotesize $s$}; 
    \node[anchor = south west] at (-.3, 0) {\footnotesize $\bar{s}$}; 
\end{tikzpicture}
} & \trimbox{-.2cm .2cm .1cm .2cm}{
\begin{tikzpicture}[scale = 0.8]
    \draw[thick] (-0.3,1.4) -- (-0.3,.5) -- (0.3,.5);
    \draw[thick] (-0.5,0) -- (-0.5,2);
    \draw[dashed, thick, looseness = 2] (-.5,2) to [out = 90, in = 90] (-2, 1);
    \draw[dashed, thick, looseness = 2] (-.5,0) to [out = -90, in = -90] (-2, 1);
    \draw[fill] (1.4,0.5) circle (2.4pt) node[below=2pt] {\scriptsize $e_\mp(y)$};
    \draw[fill] (-.3,1.4) circle (2.4pt) node[right=2pt] {\scriptsize $e_\pm(y)$};
    \draw[dashed,thick] (0.3,.5) to (1.4,.5);
    \node at (-1.2,.9) {$x$};  
\end{tikzpicture}
} & \trimbox{25pt -5pt -1pt 27pt}{
 \begin{tikzpicture}[scale = 0.8]
    \draw[thick] (-0.3,.5) -- (-0.3,-.2);
    \draw[thick] (-1,0) -- (-0.5,0) -- (-0.5,1.5);
    \draw[dashed,thick,looseness = 4] (-.5,1.5) to [out = 90, in = 180] (-1,0);
    \draw[fill] (-0.3,-1) circle (2.4pt) node[right=2pt] {\scriptsize $e_\mp(y)$};
    \draw[fill] (-.3,.5) circle (2.4pt) node[right=2pt] {\scriptsize $e_\pm(y)$};
    \draw[dashed,thick] (-0.3,-.2) to (-0.3,-1);
    \node at (-1.2,.9) {$x$};  
\end{tikzpicture}
} \\
            \hline
            \trimbox{1cm 1cm 0cm -0.2cm}{
\begin{tikzpicture}[scale = 0.8]
    \draw[thick] (-0.3,.6) -- (-0.3,2) -- (0.3,2);
    \draw[thick] (-0.5,0) -- (-0.5,2) -- (-1,2);
    \draw[dashed,thick,looseness = 3] (-1,2) to [out = 180, in = -90] (-.5,0);
    \draw[fill] (1.4,2) circle (2.4pt) node[above=2pt] {\scriptsize $e_\mp(y)$};
    \draw[fill] (-.3,.6) circle (2.4pt) node[right=2pt] {\scriptsize $e_\pm(y)$}; 
    \draw[dashed,thick] (0.3,2) to (1.4,2);
    \node at (-1.2,.8) {$x$};  
    \node[anchor = north east] at (-.5, 2) {\footnotesize $s$}; 
    \node[anchor = north west] at (-.3, 2) {\footnotesize $\bar{s}$}; 
\end{tikzpicture}
} & \trimbox{-.2cm .2cm .1cm .2cm}{
\begin{tikzpicture}[scale = 0.8]
    \draw[thick] (-0.3,.5) -- (-0.3,1.4) -- (0.3,1.4);
    \draw[thick] (-0.5,0) -- (-0.5,2);
    \draw[dashed, thick, looseness = 2] (-.5,2) to [out = 90, in = 90] (-2, 1);
    \draw[dashed, thick, looseness = 2] (-.5,0) to [out = -90, in = -90] (-2, 1);
    \draw[fill] (1.4,1.4) circle (2.4pt) node[above=2pt] {\scriptsize $e_\mp(y)$};
    \draw[fill] (-.3,0.5) circle (2.4pt) node[right=2pt] {\scriptsize $e_\pm(y)$};
    \draw[dashed,thick] (0.3,1.4) to (1.4,1.4);
    \node at (-1.2,.9) {$x$};  
\end{tikzpicture}
} & \trimbox{25pt 27pt -1pt -5pt}{
\begin{tikzpicture}[scale = 0.8]
    \draw[thick] (-0.3,1.5) -- (-0.3,2.2);
    \draw[thick] (-0.5,0.5) -- (-0.5,2) -- (-1,2);
    \draw[dashed,thick,looseness = 4] (-1,2) to [out = 180, in = -90] (-.5,0.5);
    \draw[fill] (-0.3,3) circle (2.4pt) node[right=2pt] {\scriptsize $e_\mp(y)$};
    \draw[fill] (-.3,1.5) circle (2.4pt) node[right=2pt] {\scriptsize $e_\pm(y)$}; 
    \draw[dashed,thick] (-0.3,2.2) to (-0.3,3);
    \node at (-1.2,1) {$x$};  
\end{tikzpicture}
} \\
            \hline
        \end{tabular}
        \caption{The $6$ possible configurations when the shared segment passes one of $e_\pm(y)$ exactly once.}
    \label{e- not e+}
\end{figure}
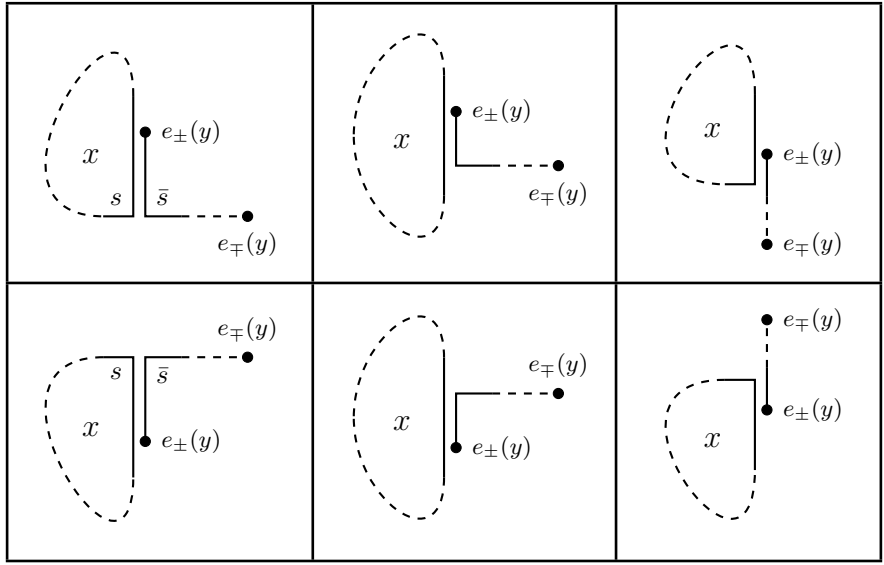

  Consider the configuration depicted in \Cref{1shared_1}. In this case, the only interaction term occurs in $\de\{x,y\}$.
  We can write $x = [w_1st_{e_-(y)}]$ and $y = t_{e_-(y)}^{-1}\bar{s}w_2$ for some words $w_1,w_2$. Then we have $\bullet_{e_+(y)}(x) = 0$ and
  \[[\bullet(x),y] = \bullet_{e_-(y)}(x)y = (-1)^{|x|}w_1st_{e_-(y)}t_{e_-(y)}^{-1}\bar{s}w_2 = (-1)^{|x|}w_1s\bar{s}w_2\,.\]
 The contribution to $\{x,y\}$ is
  \[\{s,\bar{s}\}t_{e_-(y)}^{-1}t_{e_-(y)}w_1w_2 = \{s,\bar{s}\}w_1w_2\]
  so the contribution to $\de\{x,y\}$ is
  \[(-1)^{|w_2|}\ep_1\ep_2\{s,\bar{s}\}w_1s\bar{s}w_2 = (-1)^{|x|}w_1s\bar{s}w_2 = [\bullet(x),y]\]
  since $\ep_1 = \{s,\bar{s}\}$ and $\ep_2 = (-1)^{|s|}$. (The signs $\ep_1,\ep_2$ are defined just above \Cref{string differential}.) Hence, we have shown $f(x,y) = [\bullet(x),y]$ for this configuration.
\end{proof}
\begin{figure}
    \centering
    \begin{tikzpicture}
    \node[outer sep = 0 ] (pic1) at (0,0) {\trimbox{1cm 0cm 0cm 1cm}{
\begin{tikzpicture}[scale = 0.8]
    \draw[-Straight Barb,thick] (-0.3,1.4) -- (-0.3,0) -- (0.3,0);
    \draw[{Straight Barb[reversed]}-,thick] (-1,0) -- (-0.5,0) -- (-0.5,2);
    \draw[dashed,thick,looseness = 3] (-.5,2) to [out = 90, in = 180] (-1,0);
    \draw[fill] (1.4,0) circle (2.4pt) node[below=2pt] {\scriptsize $e_+(y)$};
    \draw[fill] (-.3,1.4) circle (2.4pt) node[right=2pt] {\scriptsize $e_-(y)$};
    \node at (-.3,.9) {\scriptsize $\bigstar$}; 
    \draw[dashed,thick] (0.3,0) to (1.4,0);
    \node at (-1.2,1) {$x$};  
    \node[anchor = south east] at (-.5, 0) {\footnotesize $s$}; 
    \node[anchor = south west] at (-.3, 0) {\footnotesize $\bar{s}$}; 
\end{tikzpicture}
}};
    \node[outer sep = 0 ] (pic2) at (5,0) {\trimbox{1cm 0cm 0cm 1cm}{
\begin{tikzpicture}[scale = 0.8]
    \draw[-Straight Barb,thick] (-1,0) -- (0.3,0);
    \draw[thick,rounded corners = 3.5mm] (-0.5,2) -- (-0.4,.3) -- (-0.1,1.4);
    \draw[dashed,thick,looseness = 3] (-.5,2) to [out = 90, in = 180] (-1,0);
    \draw[fill] (1.4,0) circle (2.4pt) node[below=2pt] {\scriptsize $e_+(y)$};
    \draw[fill] (-.1,1.4) circle (2.4pt) node[right=2pt] {\scriptsize $e_-(y)$};
    \node at (-.22,1.05) {\scriptsize $\bigstar$}; 
    \draw[dashed,thick] (0.3,0) to (1.4,0);
    \node at (-1.2,1) {$x$}; 
    \draw (-.4,0) circle (2.4pt);
\end{tikzpicture}
}};
    \node[outer sep = 0 ] (pic3) at (10,0) {\trimbox{1cm 0cm 0cm 1cm}{
\begin{tikzpicture}[scale = 0.8]
    \draw[thick] (-1,0) -- (-0.5,0); 
    \draw[thick, rounded corners = 3mm] (-0.5,0) -- (-.4,.5) --(-.3,0);
    \draw[-Straight Barb,thick] (-.3,0) -- (.3,0);
    \draw[thick,rounded corners = 3.5mm] (-0.5,2) -- (-0.4,.35) -- (-0.1,1.4);
    \draw[dashed,thick,looseness = 3] (-.5,2) to [out = 90, in = 180] (-1,0);
    \draw[fill] (1.4,0) circle (2.4pt) node[below=2pt] {\scriptsize $e_+(y)$};
    \draw[fill] (-.1,1.4) circle (2.4pt) node[right=2pt] {\scriptsize $e_-(y)$};
    \node at (-.2,1.05) {\scriptsize $\bigstar$}; 
    \draw[dashed,thick] (0.3,0) to (1.4,0);
    \node at (-1.2,1) {$x$}; 
    \node[anchor = south east] at (-.5, 0) {\footnotesize $s$}; 
    \node[anchor = south west] at (-.3, 0) {\footnotesize $\bar{s}$}; 
\end{tikzpicture}
}};
    \draw[-Stealth,thick,transform canvas={yshift=-3mm}] (pic1) -- (pic2);
    \draw[-Stealth,thick,transform canvas={yshift=-3mm}] (pic2) -- (pic3);
    \node (SFT) at (2.5,0) {\footnotesize$\{\cdot,\cdot\}$};
    \node (str) at (7.5,0) {\footnotesize$\de$};
    \end{tikzpicture}
    \caption{A configuration where the shared segment passes through $e_-(y)$ exactly once. The hollow circle in the middle picture denotes the interior Reeb chord where $\de$ will be applied.}
    \label{1shared_1}
\end{figure}
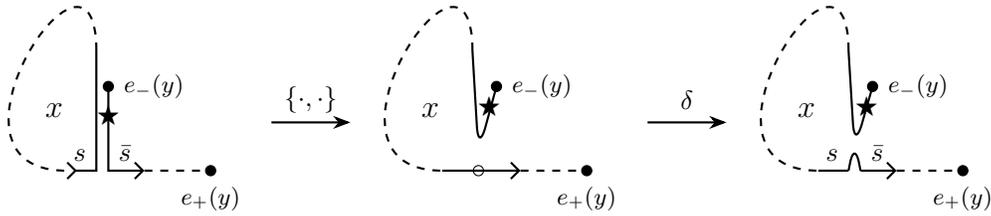


\subsection{Hamiltonian, Differential, and Curvature}
\label{LSFT differential}
 Let $D^2$ be the unit disk with the standard orientation. For $k\geq 1$, pick cyclically ordered points $\tau_1,\cdots,\tau_k\in\del D^2$, called \tit{boundary punctures}. A \tit{disk with $k$ corners} is a continuous orientation-preserving map (considered up to domain reparametrization fixing $\tau_1,\cdots,\tau_k$)
 \[\DE:D^2\to \R^2\]
satisfying
\begin{enumerate}[label = (\arabic*)]
    \item $\DE(\del D^2)\inn \LA_{xy}$;
    \item $\DE|_{D^2\setminus\{\tau_1,\cdots,\tau_k\}}$ is an immersion;
    \item $\DE$ sends each $\tau_i$ to a crossing of $\LA_{xy}$ and a neighborhood of $\tau_i$ to a quadrant at that crossing, called the \tit{$i$-th corner} of $\DE$.
\end{enumerate}
Let $\DE(\LA)$ denote the set of all disks, which has an equivalence relation given by cyclic permutation of the boundary punctures. Let $\DE^\mrm{cyc}(\LA)$ be the set of equivalence classes, and write $[\DE]\in \DE^\mrm{cyc}(\LA)$ for the equivalence class of $\DE$. A choice of representative in $[\DE]$ is the same as picking a linear ordering of the boundary punctures.
\begin{rmk}
    Outside of this subsection, ``disks'' usually mean disks up to cyclic permutation of the punctures. But here we are trying to be very careful with signs.
\end{rmk}
The boundary $\del\DE\defeq \DE|_{\del D^2\setminus\{\tau_1,\cdots,\tau_k\}}$ of a disk $\DE$ always has left-turning corners at crossings, since $\DE$ is orientation-preserving. Hence, by \Cref{holomorphic corner rmk}, $\del\DE$ is the projection of a generic cyclically composable string $\gamma_{\del\DE}$ with holomorphic corners. Now define the cyclically composable word
\[w(\DE) \defeq w(\gamma_{\del\DE})\in \mcalW_\circ\] which can be read off from the quadrants at crossings covered by $\DE$. 

\begin{figure}[H]
\centering
\begin{subfigure}{.4\textwidth}
  \centering
  \begin{tikzpicture}
      \fill[fill=gray!60] (0,0) --  (45:.8) arc(45:-225:.8) -- cycle;
      \draw[ultra thick, Stealth-] (1,-1) -- (-1,1);
      \draw[ultra thick, Stealth-] (1,1) -- (.2,.2);
      \draw[ultra thick] (-1,-1) -- (-.2,-.2);
      \node (p) at (.5,0) {$p$};
      \node (p) at (-.5,0) {$p$};
      \node (q) at (0,.5) {$q$};
      \node (q) at (0,-.5) {$q$};
  \end{tikzpicture}
  \caption{positive crossing}
  \label{Hamiltonian signs positive crossing}
\end{subfigure}%
\begin{subfigure}{.4\textwidth}
  \centering
  \begin{tikzpicture}
      \fill[fill=gray!60] (0,0) --  (45:.8) arc(45:-45:.8) -- cycle;
      \draw[ultra thick, -Stealth] (1,-1) -- (-1,1);
      \draw[ultra thick, Stealth-] (1,1) -- (.2,.2);
      \draw[ultra thick] (-1,-1) -- (-.2,-.2);
      \node (p) at (.5,0) {$p$};
      \node (p) at (-.5,0) {$p$};
      \node (q) at (0,.5) {$q$};
      \node (q) at (0,-.5) {$q$};
  \end{tikzpicture}
  \caption{negative crossing}
  \label{Hamiltonian signs negative crossing}
\end{subfigure}
\caption{Orientation signs for the Hamiltonian}
\label{Orientation signs for the Hamiltonian}
\end{figure}
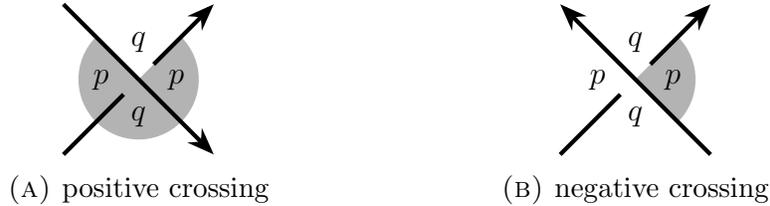

Also define signs
\[\ep_i(\DE)\defeq \begin{cases}
    +1 & i\text{-th corner of }\DE\text{ is unshaded}\\
    -1 & i\text{-th corner of }\DE\text{ is shaded}
\end{cases}\]
according to \Cref{Orientation signs for the Hamiltonian}, and $\ep'(\DE)=\pm 1$ based on if the orientation of $\del\DE$ immediately after $\tau_k$ agrees or disagrees with $\LA_{xy}$. Now set
\begin{equation}
    \label{sign of disk}
    \sgn(\DE) \defeq \ep_1(\DE)\cdots\ep_k(\DE)\ep'(\DE)
\end{equation}
and define the signed word associated to a disk
\begin{equation}
\label{w tilde}
\tilde{w}(\DE) \defeq \sgn(\DE)w(\DE)\in\mcalL\,. \end{equation}
It is shown in \cite[Lemma 2.14]{Ng23} that if $[\DE] = [\DE']$, then $[\tilde{w}(\DE)] = [\tilde{w}(\DE')]$, so we can define \[\tilde{w}([\DE]) \defeq [\tilde{w}(\DE)]\in \mcalL^\mrm{cyc}\,.\]

\begin{deflem}[{\cite{Ng23}}]
\label{Hamiltonian}
The \tit{Hamiltonian} 
\[h \defeq \sum_{[\DE]\in\DE^\mrm{cyc}(\LA)} \tilde{w}([\DE])\in \mcalL^\mrm{cyc}\]
has degree $-2$ and lies in $\mcalF^1\mcalL^\mrm{cyc}$.
\end{deflem}

\begin{prop}[{\cite{Ng23}}] 
\label{quantum master equation}
The Hamiltonian satisfies the \tit{quantum master equation}
    \[\de h = \frac{1}{2}\{h,h\}\,.\]
\end{prop}


\begin{defn}
    The \tit{LSFT differential} is the degree $-1$ $\k^m$-linear map
    \[d\defeq \{h,\cdot\} - \de:\mcalL\to \mcalL\]
    Denote the part involving disks by $d_\DE \defeq \{h,\cdot\}$.
\end{defn}

We conclude the following, generalizing the curvature formula in \cite{Ng10} (but recall that we are following the sign conventions of \cite{Ng23}).

\begin{prop}
\label{curvature formula}
    The composable LSFT algebra $(\mcalL, d)$ is a filtered curved dg $\k^m$-algebra with curvature
    \[F = -\mathord{\bullet}(h)\in \mcalF^1\mcalL\] where $\bullet$ is defined in \Cref{bullet operator}.
\end{prop}
\begin{proof}
Since $h\in\mcalF^1\mcalL^{\mrm{cyc}}$, it is easy to see that $d_\DE$ preserves the filtration by definition of $\{\cdot,\cdot\}$. By construction, $\de$ also preserves the filtration. Thus, so does $d = d_\DE - \de$.

    Given any $x\in \mcalL$, we have
    \begin{align*}
        d^2x &= \{h,\{h,x\}\}-\{h,\de x\} - \de \{h,x\}+\de^2 x\\
             &= \frac{1}{2}\{\{h,h\},x\} -\{\de h, x\} - [\bullet(h), x] = [F,x]
    \end{align*} 
where the second equality follows from the Jacobi identity, $\de^2 = 0$, and \Cref{string derivation property}, and the third follows from \Cref{quantum master equation}.  
\end{proof}

For computational purposes, we rewrite $d_\DE$ as follows. For any $s\in \mbfS$, define
\begin{equation}\label{disk with specified final corner}\DE(\LA;s)\defeq \{\DE \in \DE(\LA)\mid w_k(\gamma_{\del\DE}) = s\text{, where $k$ is the number of corners of $\DE$}\}\,.\end{equation} For any cyclically composable string $\gamma$ of length $k$, define 
\begin{equation}
\label{w'}
    w'(\gamma) \defeq w_{k1}(\gamma) w_1(\gamma)w_{12}(\gamma)\cdots w_{k-1,k}(\gamma)
\end{equation} where $w_i(\gamma)$ and $w_{i,i+1}(\gamma)$ are specified by \eqref{wi of gamma} and \eqref{w i i+1 of gamma}. For any disk $\DE$, define $w'(\DE) \defeq w'(\gamma_{\del\DE})$. Then by definitions of the Hamiltonian and $\{\cdot,\cdot\}$, we have:
\begin{prop}
\label{disk diff}
On generators, $d_\DE$ is given by
	\begin{align*}
		d_\DE  q_j &= \sum_{\DE\in \DE(\LA;p_j)}\sgn(\DE)w'(\DE)\\
		d_\DE  p_j &= -\sum_{\DE\in \DE(\LA;q_j)}
		\sgn(\DE)w'(\DE)\\
		d_\DE  t_i &= d_\DE  t_i^{-1} = 0
	\end{align*}
where $\sgn(\DE)$ is defined in \eqref{sign of disk}. 
\end{prop}
The orientation signs for the Hamiltonian (\Cref{Orientation signs for the Hamiltonian}) and the Chekanov-Eliashberg differential (\Cref{Orientation signs for CE}) are chosen to be compatible:
\begin{prop}[{\cite{Ng23}}]
\label{orientation signs are compatible}
    Let $\DE_1(\LA)\inn \DE(\LA)$ be the set of disks with a single positive puncture and let $\DE_1^\mrm{cyc}(\LA)\inn \DE^\mrm{cyc}(\LA)$ be the corresponding set of equivalence classes of disks. Define
    \[h_1 \defeq \sum_{[\DE]\in\DE_1^\mrm{cyc}(\LA)} \tilde{w}([\DE])\in \mcalL^\mrm{cyc}\,.\] Then we have
    \[\del = \{h_1,\cdot\}\,.\]
\end{prop}
\begin{rmk}
    \label{single positive puncture disks rmk} The set of disks $\DE_1(\LA;p_j)$, which appears in \Cref{CE diff def}, is more precisely given by $\DE_1(\LA;p_j)\defeq \DE_1(\LA) \cap \DE(\LA;p_j)$, using the notation \eqref{disk with specified final corner}.
\end{rmk}
Below we consider a situation where the curvature is the simplest.
\begin{defn}
    A Legendrian knot $\LA$ with a single base point $\bullet$ is \tit{cusp-pointed} if $\LA_{xy}$ is the resolution of a front diagram and $\bullet$ lies on one of the resolved rightmost cusps. The two variables associated to this resolved cusp will be denoted by $p_\bullet,q_\bullet$.
\end{defn}
\begin{figure}[H]
    \centering
    \trimbox{5pt 50pt 17pt 50pt}{
\begin{tikzpicture} 
   \begin{knot}[
     consider self intersections = true,
     end tolerance = 1pt,
     clip width = 5pt,
   ]
   \strand[ultra thick, rounded corners] (4,1).. controls (6.5,4.2) and (6.5, -.8)..(4,2.4);
   \flipcrossings{1};
    \end{knot} 
    \fill[color = gray!80, even odd rule, rounded corners] (4.68,1.7).. controls (6.25,3.1) and (6.25, 0.3)..(4.68,1.7);
   \node at (5.87,1.7) {\footnotesize$\bigstar$};
   \node at (5.6,2.1) {$\bullet$};
    \node[below = 7pt] at (4.7,1.7) {$q_\bullet$};
	\node[left = 7pt] at (4.7,1.7) {$p_\bullet$};     
\end{tikzpicture}
}
    \caption{Local picture of a cusp-pointed Legendrian knot}
    \label{cusp-pointed}
\end{figure}
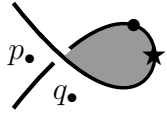
\begin{lemma}
\label{cusp-pointed curvature}
If $\Lambda$ is cusp-pointed, then the curvature of $\mcalL(\Lambda)$ is $-p_\bullet t^{\pm 1}$, depending on the orientation of $\LA$.     
\end{lemma}
\begin{proof}
    Because the base point lies on a resolved rightmost cusp, the only disk passing through the base point is the shaded disk shown in \Cref{cusp-pointed}. Hence, the Hamiltonian is $[p_\bullet t^{\pm 1}]$, depending on the orientation of $\LA$. The lemma now follows from \Cref{curvature formula}.
\end{proof}

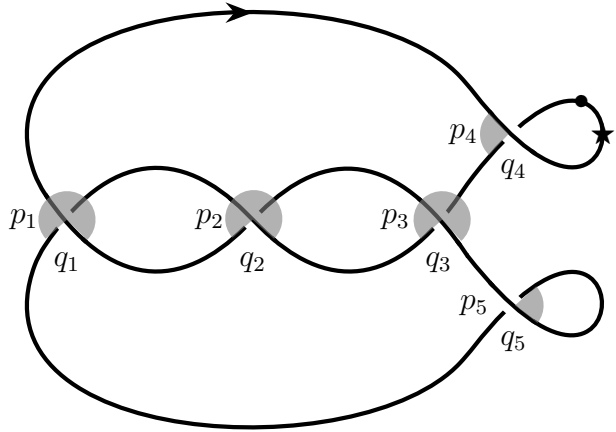
\begin{figure}[H]
    \centering
    \trimbox{17pt 15pt 15pt 15pt}{
 \begin{tikzpicture}[scale = .95]
	\node at (-1.6,-.1) {$q_1$};
	\node at (1,-.1) {$q_2$};
	\node at (3.6,-.1) {$q_3$};
	\node at (4.65,1.2) {$q_4$};
	\node at (4.65,-1.2) {$q_5$};
	
	\node at (-2.2,.5) {$p_1$};
	\node at (.4,.5) {$p_2$};
	\node at (3,.5) {$p_3$};
	\node at (3.95,1.7) {$p_4$};
	\node at (4.1,-.7) {$p_5$};
   
    \begin{knot}[
      consider self intersections = true,
      end tolerance = 1pt,
      clip width = 5pt,
    ]
    \strand[ultra thick, rounded corners] (4,0) .. controls (6.5,-3.2) and (6.5,1.8).. (4,-1.4) .. controls (2.5,-3) and (-3,-2.7) .. (-2,0) .. controls (-1,1.5) and (0,1.5) .. (1,.5) .. controls (2,-.5) and (3,-.5) .. (4,1).. controls (6.5,4.2) and (6.5, -.8)..(4,2.4) .. controls (2.5,4) and (-3,3.7) .. (-2,1) .. controls (-1,-.5) and (0,-.5) .. (1,.5)[sharp corners] .. controls (2,1.5) and (3,1.5).. cycle;
    \flipcrossings{2,4,5};
     \end{knot} 
     \fill[fill=gray, opacity = 0.6, shift =({1,.5})] (0,0) --  (-42:.4) arc (-45:222:.4) -- cycle;
     \fill[fill=gray, opacity = 0.6, shift =({3.65,.5})] (0,0) --  (-48:.4) arc (-48:222:.4) -- cycle;
     \fill[fill=gray, opacity = 0.6, shift =({-1.6,.5})] (0,0) --  (-45:.4) arc (-45:222:.4) -- cycle;
     \fill[fill=gray, opacity = 0.6, shift =({4.58,1.7})] (0,0) --  (133:.4) arc (133:227:.4) -- cycle;
     \fill[fill=gray, opacity = 0.6, shift =({4.67,-.7})] (0,0) --  (-38:.4) arc (-38:42:.4) -- cycle;     
    \node at (5.9,1.7) {\footnotesize$\bigstar$};
    \node at (5.6,2.15) {$\bullet$};
    \path[tips, -Stealth, ultra thick] (0,3.4) -- (1,3.4);
\end{tikzpicture}
}
    \caption{Lagrangian projection of the trefoil}
    \label{trefoil}
\end{figure}
\begin{exmp}
\label{disk diff of trefoil}
Let us compute $d_\DE$ for the trefoil, depicted in \Cref{trefoil}. The generators have the following gradings
\begin{align*}&|q_1|=|q_2|=|q_3|=0\,,\ |q_4|=|q_5|=1\,,\\ &|p_1|=|p_2|=|p_3|=-1\,,\ |p_4|=|p_5|=-2\,,\ |t|=0\,.\end{align*}
By \Cref{disk diff}, we have
\begin{table}[H]
    \centering
    \begin{tabular}{ll}
        $d_\DE q_1 = -p_2$ & $d_\DE p_1 = -q_2q_3p_4 - p_4 - p_5 - p_5q_3q_2$ \\
 		$d_\DE q_2 = p_1 - p_3$ &  $d_\DE p_2 = -q_3p_4q_1 - q_1p_5q_3$\\
 		$d_\DE q_3 = p_2$ &   $d_\DE p_3 = -p_4q_1q_2 - p_4 - p_5 - q_2q_1p_5$\\
 		$d_\DE q_4 = -q_1q_2q_3 - q_3 - q_1 + t$& $d_\DE p_4 = 0$\\
 		$d_\DE q_5 = q_1 + q_3 + q_3q_2q_1 + 1$ & $d_\DE p_5 = 0$\\
  		&$d_\DE t = d_\DE t^{-1} = 0$
    \end{tabular}
\end{table}
\noindent The curvature is $F=-p_4t$, by \Cref{cusp-pointed curvature}.
\end{exmp}
\subsection{Curved Augmentations} 
\label{curved augmentations}
Augmentations of the Chekanov-Eliashberg algebra provide a computable way to access the geometric information of pseudo-holomorphic disks. In the same vein, we define curved augmentations of the composable LSFT algebra, which will give us access to disks with multiple positive punctures.

Let $\LA$ be a Legendrian knot with a single marked point. The $\Z$-grading of $\mcalL(\LA)$ induces a $\Z/2\rot(\LA)$-grading, whence $|t|=|t^{-1}|=0$ by \Cref{grading prop}.
\begin{defn} Let $c\in \k$. A \tit{$c$-curved augmentation} of $\LA$ is a map of filtered curved differential $\Z/2\rot(\LA)$-graded algebras \[\tep:(\mcalL(\LA),d,F)\to (\k\bbra{\eta},0,c\eta)\] where $|\eta| = -2$ and $\k\bbra{\eta}$ is filtered by powers of $\eta$. Note that a map of curved dga's usually includes a degree $-1$ element of the target, but here the target is concentrated in even degrees. Explicitly, $\tep$ needs to satisfy the equations:
\[\begin{cases}
	\tep\circ d=0\\
	\tep(F) = c\eta
\end{cases}\hspace{-1ex}.\]
By the below lemma, the first equation is equivalent to $\tep\circ d_\DE = 0$.
\end{defn}
\begin{rmk}
    For a general pointed Legendrian link, $\mcalL$ is a $\k^m$-algebra where $m$ is the number of marked points, in which case the correct target of a curved augmentation might be $\Mat_m(\k)\bbra{\eta}$, but we will not consider this generalization in this paper.
\end{rmk}
\begin{lemma}
\label{string differential augments to zero}
    Any map of graded algebras $\tep:\mcalL\to \k\bbra{\eta}$ satisfies $\tep\circ \de = 0$, where $\de$ is the string differential.
\end{lemma}
\begin{proof}
    For any $s\in \mbfS$, every term in $\tep\circ\de(s)$ is of the form $\tep(p_i)\tep(q_i)\tep(s)$ for some $i$. Since $|p_i| + |q_i| = -1$ and $\k\bbra{\eta}$ is concentrated in even degrees, we must have $\tep(p_i)\tep(q_i) = 0$, so the lemma follows.
\end{proof}
Every $c$-curved augmentation of the LSFT algebra descends to a usual augmentation of the Chekanov-Eliashberg algebra:
\[\begin{tikzcd}
\mcalL\rar{\tep} \dar & \k\bbra{\eta}\dar\\
\mcalA\rar{\ep} & \k
\end{tikzcd}\]
since $(\mcalA,\del)$ is the zeroth associated graded of $(\mcalL,d)$. For $c = 0$, there is a canonical lift: 
\begin{prop}
\label{0-curved augmentations}
    Any augmentation $\ep:\mcalA\to \k$ can be extended trivially to a $0$-curved augmentation $\tep:\mcalL\to \k$ by setting $\tep(q_j) = \ep(q_j)$ for each $q_j$, and $\tep(p_j) = 0$ for each $p_j$.
\end{prop}
\begin{proof}
    By definition, $\tep$ vanishes on $\mcalF^1\mcalL$, so $\tep(F) = 0$ by \Cref{curvature formula}. To prove $\tep\circ d = 0$, it suffices to show that the image of $d - \del = (d_\DE - \del) - \de$ lies in $\mcalF^1\mcalL$, since $\ep\circ\del = 0$. The definition of $\de$ clearly implies $\image \de\inn \mcalF^1\mcalL$. And the difference $d_\DE - \del$ counts disks with at least two positive corners, so $\image (d_\DE - \del)\inn \mcalF^1\mcalL$ as well. 
\end{proof}
\begin{exmp}
    Consider the Legendrian unknot depicted in \Cref{unknot}. The LSFT generators have gradings $|q| = 1$, $|p| = -2$, and $|t| = 0$. By \Cref{disk diff}, we have
     \begin{align*}
     &d_\DE	q = 1+t\\
     &d_\DE p = d_\DE t = d_\DE t^{-1} = 0\,.
     \end{align*}
The curvature is $F= -pt$, by \Cref{cusp-pointed curvature}. For any $c\in\k$, it is easy to check that $\tep:\mcalL\to\k\bbra{\eta}$ given by
\[\tep(t) = - 1\,,\ \tep(q) = 0\,,\ \tep(p) = c\eta\] is the unique $c$-curved augmentation.
\end{exmp}
\begin{exmp}
\begin{figure}[H]
        \centering
        \trimbox{0pt 85pt 20pt 83pt}{
\begin{tikzpicture} 
   \begin{knot}[
     consider self intersections = true,
     end tolerance = 1pt,
     clip width = 5pt,
   ]
   \strand[ultra thick, rounded corners] (1,.5) .. controls (4.2,-4) and (4.2,4) .. (1,-.5);
   \strand[ultra thick, looseness = 1.5](1,.5) to [out = 125, in = 90] (-.5,0);
    \strand[ultra thick, looseness = 1.5](1,-.5) to [out = -125, in = -90] (-.5,0);
    \end{knot} 
       \node at (3.4,0) {\footnotesize$\bigstar$};
   \node at (3,.85) {$\bullet$};
   \path[tips, -Stealth, ultra thick] (0,.9) -- (.5,.9);
   \fill[fill=gray, opacity = 0.6, shift =({1.35,0})] (0,0) --  (129:.4) arc (129:231:.4) -- cycle;
    \node[below = 7pt] at (1.4,0) {$q$};
	\node[right = 7pt] at (1.4,0) {$p$};     
\end{tikzpicture}
}
        \caption{Max-tb unknot}
        \label{unknot}
\end{figure}
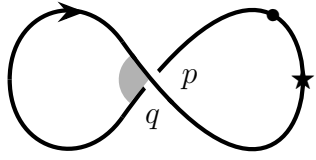
We computed $d_\DE$ of the trefoil in \Cref{disk diff of trefoil}. By killing all terms with filtration degree $\geq 1$, we can read off the Chekanov-Eliashberg differential $\del$, and one can check that $\ep:\mcalA\to\k$ given below is an augmentation:
 \[\ep(q_1)=\ep(t)=-1, \ep(q_2)=\ep(q_3)=\ep(q_4)=\ep(q_5)=0\,.\]
Let us find $c$-curved augmentations $\tep:\mcalL\to\k$ lifting $\ep$. Inserting the grading constraints \[\tilde{\ep}(q_4) = \tilde{\ep}(q_5) = \tilde{\ep}(p_1)=\tilde{\ep}(p_2)=\tilde{\ep}(p_3)=0\] into $d_\DE$ simplifies the equation $\tilde{\ep} \circ d_\DE = 0$ to 
 \[\tilde{\ep}(p_4)+\tilde{\ep}(p_5)=0\,.\]
 Hence, for any $c\in\k$, there is a unique $c$-curved augmentation $\tilde{\ep}$ lifting $\ep$, given by $\tilde{\ep}(p_4) = -\tilde{\ep}(p_5) = c\eta$.
 \end{exmp}
Computational evidence suggests the following:
 \begin{conj}
     Every augmentation lifts to a $c$-curved augmentation, for any $c\in\k$.
 \end{conj}


\section{\texorpdfstring{$\Aug_+$}{Aug+} and Its Bimodules}
\label{aug+ bimods}
This section is dedicated to the ingredients and proof of \Cref{intro A infty EES}.

\subsection{\texorpdfstring{$m$}{m}-Copies and Consistency}
\label{m-copy}
Let $\LA\inn \R^3$ be a pointed Legendrian link. Label its marked points by $\star_1,\cdots,\star_r$, base points by $\bullet_1,\cdots,\bullet_r$, and Reeb chords by $a_1,\cdots,a_n$. Let $f : \Lambda\to\R$ be a Morse function with local maxima $b_1,\cdots,b_\ell$ and local minima $c_1,\cdots,c_\ell$, all of which avoid the marked points and base points.

Pick a Maslov potential $\mu:\{\bullet_1,\cdots,\bullet_r\}\to\R$ so that $|t_i| \equiv 0 \mod 2\rot(\LA)$ for each $i$, which exists by \Cref{rotation number base point grading}. The $\Z$-grading of $\mcalA(\LA)$ is then specified by $\eqref{grading def}$, but in this section, we will consider $\mcalA(\LA)$ with the induced $\Z/2\rot(\LA)$-grading.

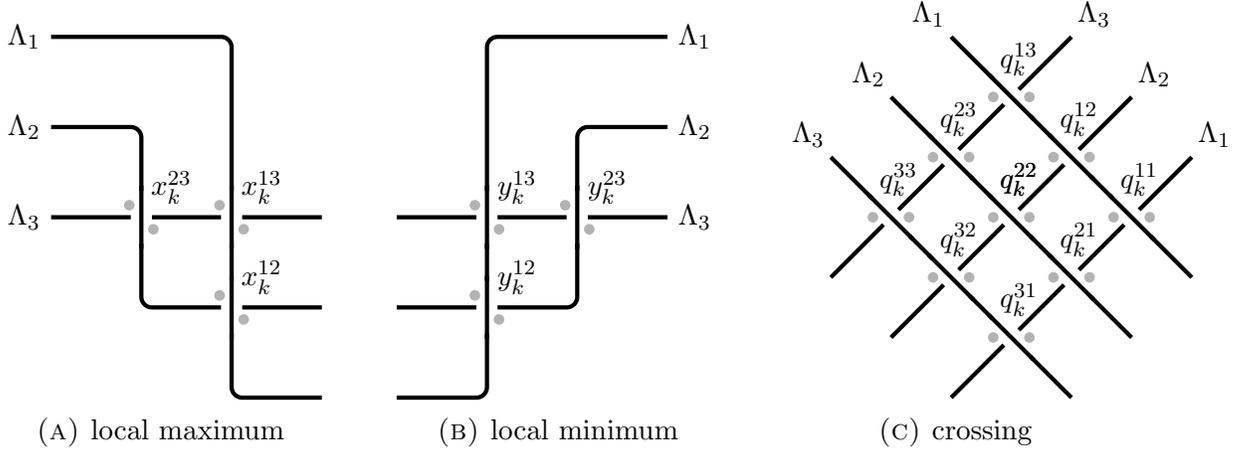
\begin{figure}[t]
    \centering
    \begin{subfigure}[t]{.28\textwidth}
        \centering
        \begin{tikzpicture}[scale = .65]
        \begin{knot}[
     consider self intersections = true,
     end tolerance = 1pt,
     clip width = 5pt,
     draft mode = off
   ]
            \strand[ultra thick, rounded corners] (-3,3) -- (0,3) -- (0,-3) -- (1.5,-3);
            \strand[ultra thick, rounded corners] (-3,1.5) -- (-1.5,1.5) -- (-1.5,-1.5) -- (1.5,-1.5);
            \strand[ultra thick] (-3,0) -- (1.5,0);
      \end{knot} 
            \node[anchor = east] at (-3,3) {\small$\Lambda_1$};
            \node[anchor = east] at (-3,1.5) {\small$\Lambda_2$};
            \node[anchor = east] at (-3,0) {\small$\Lambda_3$};
            \node[shift = {(.4,.4)}] at (0,0) {\small$x_k^{13}$};
            \node[shift = {(.4,.4)}] at (-1.5,0) {\small$x_k^{23}$}; 
            \node[shift = {(.4,.4)}] at (0,-1.5) {\small$x_k^{12}$};
            \fill[gray!60,shift={(.2,-.2))}] (0,0) circle(2.5pt);
            \fill[gray!60,shift={(-.2,.2))}] (0,0) circle(2.5pt);
            \fill[gray!60,shift={(.2,-.2))}] (-1.5,0) circle(2.5pt);
            \fill[gray!60,shift={(-.2,.2))}] (-1.5,0) circle(2.5pt); 
            \fill[gray!60,shift={(.2,-.2))}] (0,-1.5) circle(2.5pt);
            \fill[gray!60,shift={(-.2,.2))}] (0,-1.5) circle(2.5pt); 
        \end{tikzpicture}
        \caption{local maximum}
    \end{subfigure}
    \centering
    \begin{subfigure}[t]{.28\textwidth}
        \centering
                \begin{tikzpicture}[scale = .65]
        \begin{knot}[
     consider self intersections = true,
     end tolerance = 1pt,
     clip width = 5pt,
     draft mode = off
   ]
   \strand[ultra thick, rounded corners] (-1.5,-3) -- (0,-3) -- (0,3) -- (3,3);
   \strand[ultra thick, rounded corners] (-1.5,-1.5) -- (1.5,-1.5) -- (1.5,1.5) -- (3,1.5);
   \strand[ultra thick] (-1.5,0) -- (3,0);          
  \end{knot}
            \node[anchor = west] at (3,3) {\small$\Lambda_1$};
            \node[anchor = west] at (3,1.5) {\small$\Lambda_2$};
            \node[anchor = west] at (3,0) {\small$\Lambda_3$};
            \node[shift = {(.4,.4)}] at (0,0) {\small$y_k^{13}$};
            \node[shift = {(.4,.4)}] at (1.5,0) {\small$y_k^{23}$}; 
            \node[shift = {(.4,.4)}] at (0,-1.5) {\small$y_k^{12}$};
            \fill[gray!60,shift={(.2,-.2))}] (0,0) circle(2.5pt);
            \fill[gray!60,shift={(-.2,.2))}] (0,0) circle(2.5pt);
            \fill[gray!60,shift={(.2,-.2))}] (1.5,0) circle(2.5pt);
            \fill[gray!60,shift={(-.2,.2))}] (1.5,0) circle(2.5pt); 
            \fill[gray!60,shift={(.2,-.2))}] (0,-1.5) circle(2.5pt);
            \fill[gray!60,shift={(-.2,.2))}] (0,-1.5) circle(2.5pt);          
        \end{tikzpicture}
        \caption{local minimum}
    \end{subfigure}
    \centering
    \begin{subfigure}[t]{.38\textwidth}
        \centering
            
    \begin{tikzpicture}[scale = .65]
	\begin{knot}[
     consider self intersections = true,
     end tolerance = 1pt,
     clip width = 5pt,
     draft mode = off
   ]
    \strand[ultra thick](-2,2)--(2,-2);
	\strand[ultra thick,shift={(1,1)}](-2,2)--(2,-2);	
	\strand[ultra thick,shift={(-1,-1)}](-2,2)--(2,-2);
	\strand[ultra thick](-2,-2)--(2,2);
	\strand[ultra thick,shift={(-1,1)}](-2,-2)--(2,2);
	\strand[ultra thick,shift={(1,-1)}](-2,-2)--(2,2);
	\end{knot}
	\node[shift = {(-.3,.3)}] at (-1,3) {\small$\Lambda_1$};
	\node[shift = {(-.3,.3)}] at (-2,2) {\small$\Lambda_2$};
	\node[shift = {(-.3,.3)}] at (-3,1) {\small$\Lambda_3$};
	\node[shift = {(.3,.3)}] at (1,3) {\small$\Lambda_3$};
	\node[shift = {(.3,.3)}] at (2,2) {\small$\Lambda_2$};
	\node[shift = {(.3,.3)}] at (3,1) {\small$\Lambda_1$};
	\node[shift = {(.1,.5)}] at (0,0) {\small$q_k^{22}$};
	\node[shift = {(.1,.5)}] at (2,0) {\small$q_k^{11}$};
	\node[shift = {(.1,.5)}] at (-2,0) {\small$q_k^{33}$};
	\node[shift = {(.1,.5)}] at (0,0) {\small$q_k^{22}$};
	\node[shift = {(.1,.5)}] at (0,2) {\small$q_k^{13}$};
	\node[shift = {(.1,.5)}] at (0,-2) {\small$q_k^{31}$};
	\node[shift = {(.1,.5)}] at (1,1) {\small$q_k^{12}$};
	\node[shift = {(.1,.5)}] at (-1,-1) {\small$q_k^{32}$};
	\node[shift = {(.1,.5)}] at (1,-1) {\small$q_k^{21}$};
	\node[shift = {(.1,.5)}] at (-1,1) {\small$q_k^{23}$};
	\fill[gray!60,shift={(.3,0))}] (0,0) circle(2.5pt);
	\fill[gray!60,shift={(-.3,0))}] (0,0) circle(2.5pt);
	\fill[gray!60,shift={(.3,0))}] (2,0) circle(2.5pt);
	\fill[gray!60,shift={(-.3,0))}] (2,0) circle(2.5pt);
	\fill[gray!60,shift={(.3,0))}] (0,2) circle(2.5pt);
	\fill[gray!60,shift={(-.3,0))}] (0,2) circle(2.5pt);
	\fill[gray!60,shift={(.3,0))}] (-2,0) circle(2.5pt);
	\fill[gray!60,shift={(-.3,0))}] (-2,0) circle(2.5pt);
	\fill[gray!60,shift={(.3,0))}] (0,-2) circle(2.5pt);
	\fill[gray!60,shift={(-.3,0))}] (0,-2) circle(2.5pt);
	\fill[gray!60,shift={(.3,0))}] (1,1) circle(2.5pt);
	\fill[gray!60,shift={(-.3,0))}] (1,1) circle(2.5pt);
	\fill[gray!60,shift={(.3,0))}] (1,-1) circle(2.5pt);
	\fill[gray!60,shift={(-.3,0))}] (1,-1) circle(2.5pt);
	\fill[gray!60,shift={(.3,0))}] (-1,1) circle(2.5pt);
	\fill[gray!60,shift={(-.3,0))}] (-1,1) circle(2.5pt);
	\fill[gray!60,shift={(.3,0))}] (-1,-1) circle(2.5pt);
	\fill[gray!60,shift={(-.3,0))}] (-1,-1) circle(2.5pt);
\end{tikzpicture}
        \caption{crossing}
    \end{subfigure}
    \caption{Depicted here are local pictures of $\LA_f^m$ near critical points of $f$ and crossings of $\LA$.  At each local maximum of $f$, there is a right-handed half twist. At each local minimum of $f$, there is a left-handed half twist. At each crossing of $\LA$, there is a grid of crossings of $\LA_f^m$. The gray dots label positive quadrants at the crossings.}
\label{half twists}
\end{figure}
\begin{figure}[t]
\centering
\picw{.6}{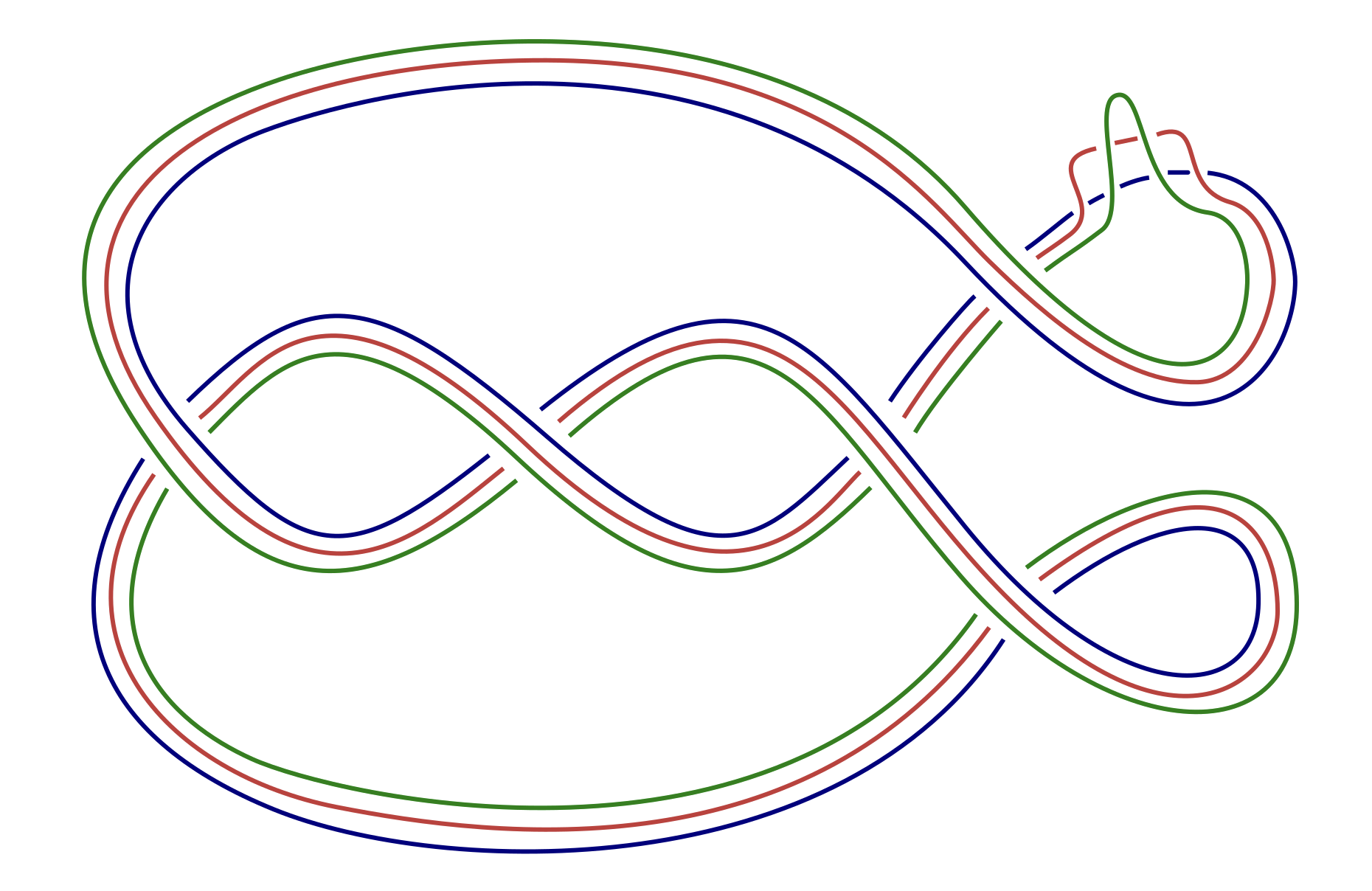}
\caption{Lagrangian $3$-copy of the trefoil}
\label{3-copy trefoil}
\end{figure}

For each $m > 0$, define the \tit{Lagrangian $m$-copy} $\Lambda^m_f$ as follows. Translate $\LA$ in the Reeb ($+z$) direction by small amounts to obtain $m$ parallel copies of $\LA$, labelled by $\LA_1,\cdots,\LA_m$ from high to low. Pick a small neighborhood of $\LA\inn \R^3$ contactomorphic to a neighborhood of the zero section of $J^1(\LA)$ such that its translates are disjoint. Now for each $1\leq \alpha\leq m$, replace $\LA_\alpha$ by the $1$-jet $j^1((m-\alpha)\ep f)$ in its neighborhood, where $\ep > 0$ is small. Away from critical points of $f$, the Lagrangian projection of the union of these $1$-jets looks like $m$ parallel copies of $\LA_{xy}$. At each critical point, the projection is $m$-to-$1$, with $\LA_\alpha$ sitting on top of $\LA_{\alpha+1}$. Now use small Legendrian isotopies to perturb $\LA_\alpha$'s near critical points of $f$ so that at each local maximum/minimum, the projection appears as a right-/left-handed half twist (see \Cref{half twists}). Let $\LA_f^m$ be the disjoint union of these perturbed copies, which we still denote by $\LA_1,\cdots,\LA_m$ by abuse of notation. See \Cref{3-copy trefoil} for the Lagrangian $3$-copy of the trefoil, where $f$ has exactly two critical points and they are adjacent.

The Lagrangian $m$-copy $\LA_f^m$ has two kinds of Reeb chords:
\begin{outline}
    \1 \tit{long} Reeb chords: $a^{\alpha\beta}_j$, for each Reeb chord $a_j$ of $\LA$ and $1\leq \alpha,\beta\leq m$; 
    \1 \tit{short} Reeb chords: 
        \2 $b^{\alpha\beta}_k$, for each local maximum $b_k$ of $f$ and $1\leq \alpha < \beta\leq m$;
        \2 $c^{\alpha\beta}_k$, for each local minimum $c_k$ of $f$ and $1\leq \alpha < \beta\leq m$;
\end{outline}
where the superscript ${}^{\alpha\beta}$ means the Reeb chord starts on $\LA_\beta$ and ends on $\LA_\alpha$. See \Cref{half twists} for illustrations.

For each marked point $\star_i$ on $\LA$, we place a marked point $\star_i^\alpha$ on $\LA_\alpha$ at the image of $\star_i$ under the (perturbed) $1$-jet defining $\LA_\alpha$. Do the same for the base points $\bullet_i$'s. Equip $\LA_f^m$ with the Maslov potential $\mu^m$ given by $\mu^m(\bullet_i^\alpha)\defeq \mu(\bullet_i)$. Note that $\mu^m$ is indeed a Maslov potential because the Lagrangian projections of $\LA_\al$'s are parallel to $\LA_{xy}$ away from the critical points.

We are now ready to describe the generators of the $m$-copy Chekanov-Eliashberg algebra:
\begin{outline}
    \1 $q^{\alpha\beta}_j$, for each long Reeb chord $a^{\alpha\beta}_j$;
    \1 $x^{\alpha\beta}_k$ (resp. $y^{\alpha\beta}_k$), for each short Reeb chord $b^{\alpha\beta}_k$ (resp. $c^{\alpha\beta}_k$);
    \1 $t^{\alpha}_i$ and $(t^{\alpha}_i)^{-1}$, for each marked point $\star^{\alpha}_i$\,.
\end{outline}
The $\Z$-gradings of these generators are specified by $\mu^m$ and \eqref{grading def}. We will again consider the induced $\Z/2\rot(\LA)$-grading, where we have
\[|q^{\alpha\beta}_j| = |q_j|,\ |x^{\alpha\beta}_k| = 0,\ |y^{\alpha\beta}_k| = -1,\ |t^\alpha_i| = 0\,.\]
Note that $|t^\alpha_i|= |t_i| = 0$ by choice of $\mu$. 

Write $\mbfQ_m = \mbfQ(\LA_f^m) = \{q^{\alpha\beta}_j,x^{\alpha\beta}_k,y^{\alpha\beta}_k\}$, $\mbfT_m = \mbfT(\LA_f^m) = \{(t^\alpha_i)^{\pm 1}\}$, and \[\mbfS_m^0 = \mbfS^0(\LA_f^m) = \{q^{\alpha\beta}_j,x^{\alpha\beta}_k,y^{\alpha\beta}_k,(t^\alpha_i)^{\pm 1}\}\,.\] Let $\mbfQ_m^\mrm{long} = \{q^{\alpha\beta}_j\}$, $\mbfQ_m^\mrm{short} = \{x^{\alpha\beta}_k,y^{\alpha\beta}_k\}$, $\mbfQ_m^+ = \mbfQ_m^\mrm{short}\cup\{q^{\alpha\beta}_j\mid \alpha < \beta\}$, $\mbfQ_m^- = \{q^{\alpha\beta}_j\mid \alpha > \beta\}$. 
\begin{defn}
\label{m-copy link grading def}
The \tit{$m$-copy link grading} of the Chekanov-Eliashberg generators of $\LA_f^m$ is the pair of functions \[e_-^m,e_+^m:\mbfS_m^0\to \{1,\cdots, m\}\] given by \[\begin{cases}
    e_-^m(q^{\alpha\beta}_j) = e_-^m(x^{\alpha\beta}_k) = e_-^m(y^{\alpha\beta}_k) = e_-^m(t^\alpha_i) = e_+^m(t^\alpha_i) = \alpha\\
    e_+^m(q^{\alpha\beta}_j) = e_+^m(x^{\alpha\beta}_k) = e_+^m(y^{\alpha\beta}_k) = \beta
\end{cases}\hspace{-1ex}.\]  
\end{defn}
\begin{rmk}
\label{link grading diff rmk}
    Note the difference between the $m$-copy link grading and the internal link grading of $\LA_f^m$, as in \Cref{internal link grading def}. However, if $\LA$ is a Legendrian knot with single base point, then the two link gradings agree.
\end{rmk}

\begin{defn}
The \tit{composable $m$-copy Chekanov-Eliashberg algebra} $\mcalA_m$ of $\LA$ with Morse function $f$ is defined as the tensor algebra of the $\Z/2\rot(\LA)$-graded $\k^m$-bimodule $V(\mbfS_m^0;e_-^m,e_+^m)$, modulo the relations \[t^\alpha_i(t^\alpha_i)^{-1} = (t^\alpha_i)^{-1}t^\alpha_i=e_\alpha\,.\] The differential $\del_m$ is defined just like the Chekanov-Eliashberg differential $\del$ of $\LA_f^m$, with one modification: for each $s\in \mbfS_m^0$ with $e_-^m(s) = e_+^m(s) = \alpha$,  replace $1$ with $e_\alpha$ in $\del s$. Also set $\del_m e_\alpha = 0$ for each $1\leq \alpha \leq m$. 
\end{defn}
\begin{rmk}
The differential $\del_m$ is well-defined on $\mcalA_m$ because any word $s_1\cdots s_k$ in $\del s$ is composable in the sense of \Cref{composability rmk} and satisfies $e_-(s_1) = e_-(s)$ and $e_+(s_k) = e_+(s)$. Also see the discussion in \cite[Section 3.2]{NRSSZ}. 
\end{rmk}

It turns out the $m$-copy differential $\del_m$ is determined by the $1$-copy differential $\del$, the positions of the critical points of $f$, and the internal link grading $e_-,e_+$ of $\LA$ (see \Cref{internal link grading def}).
 \begin{prop}[{\cite[Proposition 4.14]{NRSSZ}}]
 \label{matrix formula prop} Assume there are as many marked points as local maxima, and each $\star_k$ appears immediately before $b_k$, which appears immediately before $c_k$, according to the orientation of $\LA$. Here, ``immediately'' means there are no crossing points on the oriented segment connecting the two points. Define the following elements of $\mcalA_m$: \[A_j = \sum_{\alpha,\beta} q^{\alpha\beta}_j\,,\ X_k = \sum_{\alpha < \beta} x_k^{\alpha\beta}\,,\ Y_k = \sum_{\alpha < \beta} y_k^{\alpha\beta}\,,\ T_k = \sum_\alpha t_k^\alpha\,.\]
     Then we have
     \begin{equation}
         \label{matrix formulas}
         \begin{aligned}
             &\del_mA_j = \Phi(\del q_j)+Y_{e_-(q_j)}A_j - (-1)^{|q_j|}A_jY_{e_+(q_j)}\\
             &\del_mX_k = T_k^{-1}Y_{k-1}T_k(1+X_k)-(1+X_k)Y_{k}\\
             &\del_mY_k = Y_k^2
         \end{aligned}
     \end{equation}
     where $\Phi:\mcalA\to \mcalA_m$ is the $\k$-algebra map given by \[\Phi(q_j) = A_j\,,\ \Phi(t_k) = T_k(1+X_k)\,,\ \Phi(t_k^{-1}) = (1+X_k)^{-1}T_k^{-1}\,.\]
     Note that to extract $\del_ms$ for any $s\in \mbfS_m^0$, one multiplies the appropriate equation by $e_{e_-^m(s)}$ on the left and $e_{e_+^m(s)}$ on the right.
 \end{prop}

We recall the consistency results needed to construct $\Aug_+$. Observe that if we delete $k$ copies from $\LA_f^m$, the resulting link looks exactly like $\LA_f^{m-k}$. This observation defines algebra maps between $\mcalA_m$'s, as follows. Given any $m'$-element subset $I\inn [m] = \{1,\cdots,m\}$, let $J_I$ be the two-sided ideal generated by \[\{s\in \mbfS_m^0\mid e_-^m(s)\notin I \text{ or } e_+^m(s)\notin I\}\,.\] Since $\del_m$ is $\k^m$-linear, it preserves $J_I$, so we get a quotient dga $\mcalA_m^I\defeq \mcalA_m/J_I$. The unique order-preserving bijection $\iota:[m']\to I\inn [m]$ defines an injection $h_I:\mbfS_{m'}^0\to\mbfS_m^0$ given by \begin{equation}
    \label{hI}
    h_I(q^{\alpha\beta}_j) = q^{\iota(\alpha)\iota(\beta)}_j\,,\ h_I(x^{\alpha\beta}_k) = x^{\iota(\alpha)\iota(\beta)}_k\,,\ h_I(y^{\alpha\beta}_j) = y^{\iota(\alpha)\iota(\beta)}_k\,,\ h_I((t^{\alpha}_i)^{\pm1}) = (t^{\iota(\alpha)}_i)^{\pm1}\,.
\end{equation} Hence, we get an inclusion of $\k^{m'}$-algebras $h_I:\mcalA_{m'}\to \mcalA_m$, where $\mcalA_m$ is viewed as an $\k^{m'}$-algebra via the map $\k^{m'}\to \k^m$ induced by $\iota$. This inclusion does not respect the differentials, but we have:

\begin{lemma}[{\cite[Proposition 4.2]{NRSSZ}}]
\label{m-copy CE consistency}
    The composite $\mcal{A}_{m'}\xrar{h_I} \mcal{A}_m\to \mcal{A}_m^I$ is an isomorphism of dg $\k^{m'}$-algebras. 
\end{lemma}
The above lemma allows us to assemble $m$ augmentations of $\mcalA$ into one for $\mcalA_m$:
\begin{lemma}[{\cite[Lemma 3.15]{NRSSZ}}]
Given any $m$-tuple $(\ep_1,\cdots,\ep_m)$ of augmentaions of $\mcalA$, the map $\ep: \mcalA_m\to\k^m$ given by \[\ep(s) =\begin{cases}
    \ep_\alpha(s)e_\al & e_-^m(s) = e_+^m(s) = \al\\
    0 & e_-^m(s)\neq e_+^m(s)
\end{cases}\] for any $s\in\mbfS_m^0$, is an augmentation.
\end{lemma}
We say $\ep$ is \tit{pure}, in the sense that it vanishes on mixed generators. By abuse of notation, we write $\ep = (\ep_1,\cdots,\ep_m)$ for the augmentation of $\mcalA_m$ defined above, which is called the \tit{diagonal augmentation} in \cite{NRSSZ}.

Define a $\k^m$-algebra automorphism $\phi_\ep: \mcalA_m\to \mcalA_m$ by 
 \[\phi_\ep(s) \defeq 
 \begin{cases}
     s+\ep(s) & s\in\mbfQ_m\\
     s & s\in\mbfT_m
 \end{cases}
 \]for generators $s\in \mbfS_m^0$. Define the \tit{twisted $m$-copy differential} \[\del_m^\ep \defeq \phi_\ep\circ  \del_m\circ \phi_\ep^{-1}\] which clearly squares to zero. Moreover, it descends to the quotient \[\mcalA_m^\ep \defeq \mcalA_m/\agl{t-\ep(t)\mid t\in\mbfT_m}\] since $\del_m^\ep t = \del_m t = 0$ for any $t\in\mbfT_m$. 
 
\begin{lemma}
\label{twisted m-copy CE no constant}
    For any $q\in\mbfQ_m$, the differential $\del_m^\ep q$ has no constant terms (i.e., elements of $\k^m$).
\end{lemma}
\begin{proof}
    The constant term in $\del_m^\ep q$ is precisely $\ep\circ\del_m (q)$, which is zero since $\ep:\mcalA_m\to\k^m$ is an augmentation.
\end{proof}

The dgas $\{\mcal{A}_m^\ep\}$ with twisted differentials inherit a similar structure: given $\ep = (\ep_1,\cdots,\ep_m)$, we can augment $\mcal{A}_{m'}$ by $\ep_I = \{\ep_\alpha\mid \alpha\in I\}$, then we have:

\begin{lemma}
\label{augmented m-copy CE consistency}
 The above isomorphism $\mcal{A}_{m'}\iso \mcal{A}_m^I$ induces an isomorphism $\mcal{A}_{m'}^{\ep_I}\iso (\mcal{A}_m^I)^{\ep_I}$. Moreover, we can identify $(\mcal{A}_m^I)^{\ep_I} = (\mcal{A}_m^\ep)^I$. 
\end{lemma}
\begin{proof}
    The first statement is clear from the definition of the twisted differential. For the second, note that the two dgas $(\mcal{A}_m^I)^{\ep_I}$ and $(\mcal{A}_m^\ep)^I$,  have the same underlying algebra. Since $\ep:\mcal{A}_m\to \k^m$ is pure, they also have the same differential.
\end{proof}

\subsection{Dualizing Derivations}
\label{dualizing derivations section}
In \cite{NRSSZ}, it was claimed without proof that dualizing using the pairing \eqref{dualizing pairing} takes derivations to coderivations. Since this claim will be crucial to our constructions, we include a proof in this subsection. 

Let $V = V(\mbfS;e_-,e_+)$ as in \Cref{km bimodule}. Make its $\k$-linear dual $V^* = \Hom_\k(V,\k)$ into a graded $\k^m$-bimodule by setting $|s^*| = |s|$ and \[e_\alpha s^*e_\beta = \de_{\alpha,e_+(s)}\de_{\beta,e_-(s)}s^*\] for $s\in\mbfS$. (Note the grading difference from $V^\dagger$, where $|s^\dagger| = -|s|$, as defined just before \Cref{linear dual bimodule def}.) 

For this subsection, write $\otimes = \otimes_{\k^m}$, $TV = T_{\k^m}(V)$, and $TV^* = T_{\k^m}(V^*)$. Define a (non-homogeneous) $\k$-linear perfect pairing
$\agl{\cdot,\cdot}:TV^*\otimes TV\to \k$ by
\begin{equation}
    \label{dualizing pairing}
    \agl{s_k^*\cdots s_1^*,s_1\cdots s_k} = (-1)^{\sum_{i<j}|s_i||s_j|}
\end{equation} and otherwise zero. Note that $\agl{e_\alpha\cdot{}, \cdot} = \agl{\cdot,{}\cdot e_\alpha}$ for any $e_\alpha$. In particular, for any homogeneous $\k^m$-linear map $d:TV\to TV$, its $\k$-linear dual $d^*: TV^*\to TV^*$, defined via \[\agl{\cdot,d(\cdot)} = \agl{d^*(\cdot),\cdot}\] is in fact $\k^m$-linear and homogeneous of degree $-|d|$. 
\begin{lemma}
\label{dualizing derivations}
Under the above pairing, the dual of a ($\k^m$-linear) derivation on $TV$ is a ($\k^m$-linear) coderivation on the coalgebra $TV^*$. (The comultiplication of the tensor coalgebra is specified in \eqref{tensor comultiplication}.)
\end{lemma}
\begin{proof}
Define an auxiliary perfect pairing $(TV^*)^{\otimes 2}\otimes TV^{\otimes 2}\to\k$ as the composition
\[
\tikzcdset{every label/.append style = {scale=0.7,yshift=0.6ex}}
\begin{tikzcd}[column sep = 60]
   TV^*\otimes TV^*\otimes TV\otimes TV\rar{1\otimes \agl{\cdot,\cdot}\otimes 1} & TV^*\otimes TV\rar{\agl{\cdot,\cdot}} & \k
\end{tikzcd}
\]which we will also denote by $\agl{\cdot,\cdot}$. To prove the lemma, we first show \eqref{paring coproduct and product} and \eqref{dualizing tensor product of maps} below.

Under this pairing, we have 
\begin{align*}
    &\quad\, \agl{\DE(s_n^*\cdots s_1^*),s_1\cdots s_k\otimes s_{k+1}\cdots s_n}\\
    &= \agl{s_n^*\cdots s_{k+1}^*\otimes s_k^*\cdots s_1^*, s_1\cdots s_k\otimes s_{k+1}\cdots s_n}= (-1)^{\sum_{1\leq i < j\leq k}|s_i||s_j|+\sum_{k+1\leq i < j\leq n}|s_i||s_j|}\,.
\end{align*}
On the other hand, the multiplication map $m:TV\otimes TV\to TV$ satisfies
\[
    \agl{s_n^*\cdots s_1^*, m(s_1\cdots s_k\otimes s_{k+1}\cdots s_n)} = \agl{s_n^*\cdots s_1^*, s_1\cdots s_n} = (-1)^{\sum_{1\leq i < j\leq n}|s_i||s_j|}\,.
\]
Since all nonzero pairings must be of the above form, we get
\begin{equation}
\label{paring coproduct and product}
\agl{\DE(\cdot),s_1\cdots s_k\otimes s_{k+1}\cdots s_n} =  (-1)^{(\sum_{i=1}^k |s_i|)(\sum_{i=k
 +1}^n |s_i|)}\agl{\cdot, m(s_1\cdots s_k\otimes s_{k+1}\cdots s_n)}\,.
 \end{equation}

 For homogeneous $f,g: TV\to TV$, $a,b\in TV$, $c,d\in TV^*$, we have
\[\begin{cases}
    \agl{c\otimes d, f\otimes g(a\otimes b)}= (-1)^{|g||a|}\agl{c\otimes d, f(a)\otimes g(b)}= (-1)^{|g||a|}\agl{d,f(a)}\agl{c,g(b)}\\\agl{g^*\otimes f^*(c\otimes d), a\otimes b}
    = (-1)^{|f^*||c|}\agl{f^*(d),a}\agl{g^*(c),b}= (-1)^{|f||c|}\agl{d,f(a)}\agl{c,g(b)}\end{cases}\hspace{-1ex}.\]
For these pairings to be nonzero, we must have $|d| = |f(a)|$ and $|c| = |g(b)|$. Thus, we get
$(-1)^{|f||c|} = (-1)^{|f||g| + |f||b|}$,
giving us
 \begin{equation}
 \label{dualizing tensor product of maps}
     \agl{\cdot, f\otimes g(a\otimes b)} = (-1)^{|g||a|+|f||b| + |f||g|}\agl{g^*\otimes f^*(\cdot),a\otimes b}\,.
 \end{equation}

 Now given any derivation $d$ on $TV$, we have
\begin{align*}&\quad\,\agl{\DE d^*(\cdot),s_1\cdots s_k\otimes s_{k+1}\cdots s_n} \\ &=(-1)^{(\sum_{i=1}^k |s_i|)(\sum_{i=k
 +1}^n |s_i|)}\agl{\cdot,dm(s_1\cdots s_k\otimes s_{k+1}\cdots s_n)}\\ &= (-1)^{(\sum_{i=1}^k |s_i|)(\sum_{i=k
 +1}^n |s_i|)}\agl{\cdot,m(d\otimes 1+1\otimes d)(s_1\cdots s_k\otimes s_{k+1}\cdots s_n)}\\ &=(-1)^{\sum_{i=k
 +1}^n |s_i|}\agl{\DE(\cdot),d\otimes 1(s_1\cdots s_k\otimes s_{k+1}\cdots s_n)}\\ &\quad+(-1)^{\sum_{i=1}^k |s_i|}\agl{\DE(\cdot),1\otimes d(s_1\cdots s_k\otimes s_{k+1}\cdots s_n)}\\ 
 &=\agl{(1\otimes d^*)\DE(\cdot),s_1\cdots s_k\otimes s_{k+1}\cdots s_n}+\agl{(d^*\otimes 1)\DE(\cdot),s_1\cdots s_k\otimes s_{k+1}\cdots s_n}\end{align*}
where we used \eqref{paring coproduct and product} for the first and third equalities and \eqref{dualizing tensor product of maps} for the last equality. Since the pairing is perfect, we have shown $\DE d^* = (d^*\otimes 1+1\otimes d^*)\DE$, and so $d^*$ is a coderivation on $TV^*$.
\end{proof}
\begin{rmk}
\label{completed dualizing}
    The above lemma also holds if $TV$ is replaced by $TV$ completed with respect to any subset of the generators, since the pairing is still well-defined and the same proof works. We will need this generalization in \Cref{twisted m-copy LSFT}.
\end{rmk}
\begin{coro}
\label{A infinity km algebra}
    Let $\br{T}V = \bigoplus_{k >0}V^{\otimes k}\inn TV$. If $d$ is a degree $-1$ $\k^m$-linear derivation on $TV$ such that $d^2 = 0$ and $d(\br{T}V)\inn \br{T}V$, then \[V^\vee\defeq V^*[-1]\]
    equipped with the codifferntial $d^*$ on $BV^\vee = TV^*$ in an $A_\infty$ $\k^m$-algebra.
\end{coro}
\begin{proof}
    By \Cref{dualizing derivations}, $d^*$ is a degree $1$ $\k^m$-linear coderivation on $BV^\vee = TV^*$. Since $d^2 = 0$ and $d(\br{T}V)\inn \br{T}V$, we get $(d^*)^2 = 0$ and $d^*(1) = 0$.
\end{proof}
\subsection{\texorpdfstring{$\Aug_+$}{Aug+} and the Short Chord Category \texorpdfstring{$\mcalC$}{C}}
\label{aug+ and C}
In this subsection, we recall the construction of $\Aug_+$ from \cite{NRSSZ}, phrased in terms of the composable $m$-copy dgas. Then we define the \tit{short chord category} $\mcalC$ as a wide quotient of $\Aug_+$ and exhibit the diagonal bimodule of $\Aug_+$ as the mapping cone of a morphism of $\Aug_+$-bimodules. 

The underlying algebra of $\mcalA_m^\ep$ is $T_{\k^m}V(\mbfQ_m;e_-^m,e_+^m)$. Hence, by \Cref{twisted m-copy CE no constant} and \Cref{A infinity km algebra}, dualizing the differential $\del^\ep_m$ on $\mcalA_m^\ep$ gives us an $A_\infty$ $\k^m$-algebra structure on
\[Q_m\defeq  V(\mbfQ_m;e_-^m,e_+^m)^\vee\] with codifferential $b^\ep_m\defeq  (\del^\ep_m)^*$. Let $I^\geq_m\inn \mcal{A}_m^{\ep}$ be the two-sided ideal generated by $\{q^{\alpha\beta}_j\mid \alpha \geq \beta\}$. Since $\del_m^\ep$ is $\k^m$-linear, it preserves $I^\geq_m$ and so descends to the quotient $\mcalA_m^\ep/I^\geq_m$, whose underlying algebra is $T_{\k^m}V(\mbfQ_m^+;e_-^m,e_+^m)$. Therefore, \[Q_m^+\defeq V(\mbfQ_m^+;e_-^m,e_+^m)^\vee\] is an $A_\infty$ $\k^m$-subalgebra of $Q_m$, meaning the obvious inclusion is a strict morphism.

The injections $h_I:\mbfS_n\to\mbfS_m$ defined in \eqref{hI} induce linear maps $h_I: Q_n\to Q_m$ by letting $h_I(s^\vee) = h_I(s)^\vee$. For any $1\leq \alpha,\beta\leq m$, write $h_{\alpha\beta} \defeq h_{\{\alpha,\beta\}}$ and let $Q_{\alpha\beta}\inn Q_m$ be the $\k$-submodule spanned by $\{(q_j^{\alpha\beta})^\vee\mid 1\leq j\leq n\}$.
\begin{deflem}
The \tit{positive augmentation category} of $\LA$ with Morse function $f$ is the $\Z/2\rot(\LA)$-graded $A_\infty$ category $(\Aug_+,b^+)$, defined as follows. The objects of $\Aug_+$ are the augmentations of $\LA$. For any two augmentations $\ep,\ep'$, set \[\Aug_+(\ep,\ep')\defeq Q^+\defeq Q^+_2\] which is spanned by $\{(q^{12}_j)^\vee,(x^{12}_k)^\vee, (y^{12}_k)^\vee\mid 1\leq j\leq n, 1\leq k\leq \ell\}$. Following \cite{NRSSZ}, we write $q_j^+, x_k^+,y_k^+$, when viewed as basis elements of $\Aug_+(\ep,\ep')$. Their degrees are given by
\[|q_j^+| = |q_j|+1\,,\quad |x_k^+| = 1\,, \quad |y_k^+| = 0\,.\]
For any $m\geq 1$ and augmentations $\ep_1,\cdots,\ep_{m+1}$, define $(b^+_m)^{\ep_1,\cdots,\ep_{m+1}}$ to be the composition 
\[
\begin{tikzcd}[column sep = 50,row sep = 45]
Q^+[1]^{\otimes m}\dar{h_{m,m+1}\otimes\cdots\otimes h_{12}}\\ Q_{m,m+1}[1]\otimes\cdots\otimes Q_{12}[1]\rar{(b^{\ep}_{m+1})_m} &Q_{1,m+1}[1]\arrow{r}{h_{1,m+1}^{-1}} & Q^+[1]\end{tikzcd}
\]where $\ep = (\ep_1,\cdots,\ep_{m+1})$.
\end{deflem}

\begin{proof}
Given any $m\geq 1$ and tuple $\ep = (\ep_1,\cdots,\ep_{m+1})$, we need to check the corresponding $A_\infty$ relation \eqref{A_infty relations}. To do this, note that for any $m'\leq m$ and $1\leq \alpha_1 <\cdots < \alpha_{m'+1} \leq m+1$, the analogous composition
\[
\begin{tikzcd}[column sep = 50, row sep = 45]
Q^+[1]^{\otimes k}\dar{h_{\alpha_{m'}\alpha_{m'+1}}\otimes\cdots\otimes h_{\alpha_1\alpha_2}}\\ Q_{\alpha_{m'}\alpha_{m'+1}}[1]\otimes\cdots\otimes Q_{\alpha_1\alpha_2}[1]\rar{(b^{\ep}_{m+1})_{m'}} &Q_{\alpha_1\alpha_{m'+1}}[1]\arrow{r}{h_{\alpha_1\alpha_{m'+1}}^{-1}} & Q^+[1]\end{tikzcd}
\] agrees with $(b^+_{m'})^{\ep_{\alpha_1},\cdots,\ep_{\alpha_{m'+1}}}$, by \Cref{augmented m-copy CE consistency}. Hence, the desired relation is equivalent to the $m$-input relation for the $A_\infty$ algebra $(Q^+_{m+1},b^\ep_{m+1})$, since the first and last maps in the above compositions are isomorphisms.
\end{proof}

Explicitly, for any $m\geq 1$, augmentations $\ep_1,\cdots,\ep_{m+1}$, and $s_1,\cdots,s_m\in\mbfQ\cup\{x_k,y_k\}_{k=1}^\ell$, we have
\begin{equation}
\label{b+ formula}
\begin{aligned}
&\quad\,s^{-1}(b^+_m)^{\ep_1,\cdots,\ep_{m+1}}[s^+_m|\cdots| s^+_1]\\ 
&= (-1)^{\sum_{\alpha < \beta}|s_\al||s_\beta|}\sum_{q\in \mbfQ\cup\{x_k,y_k\}_{k=1}^\ell}\Coeff_{s_1^{12}\cdots s_k^{m,m+1}}(\del^\ep_{m+1}q^{1,m+1})\cdot q^+
\end{aligned}
\end{equation}
where $\ep = (\ep_1,\cdots,\ep_{m+1})$.
\begin{prop}[{\cite[Proposition 3.28]{NRSSZ}}]
    The $A_\infty$ category $\Aug_+$ is strictly unital, and the unit $e:\k_{\ob\Aug_+}\to\Aug_+$ is given by \[e_\ep = -\sum_{k=1}^\ell y_k^+\in\Aug_+(\ep,\ep)\] for any $\ep$.
\end{prop}
\begin{lemma}
    The $\k$-subquiver $\mcalQ_+^\mrm{long}\inn \Aug_+$ given by \[\mcalQ_+^\mrm{long}(\ep,\ep')\defeq \k\{q_j^+\mid 1\leq j\leq n\}\] is a wide ideal (see \Cref{wide ideal}).
\end{lemma}
\begin{proof}
    Since $\mcalA_m$ is filtered by the length of Reeb chords, the twisted differential $\del^\ep_m$ preserves the $\k^m$-subalegbra of $\mcalA_m^\ep$ generated by the short Reeb chords $x_k^{\alpha\beta},y_k^{\alpha\beta}$. Hence, the coefficient of $s_1^{12}\cdots s_m^{m,m+1}$ in $d^\ep_{m+1}q^{1,m+1}$ is zero if $q^{1,m+1}$ is a short Reeb chord and $s_\alpha^{\alpha,\alpha+1}$ is a long Reeb chord for some $1\leq \alpha\leq m$. Now by the formula \eqref{b+ formula} of $b^+_m$, we see that $b^+_m$ sends $\Aug_+[1]^{\otimes i}\otimes \mcalQ_+^\mrm{long}[1]\otimes \Aug_+[1]^{\otimes j}$ to $\mcalQ_+^\mrm{long}[1]$, for any $i+j = m-1$. 
\end{proof}

\begin{defn}
\label{short chord category}
    The \tit{short chord category} of $\LA$ with Morse function $f$ is the wide quotient\[\mcalC\defeq \Aug_+/\mcalQ_+^\mrm{long}\] as in \Cref{wide ideal correspondence}. Denote its codifferential by $b^\mcalC$ and the quotient map by
    \[\pi:\Aug_+\to\mcalC\] which is a strict $A_\infty$ functor.
\end{defn}
 For any $\ep,\ep'$, we can identify 
 \[\mcalC(\ep,\ep')= \k\{x^+_k,y^+_k\mid1\leq k\leq \ell\}\] as graded $\k$-modules. Hence, we have a splitting of $\k$-quivers
 \[\Aug_+ = \mcalC\oplus \mcalQ_+^\mrm{long}\,.\]
With respect to this splitting, the codifferential of the diagonal bimodule of $\Aug_+$ can be written as \[b^+ = \begin{bmatrix}
    b^{\mcalQ_+^\mrm{long}} &0 \\
    * & b^{\pi^*\mcalC}
\end{bmatrix}\]where $\pi^*\mcalC$ is the pullback bimodule, as in \Cref{pullback def}. Therefore, we have:
\begin{lemma}
\label{cone for long}
The diagonal bimodule $\Aug_+$ is the mapping cone of a morphism of $\Aug_+$-bimodules $\pi^*\mcalC[-1]\to \mcalQ_+^\mrm{long}$.
\end{lemma}

Under the assumptions of \Cref{matrix formula prop}, we can explicitly write out the structure maps of the short chord category $\mcalC$. Let $\{\beta^\mcalC_m\}$ denote the un-suspended components of $b^\mcalC$, then we have
\begin{equation}
\label{bC formula}
    \begin{cases}\beta_1^\mcalC(x_k^+) = 0,\ (\beta_1^\mcalC)^{\ep_1,\ep_2}(y_k^+) = \ep_1(t_{k+1})^{-1}\ep_2(t_{k+1})x_{k+1}^+ - x_k^+\\ \beta_2^\mcalC(x_k^+,x_k^+) = 0,\ (\beta_2^\mcalC)^{\ep_1,\ep_2,\ep_3}(x_k^+, y_{k-1}^+) = -\ep_1(t_k)^{-1}\ep_2(t_k)x_k^+\\\beta_2^\mcalC(y_{k}^+,x_k^+) = -x_k^+,\ \beta_2^\mcalC(y_k^+,y_k^+) = -y_k^+
\end{cases}
\end{equation}
and $\beta_m^\mcalC = 0$ for $m\geq 3$. Thus, $\mcalC$ is in fact a dg category.
\begin{prop}
\label{endo computes cohomology}
    For any augmentation $\ep$, the complex $(\mcalC(\ep,\ep), \beta^\mcalC_1)$ is the Morse cochain complex of $f$.
\end{prop}
\begin{proof}
The basis elements of $\mcalC(\ep,\ep)$ correspond to critical points of $f$, and their degrees agree with the indices of the corresponding critical points. By \eqref{bC formula}, we have \[\beta_1^\mcalC(x_k^+) = 0\,,\quad \beta_1^\mcalC(y_k^+) = x_{k+1}^+ - x_k^+\] on $\mcalC(\ep,\ep)$, which agrees with the Morse differential of $f$.
\end{proof}
If $\LA$ is a Legendrian knot with a single base point, then by Leverson's theorem \cite{Lev16}, $\ep(t) = -1$ for every augmentation $\ep$. So in this case, the structure maps simplify to \begin{equation}
    \label{C for knot}
    \beta_2^\mcalC(x^+,x^+) = 0,\ \beta_2^\mcalC(x^+, y^+) = \beta_2^\mcalC(y^+,x^+) = -x^+,\ \beta_2^\mcalC(y^+,y^+) = -y^+
\end{equation} 
and $\beta_m^\mcalC = 0$ for $m\neq 2$. 

\subsection{Negative Augmentation Bimodule}
\label{negative aug bimod}
In \cite[Proposition 3.23]{NRSSZ}, it is claimed that the dg $\k$-quiver $\Hom_-$ extends to an $\Aug_+$-bimodule, but the full construction is omitted. In this subsection, we complete this construction by defining the \tit{negative augmentation bimodule} $\mcalQ_-$, which turns out to be strictly isomorphic to $\mcalQ_+^\mrm{long}$. As a result, we prove \Cref{intro A infty EES}.

Since $Q_m^+\inn Q_m$ is an $A_\infty$ subalgebra, we can regard $Q_m$ as a $Q_m^+$-bimodule, which has a sub-bimodule spanned by $\mbfQ_m^+\cup \{q_j^{\alpha\alpha}\mid j,\alpha\}$. The quotient bimodule can be identified with \[Q_m^-\defeq V(\mbfQ_m^-;e_-^m,e_+^m)^\vee\] as graded $\k^m$-bimodules. Denote its codifferential by $b^{Q_m^-}_\ep$.
We assemble these into an $\Aug_+$-bimodule:
\begin{deflem}
\label{Q- def}
The \tit{negative augmentation bimodule} of $\LA$ with Morse function $f$ is the $\Aug_+$-bimodule $(\mcalQ_-, b^{\mcalQ_-})$, defined as follows. For any two augmentations $\ep,\ep'$, set \[\mcalQ_-(\ep,\ep')\defeq Q^-\defeq Q^-_2\] which is spanned by $(q^{21}_1)^\vee, \cdots, (q^{21}_n)^\vee$. Following \cite{NRSSZ}, we write $q_1^-, \cdots, q_n^-$, when viewed as basis elements of $\mcalQ_-(\ep,\ep')$. Their degrees are given by \[|q_j^-| = |q_j| + 1\,.\]
For any $i,j\geq 0$ and augmentations $\ep_1,\cdots,\ep_{j+1},\ep_{j+2},\cdots,\ep_{m+1}$ where $m=i+j+1$, define the component $(b^{\mcalQ_-}_{i|j})^{\ep_1,\cdots,\ep_{j+1}}_{\ep_{j+2},\cdots,\ep_{m+1}}$ to be the composition
\[
\begin{tikzcd}[column sep = 60,row sep = 40,ampersand replacement=\&]
Q^+[1]^{\otimes i}\otimes Q^-[1]\otimes Q^+[1]^{\otimes j} \dar{h_{i,i+1}\otimes\cdots \otimes h_{12}\otimes h_{m+1,1}\otimes h_{m,m+1}\otimes\cdots\otimes h_{i+2,i+3}}\\ \begin{multlined}
 Q_{i,i+1} [1]\otimes\cdots\otimes Q_{12}[1]
 \otimes Q_{m
+1,1}[1]\\ \otimes Q_{m,m+1}[1]\otimes\cdots \otimes Q_{i+2,i+3} [1]   
\end{multlined}
\rar{(b^{Q_{m+1}^-}_\ep)_{i|j}} \&Q_{i+2,i+1}[1]\arrow{r}{h_{i+2,i+1}^{-1}} \& Q^-[1]\end{tikzcd}
\]where $\ep=(\ep_{j+2},\cdots,\ep_{m+1},\ep_1,\cdots,\ep_{j+1})$. Note that the two tuples are swapped.
\end{deflem}
\begin{proof}
    Given any $i,j\geq 0$ and augmentations $\ep_1,\cdots,\ep_{j+1},\ep_{j+2},\cdots,\ep_{m+1}$ where $m = i+j+1$, we need to check the $(i,j)$-th relation \eqref{A_infty bimodule relations}, restricted to \[\Aug_+[1](\ep_{j+2},\cdots,\ep_{m+1})\otimes \mcalQ_-[1](\ep_{j+1},\ep_{j+2})
         \otimes \Aug_+[1](\ep_1,\cdots,\ep_{j+1})\,.\] To do this, note that for any $0\leq i'\leq i$, $0\leq j'\leq j$, $1\leq \alpha_1 <\cdots < \alpha_{j'+1} \leq j+1$, and $j+2\leq \al_{j'+2} <\cdots < \al_{m'+1} \leq m+1$, where $m' = i'+j'+1$, the analogous composition
\[
\begin{tikzcd}[column sep = 60,row sep = 40,ampersand replacement=\&]
Q^+[1]^{\otimes i'}\otimes Q^-[1]\otimes Q^+[1]^{\otimes j'} \arrow{d}{h_{\al_{i'}\al_{i'+1}}\otimes\cdots \otimes h_{\al_1\al_2}\otimes h_{\al_{m'+1}\al_1}\otimes h_{\al_{m'}\al_{m'+1}}\otimes\cdots\otimes h_{\al_{i'+2}\al_{i'+3}}}\\ \begin{multlined}
Q_{\al_{i'}\al_{i'+1}} [1]\otimes\cdots\otimes Q_{\al_1\al_2}[1] 
 \otimes Q_{\al_{m'
+1}\al_1}[1]\\ \otimes  Q_{\al_{m'}\al_{m'+1}}[1]\otimes\cdots \otimes Q_{\al_{i'+2}\al_{i'+3}}[1]
\end{multlined}
\rar{(b^{Q_{m+1}^-}_\ep)_{i'|j'}} \&Q_{\al_{i'+2}\al_{i'+1}}[1]\rar{h_{\al_{i'+2}\al_{i'+1}}^{-1}} \& Q^-[1]\end{tikzcd}
\] where $\ep = (\ep_{\al_{j'+2}},\cdots,\ep_{\al_{m'+1}}, \ep_{\al_1},\cdots,\ep_{\al_{j'+1}})$, agrees with $(b^{\mcalQ_-}_{i'|j'})^{\ep_{\al_1},\cdots,\ep_{\al_{j'+1}}}_{\ep_{\al_{j'+2}},\cdots,\ep_{\al_{m'+1}}}$, by \Cref{augmented m-copy CE consistency}. Hence, the desired relation is equivalent to the $(i,j)$-th relation for the $Q^+_{m+1}$-bimodule $(Q^-_{m+1},b^{Q_{m+1}^-}_\ep)$, since the first and last maps in the above compositions are isomorphisms. 
\end{proof}

Explicitly, for any $i,j\geq 0$, augmentations $\ep_1,\cdots,\ep_{j+1},\ep_{j+2},\cdots,\ep_{m+1}$ where $m = i+j+1$, $s_1,\cdots,s_i,s_{i+2},\cdots, s_m\in\mbfQ\cup\{x_k,y_k\}_{k=1}^\ell$, and $s_{m+1}\in\mbfQ$, we have
\begin{equation}
\label{bQ- formula}
\begin{aligned}
   &\quad\, s^{-1}(b^{\mcalQ_-}_{i|j})_{\ep_{j+2},\cdots,\ep_{m+1}}^{\ep_{1},\cdots,\ep_{j+1}}[s_i^+|\cdots| s_1^+| s_{m+1}^-| s_m^+|\cdots | s_{i+2}^+]\\
   &= (-1)^{\sum_{\alpha < \beta}|s_\al||s_\be|}\sum_{q\in \mbfQ}\Coeff_{s_{i+2}^{i+2,i+3}\cdots  s_m^{m,m+1}s_{m+1}^{m+1,1}s_1^{12}\cdots s_i^{i,i+1}}\left(\del^\ep_{m+1}q^{i+2,i+1}\right)\cdot q^- 
\end{aligned}  
\end{equation} where $\ep=(\ep_{j+2},\cdots,\ep_{m+1},\ep_1,\cdots,\ep_{j+1})$.

\begin{rmk}
    The $\Aug_+$-bimodule $\mcalQ_-$ defined above has the same underlying dg $\k$-quiver as the negative augmentation category $\Aug_-$, defined in \cite{BC14}\cite[Definition 4.6]{NRSSZ}.
\end{rmk}

\begin{lemma}
\label{long vs -}
	There is a strict isomorphism of $\Aug_+$-bimodules $\mcalQ_-\simeq \mcalQ_+^\mrm{long}$, given by \[q^-_j\mapsto q_j^+\,.\]	
\end{lemma}
\begin{proof}
    Using the same notations from the the above formula \eqref{bQ- formula}, we have
    \begin{align*}
   &\quad\, s^{-1}(b^{\mcalQ_+^\mrm{long}}_{i|j})^{\ep_1,\cdots,\ep_{j+1}}_{\ep_{j+2},\cdots,\ep_{m+1}}[s_i^+|\cdots| s_1^+| s_{m+1}^+| s_m^+|\cdots | s_{i+2}^+]\\
   &= (-1)^{\sum_{\al < \be}|s_\al||s_\be|}\sum_{q\in \mbfQ}\Coeff_{s_{i+2}^{12}\cdots s_m^{j,j+1}s_{m+1}^{j+1,j+2}s_1^{j+2,j+3}\cdots s_i^{m,m+1}}\left(\del^{\ep'}_{m+1}q^{1,m+1}\right)\cdot q^+
\end{align*} where $\ep' = (\ep_1,\cdots,\ep_{j+1},\ep_{j+2},\cdots,\ep_{m+1})$, by \eqref{b+ formula} and the fact that $\mcalQ_+^\mrm{long}\inn\Aug_+$ is a sub-bimodule.

 Now we compare the two formulas. The signs agree, so it remains to compare the coefficients, which are determined by the first line of \eqref{matrix formulas} and $\ep,\ep'$. By definition of $\Phi$, for any $q\in\mbfQ$ and $s_1,\cdots,s_N\in \mbfQ\cup \{x_k,y_k\}_{k=1}^\ell$, the term $\Phi(\del q)$ makes equal contributions to $\Coeff_{s_1^{\al_1\al_2}\cdots s_N^{\al_N \al_{N+1}}}(\del_mq^{\al_1\al_{N+1}})$ for any $1\leq \al_1,\cdots,\al_{N+1}\leq m$, as long as each $s_r^{\al_r \al_{r+1}}$ is a valid generator. In particular, since the subscripts of the $s$'s in the two coefficients have the same order, there is a bijection between the terms in $\del_{m+1}q^{i+2,i+1}$ and $\del_{m+1}q^{1,m+1}$ which will contribute to the two coefficients. Moreover, the superscripts of the $s$'s and the pure augmentations $\ep,\ep'$ are such that each pair of these corresponding terms make equal contributions to the coefficients in the twisted differentials $\del_{m+1}^\ep q^{i+2,i+1}$ and $\del_{m+1}^{\ep'}q^{1,m+1}$. Similarly, the term $YA_j- (-1)^{|a_j|}A_jY$ also make the same contributions to both.
\end{proof}

By \Cref{cone for long} and \Cref{long vs -}, we get:  
\begin{thm}
\label{Q+- triangle}
    The diagonal bimodule $\Aug_+$ is the mapping cone of a morphism \[a:\pi^*\mcalC[-1]\to\mcalQ_-\,.\] If $\rot(\LA) = 0$, then this gives an exact triangle 
    \[\begin{tikzcd}
        \mcalQ_-\rar &\Aug_+\rar{\pi} &\pi^*\mcalC\rar{a} & {}
    \end{tikzcd}\]
    in the homotopy category $H^0(\bimod{\Aug_+})$.
\end{thm}

By \eqref{b+ formula} and \Cref{long vs -}, we can write out the morphism $a$ explicitly: for any $i,j\geq 0$, augmentations $\ep_1,\cdots,\ep_{j+1},\ep_{j+2},\cdots,\ep_{m+1}$ where $m=i+j+1$,   $s_1,\cdots,s_j,s_{j+2},\cdots,s_m\in \mbfQ\cup\{x_k,y_k\}_{k=1}^\ell$, and $s_{j+1}\in \{x_k,y_k\}_{k=1}^\ell$, we have
\begin{equation}
\label{h formula}
\begin{aligned}
&\quad\, s^{-1}(a_{i|j})^{\ep_1,\cdots,\ep_{j+1}}_{\ep_{j+2},\cdots,\ep_{m+1}}[s^+_m|\cdots|s^+_{j+2}|s^+_{j+1}|s^+_j|\cdots | s^+_1]\\ 
&= (-1)^{\sum_{\al < \be}|s_\al||s_\be|}\sum_{q\in \mbfQ}\Coeff_{s_1^{12}\cdots s_m^{m,m+1}}(\del^\ep_{m+1}q^{1,m+1})\cdot q^-
\end{aligned}
\end{equation}
where $\ep = (\ep_1,\cdots,\ep_{j+1},\ep_{j+2},\cdots,\ep_{m+1})$.

\section{\texorpdfstring{$A_\infty$}{A-infinity} Sabloff Duality and Higher Homotopies}
\label{A infty Sabloff and higher homotopies section}
This section contains the proofs of \Cref{intro A infty Sabloff} and \Cref{intro homotopies}.
\subsection{Composable \texorpdfstring{$m$}{m}-Copy LSFT Algebras} This subsection studies the composable $m$-copy LSFT algebra of Legendrian knots, whose zeroth associated graded is the composable $m$-copy Chekanov-Eliashberg algebra. We extend \Cref{matrix formula prop} for Legendrian knots, by writing down the $m$-copy Hamiltonian in terms of the $1$-copy Hamiltonian plus thin disks. We also prove a consistency result analogous to \Cref{m-copy CE consistency}.

Let $\LA$ be a Legendrian knot with a single marked point $\star$, a base point $\bullet$, a Maslov potential $\mu: \{\bullet\}\to\R$, and a Morse function $f: \LA\to \R$. Label Reeb chords of $\LA$ by $a_1,\cdots,a_n$, local maxima of $f$ by $b_1,\cdots,b_\ell$, and local minima of $f$ by $c_1,\cdots,c_\ell$. Let $\Lambda^m_f$ be its Lagrangian $m$-copy, with base points, marked points, and Maslov potential as specified in \Cref{m-copy}, where the reader can also find a description of long and short Reeb chords of $\LA_f^m$.
\begin{defn}
\label{composable m-copy LSFT def}
 The \tit{composable m-copy LSFT algebra} $\mcalL_m$ of $\LA$ is the composable LSFT algebra $\mcalL(\LA_f^m)$. Denote its Hamiltonian by $h_{\mcalL_m}$,  its curvature by $F_m$, and the internal link grading of its generators by $e_-^m,e_+^m$ (see \Cref{internal link grading def}).
\end{defn}
\begin{rmk}
Since $\LA$ is a Legendrian knot with a single base point, the $m$-copy link grading of $\LA^m_f$ agrees with the internal link grading of $\LA^m_f$ as a Legendrian link itself. \Cref{composable m-copy LSFT def} takes advantage of this fact, but for Legendrian links with possibly multiple marked points on some components, one needs to reconcile three things: the internal link grading of $\LA$, the $m$-copy link grading of $\LA_f^m$, and the internal link grading of $\LA_f^m$. 
\end{rmk}
 The generators of $\mcalL_m$ are:
\begin{outline}
    \1 $q^{\al\be}_j$, $p^{\al\be}_j$, for each long Reeb chord $a^{\al\be}_j$;
    \1 $x^{\al\be}_k$, $\bar{x}^{\al\be}_k$ (resp. $y^{\al\be}_k$, $\bar{y}^{\al\be}_k$), for each short Reeb chord $b^{\al\be}_k$ (resp. $c^{\al\be}_k$);
    \1 $t^\al$ and $(t^\al)^{-1}$, for each marked point $\star^\al$\,.
\end{outline}
The gradings, which we consider to take values in $\Z/2\rot(\LA)$, are \begin{align*}
    &|q_j^{\al\be}| = |q_j|\,,\ |x_k^{\al\be}| = 0\,,\ |y_k^{\al\be}| = -1\,,\ |t^\al| = |(t^\al)^{-1}| = 0\,,\\
    &|p_j^{\al\be}| = |p_j| = -1 - |q_j|\,,\ |\bar{x}^{\al\be}_k| = -1\,,\ |\bar{y}^{\al\be}_k| = 0\,.
\end{align*}
The filtration degrees are
\begin{equation}
\label{composable m-copy LSFT filt degrees}
    \begin{cases}
        \mcalF(q_j^{\al\be}) = \mcalF(x_k^{\al\be}) = \mcalF(y_k^{\al\be}) = \mcalF((t^\al)^{\pm1}) = 0\\
        \mcalF(p_j^{\al\be}) = \mcalF(\bar{x}_k^{\al\be}) = \mcalF(\bar{y}_k^{\al\be}) = 1
    \end{cases}\hspace{-1ex}.
\end{equation}
Recall from \Cref{m-copy} that $\mbfQ_m=\{q^{\al\be}_j,x^{\al\be}_k,y^{\al\be}_k\}$. Write $\mbfP_m = \mbfP(\LA_f^m) = \{p^{\al\be}_j,\bar{x}^{\al\be},\bar{y}^{\al\be}\}$, $\mbfR_m = \mbfR(\LA_f^m) = \{q^{\al\be}_j,x_k^{\al\be},y_k^{\al\be},p^{\al\be}_j,\bar{x}^{\al\be}_k,\bar{y}^{\al\be}_k\}$, and \[\mbfS_m = \mbfS(\LA_f^m) = \{q_j^{\al\be},x_k^{\al\be},y_k^{\al\be},p^{\al\be}_j,\bar{x}_k^{\al\be},\bar{y}_k^{\al\be},(t^\al)^{\pm 1}\}\,.\] Let $\mbfR_m^+ = \{s\in\mbfR_m\mid e_-^m(s) < e_+^m(s)\}$ and $\mbfS_m^\geq = \{s\in\mbfS_m\mid e_-^m(s)\geq e_+^m(s)\}$. Note that for generators in $\mbfQ_m$ with superscripts ${}^{\al\be}$, we have $e_-^m=\al$ and $e_+^m=\be$. For generators in $\mbfP_m$, with superscripts ${}^{\al\be}$, we have $e_-^m=\be$ and $e_+^m=\al$.
\begin{defn}
Let $s\in\mbfS_m$. If $e_-^m(s) = e_+^m(s)$, we say $s$ is \tit{pure}; if $e_-^m(s)\neq e_+^m(s)$, then $s$ is \tit{mixed}; if $e_-^m(s) > e_+^m(s)$, then $s$ is \tit{positive}; if $e_-^m(s) < e_+^m(s)$, then $s$ is \tit{negative}. 
\end{defn}



Just like $\del_m$ is determined by $\del$, we will show $h_{\mcalL_m}$ is determined by $h$. We first fix some notations. For the rest of this subsection, assume $f$ has exactly two critical points. Define the following elements of $\mcalL_m$: \begin{align*}
 &A_j = \sum_{\al,\be} q^{\al\be}_k\,,\ X = \sum_{\al<\be} x^{\al\be}\,,\ Y = \sum_{\al<\be} y^{\al\be}\,,\ T = \sum_\al t^\al\\ 
 &\bar{A}_j = \sum_{\al,\be} p^{\al\be}_k\,,\ \bar{X} = \sum_{\al<\be} \bar{x}^{\al\be}\,,\ \bar{Y} = \sum_{\al<\be} \bar{y}^{\al\be}\,. \end{align*} Let $\Phi:\mcalL\to\mcalL_m$ be the $\k$-algebra map given by \[\Phi(q_j) = A_j\,,\ \Phi(p_j) = \bar{A}_j\,,\ \Phi(t) = T(1+X)(1+\bar{Y})\,,\ \Phi(t^{-1}) = (1+\bar{Y})^{-1}(1+X)^{-1}T^{-1}\,,\] which clearly descends to the cyclic quotient, denoted by $\Phi^\mrm{cyc}:\mcalL^\mrm{cyc} \to \mcalL_m^\mrm{cyc}$. 

\begin{prop}
 \label{LSFT matrix formula prop}
     Assume the single marked point $\star$ and the two critical points of $f$ are positioned the same way as in \Cref{matrix formula prop}. Then the Hamiltonion of $\mcalL_m$ is given by \[h_{\mcalL_m} = h^\text{thick} + h^\text{thin}\] where $h^\text{thin} = h_\text{triv} + h_\text{crit} + h_\text{crit-cross} + h_\text{cross}$ and
     \begin{align*}
         &h^\text{thick} = \Phi^\mrm{cyc}(h)\\
         &h_\text{triv} = -[\bar{X}X\bar{X}] - \sum_{\be<\al<\al'<\be'} [x^{\al\al'}\bar{x}^{\be\al'}x^{\be\be'}\bar{x}^{\al\be'}]-[Y\bar{Y}Y] - \sum_{\al<\be<\be'<\al'}y^{\al\al'}\bar{y}^{\be\al'}y^{\be\be'}\bar{y}^{\al\be'}\\
         &\qquad\quad + \sum_{\substack{a_j\text{ positive}\\\al<\be,\al'<\be'}} q_j^{\al\al'}p_j^{\be\al'}q_j^{\be\be'}p_j^{\al\be'} + + \sum_{\substack{a_j\text{ negative}\\\al<\be,\al'>\be'}} q_j^{\al\al'}p_j^{\be\al'}q_j^{\be\be'}p_j^{\al\be'}\\
         & h_\text{crit} = [\bar{X}(1+X)\kappa^+((1+\bar{Y})Y) ] - [ T(1+X)\bar{X}T^{-1}\kappa^+(Y(1+\bar{Y})) ]\\ 
         & h_\text{cross} = -\sum_{a_i^+\to a_j^-} [\kappa^+(A_i\bar{A}_i)\bar{A}_jA_j]-\sum_{a_i^-\to a_j^+} [\kappa^+(\bar{A}_iA_i)A_j\bar{A}_j]\\
         &\quad\qquad\ + \sum_{a_i^+\to a_j^+,\,i\neq j}[\kappa^+(A_i\bar{A}_i)A_j\bar{A}_j]+ \sum_{a_i^-\to a_j^-,\,i\neq j}[\kappa^+(\bar{A}_iA_i)\bar{A}_jA_j]\\
         & h_\text{crit-cross} = \sum_j - [ T(1+X)\bar{X}T^{-1}\kappa^+(A_j\bar{A}_j- \bar{A}_jA_j) ]+ [Y(1+\bar{Y})\kappa^-(A_j\bar{A}_j-\bar{A}_jA_j)]\,.
      \end{align*} where $[\cdot]:\mcalL\to \mcalL^\mrm{cyc}$ denotes the quotient map, $\kappa^+,\kappa^-:\mcalL\to\mcalL$ are given by
      \[\kappa^+(w) = \sum_{\al<\be}e_\al we_\be \,,\ \kappa^-(w) = \sum_{\al >\be }e_\al we_\be \] and $a_i^\pm\to a_j^\pm$ means $a_i^\pm$ appears before $a_j^\pm$ according to the orientation of $\LA$.      
 \end{prop}
 \begin{proof}
    Disks of $\LA_f^m$ come in two kinds: \tit{thin} disks lie in a small neighborhood of $\LA_{xy}$ containing $\pi_{xy}(\LA_f^m)$; all other disks are \tit{thick}. As the perturbation parameter $\ep$ tends to zero, thin disks limit to gradient trajectories of $f$ on $\LA$, and thick disks limit to disks of $\LA_{xy}$. (See \cite{Mis03} for a more detailed discussion.) Denote the sum of thick (resp. thin) disks by $h^\mrm{thick}$ (resp. $h^\mrm{thin}$), so that $h_{\mcalL_m} = h^\text{thick} + h^\text{thin}.$

    We will look at thick disks first. Given any disk $\DE\in \DE(\LA)$, let \[W(\DE) = \sum_{\substack{\bar{\DE}\in\DE(\LA_f^m)\\ \bar{\DE} \text{ limits to } \DE}}w(\bar{\DE})\in \mcalL_m\,.\] We are going to show that $W(\DE) = \Phi(w(\DE))$. 
    
    To compute $W(\DE)$ for a given $\DE$, we consider all potential paths $\del\bar{\DE}$ could take as we traverse $\del\DE$. Every time $\DE$ has a positive/negative corner at some Reeb chord $a_j$, the thick disk $\bar{\DE}$ has a positive/negative corner at $a_j^{\al\be}$ for some $1\leq \al,\be\leq m$. Thus, each appearance of $q_j$ (resp. $p_j$) in $w(\DE)$ contributes a factor of $A_j = \Phi(q_j)$ (resp. $\bar{A}_j = \Phi(p_j)$) to $W(\DE)$. Every time $\del\DE$ passes through $\star$ positively, $\del\bar{\DE}$ passes through $\star^\al$ positively for some $\al$, contributing $t^\al = e_\al Te_\al$; then it has either zero or one negative corner at some $x$, contributing \[e_\al(1+X)e_\be = \begin{cases}
        e_\al & \be = \al\\
        x^{\al\be} & \be > \al
    \end{cases}\] for some $\be\geq \al$; finally it has either zero or one positive corner at some $y$, contributing \[e_\be(1+\bar{Y})e_\gamma = \begin{cases}
        e_\be & \gamma = \be\\
        \bar{y}^{\gamma\be} & \gamma < \be
    \end{cases}\] for some $\gamma\leq \be$. Thus, each appearance of $t$ in $w(\DE)$ contributes a factor of $T(1+X)(1+\bar{Y}) = \Phi(t)$ to $W(\DE)$.
    Similarly, every time $\del\DE$ passes through $\star$ negatively, $\del\bar{\DE}$ has $0\leq k_1 < m$ positive corners at $y$'s, contributing $e_\al$ or $(-\bar{y}^{\al_1\al})(-\bar{y}^{\al_2\al_1})\cdots (-\bar{y}^{\be\al_{k_1-1}})$ for some $\al > \al_1 > \cdots > \al_{k_1-1} > \be$, which is a summand of
    \[e_\al(1+\bar{Y})^{-1}e_\be = \begin{cases}
        e_\al & \be = \al\\
        \sum_{\al > \al_1 >\cdots > \al_{k_1-1} > \be}(-\bar{y}^{\al_1\al})(-\bar{y}^{\al_2\al_1})\cdots (-\bar{y}^{\be\al_{k_1-1}}) & \be < \al
    \end{cases}\]
    for some $\be\leq \al$; then it has $0\leq k_2 < m$ negative corners at $x$'s, contributing $e_\be$ or $(-x^{\be\be_1})(-x^{\be_1\be_2})\cdots (-x^{\be_{k_2-1}\gamma})$ for some $\be < \be_1 < \cdots < \be_{k_2-1} < \gamma$, which is a summand of
    \[e_\be(1+X)^{-1}e_\gamma = \begin{cases}
        e_\be & \gamma = \be\\
        \sum_{\be < \be_1 < \cdots < \be_{k_2-1} < \gamma}(-x^{\be\be_1})(-x^{\be_1\be_2})\cdots (-\bar{y}^{\be_{k_2-1}\gamma}) & \gamma > \be
    \end{cases}\]
    for some $\gamma\geq \be$; finally it passes through $\star^\gamma$ negatively for some $\gamma$, contributing $(t^\gamma)^{-1} = e_\gamma T^{-1}e_\gamma$. Thus, each appearance of $t^{-1}$ in $w(\DE)$ contributes a factor of $(1+\bar{Y})^{-1}(1+X)^{-1}T^{-1}= \Phi(t^{-1})$ to $W(\DE)$. All told, we have shown $W(\DE) = \Phi(w(\DE))$. 
    
    To complete the discussion of thick disks, we need to sort out the signs, which are defined just before \Cref{Hamiltonian}. We have already determined the orientation signs of $\bar{\DE}$ at short Reeb chords. The orientation signs at long Reeb chords clearly agree with those of $\DE$. Finally, since the orientation of $\DE_f^m$ copies that of $\DE$, we have $\ep'(\bar{\DE}) = \ep'(\DE)$. Hence, $\sgn(\bar{\DE}) = \sgn(\DE)$ if $\bar{\DE}$ limits to $\DE$. In conclusion, we have
    \begin{align*}h^\mrm{thick} = \sum_{[\bar{\DE}]\in \DE^\mrm{cyc}(\LA_f^m),\,\bar{\DE}\text{ thick}}[\tilde{w}(\bar{\DE})] &=  \sum_{[\DE]\in\DE^\mrm{cyc}(\LA)} [\sgn(\DE)W(\DE)]\\ &= \sum_{[\DE]\in\DE^\mrm{cyc}(\LA)} \Phi^{cyc}([\tilde{w}(\DE)]) = \Phi^\mrm{cyc}(h)\end{align*} as claimed.

 	Now onto thin disks. We group them according to the types of gradient trajectories they limit to: $h_\text{triv}$ includes trajectories of length zero; $h_\text{crit}$ includes those going from the maximum to the minimum; $h_\text{cross}$ includes those connecting distinct $a_j^\pm$'s; $h_\text{crit-cross}$ includes trajectories connecting a critical point to some $a_j^\pm$. See \Cref{thin disks} for illustrations.

    For $h_\text{triv}$, there is one trivial trajectory for each critical point and each crossing. (Since $a_j^+$ and $a_j^-$ project to the same point in $\R^2$, they only get counted once together.) Thin disks limiting to the trivial trajectory of a critical point are triangles and quadrilaterals in the half twist corresponding to the critical point. A triangle at the maximum contributes $[(-\bar{x}^{\al\al'})(-x^{\al\be})(-\bar{x}^{\al'\be})]$ for some $\al<\al'<\be$, which is a summand of $-[\bar{X}X\bar{X}]$. A quadrilateral at the maximum contributes 
   \[[x^{\al\al'}(-\bar{x}^{\be\al'})(-x^{\be\be'})(-\bar{x}^{\al\be'})]\] for some $\be<\al<\al'<\be'$. Similarly, triangles and quadrilaterals at the minimum contribute $-[Y\bar{Y}Y] - \sum_{\al < \be <\be'<\al'}y^{\al\al'}\bar{y}^{\be\al'}y^{\be\be'}\bar{y}^{\al\be'}$. For each $a_j$, we have quadrilaterals at the grid of crossing corresponding to $a_j$. If the crossing $a_j$ is positive, then each of them contributes \[-[(-q_j^{\al\al'})(-p_j^{\be\al'})q_j^{\be\be'}(-p_j^{\al\be'})]\] for some $\al < \be $ and $\al'<\be'$. If the crossing $a_j$ is negative, then a quadrilateral contributes $-[q_j^{\al\al'}p_j^{\be\al'}q_j^{\be\be'}(-p_j^{\al\be'})]$ for some $\al < \be $ and $\al' > \be'$.

    For $h_\text{crit}$, there are two trajectories flowing from the maximum to the minimum. Correspondingly, there are bigons, triangles, and quadrilaterals lying in the dip and ones that go around the knot. A bigon in the dip contributes $-[(-\bar{x}^{\al\be})y^{\al\be}]$ for some $\al<\be$, which is a summand of $\bar{X}Y$. A triangle in the dip contributes either $-[(-\bar{x}^{\al\al'})x^{\al\be}y^{\be\al'}]$ for some $\al < \be <\al'$, which is a summand of  $[\bar{X}XY]$, or $-[(-\bar{x}^{\al\al'}\bar{y}^{\be\al}y^{\be\al'}]$ for some $\be<\al<\al'$, which is a summand of $[\bar{X}\bar{Y}Y]$. A quadrilateral in the dip contributes \[-[(-\bar{x}^{\al\al'})x^{\al\be}\bar{y}^{\be'\be}y^{\be'\al'}]\] for some $\al,\be' < \be < \al'$, which is a summand of $[\bar{X}X\kappa^+(\bar{Y}Y)]$. Adding these up gives $[\bar{X}(1+X)\kappa^+((1+\bar{Y})Y) ]$. Similarly, one can check that the thin disks going around the knot contribute $- [ T(1+X)\bar{X}T^{-1}\kappa^+(Y(1+\bar{Y})) ]$.
    
     For $h_\text{cross}$, there is a unique trajectory connecting $a_i^+$ and $a_j^-$ for any $i, j$, one connecting $a_i^+$ and $a_j^+$ for any $i < j$, and one connecting $a_i^-$ and $a_j^-$ for any $i < j $. Each of these trajectories corresponds to a set of quadrilaterals with two consecutive corners at the grid of crossings corresponding to $a_i$ and the other two at the grid corresponding  to $a_j$. A quadrilateral limiting to the trajectory $a_i^+\to a_j^-$ contributes \[-[q_i^{\al_1\al_2}p_i^{\al_3\al_2}p_j^{\al_4\al_3}q_j^{\al_4\al_1}]\] for some $\al_1<\al_3$ and $1\leq \al_2,\al_4\leq m$, which is a summand of $-[\kappa^+(A_i\bar{A}_i)\bar{A}_jA_j]$. Similarly, quadrilaterals limiting to $a_i^-\to a_j^+$ contribute $-[\kappa^+(\bar{A}_iA_i)A_j\bar{A}_j]$. A quadrilateral limiting to the trajectory $a_i^+\to a_j^+$, contributes \[[q_i^{\al_1\al_2}p_i^{\al_3\al_2}q_j^{\al_3\al_4}p_j^{\al_1\al_4}]\] for some $\al_1 < \al_3$ and $1\leq \al_2,\al_4\leq m$, which is a summand of $[\kappa^+(A_i\bar{A}_i)A_j\bar{A}_j]$. Similarly, quadrilaterals limiting to $a_i^-\to a_j^-$ contribute $[\kappa^+(\bar{A}_iA_i)\bar{A}_jA_j]$. It is straightforward to verify that these contributions do not depend on the crossing types of $a_i$, $a_j$. 

    For $h_\text{crit-cross}$, there is a unique trajectory connecting each critical point to each $a_j^\pm$. Thin disks corresponding to such a trajectory are triangles and quadrilaterals with one or two corners at the critical point and two consecutive corners at the grid of crossings corresponding to $a_j$. A triangle flowing from the maximum to $a_j^+$ contributes \[- [t^\be \bar{x}^{\al\be}(t^\al)^{-1}q_j^{\al\al'}p_j^{\be\al'}]\] for some $\al<\be$ and $1\leq \al'\leq m$, which is a summand of $-[T\bar{X}T^{-1}A_j\bar{A}_j]$. Similarly, triangles flowing from the maximum to $a_j^-$ contributes $[T\bar{X}T^{-1}\bar{A}_jA_j]$. A quadrilateral flowing from the maximum to $a_j^+$ contributes \[- [t^\al x^{\al\al'}\bar{x}^{\be\al'}(t^\be)^{-1}q_j^{\be\be'}p_j^{\al\be'}]\] for some $\be < \al < \al'$ and $1\leq \be'\leq m$, which is a summand of $- [ TX\bar{X}T^{-1}\kappa^+(A_j\bar{A}_j)]$. Similarly, quadrilaterals flowing from the maximum to $a_j^-$ contribute $[ TX\bar{X}T^{-1}\kappa^+(\bar{A}_jA_j)]$. We leave to the reader to verify that triangles and quadrilaterals flowing from $a_j^\pm$ to the minimum contribute $[Y(1+\bar{Y})\kappa^-(A_j\bar{A}_j-\bar{A}_jA_j)]$. The contributions do not depend on the crossing type of $a_j$.
\end{proof}
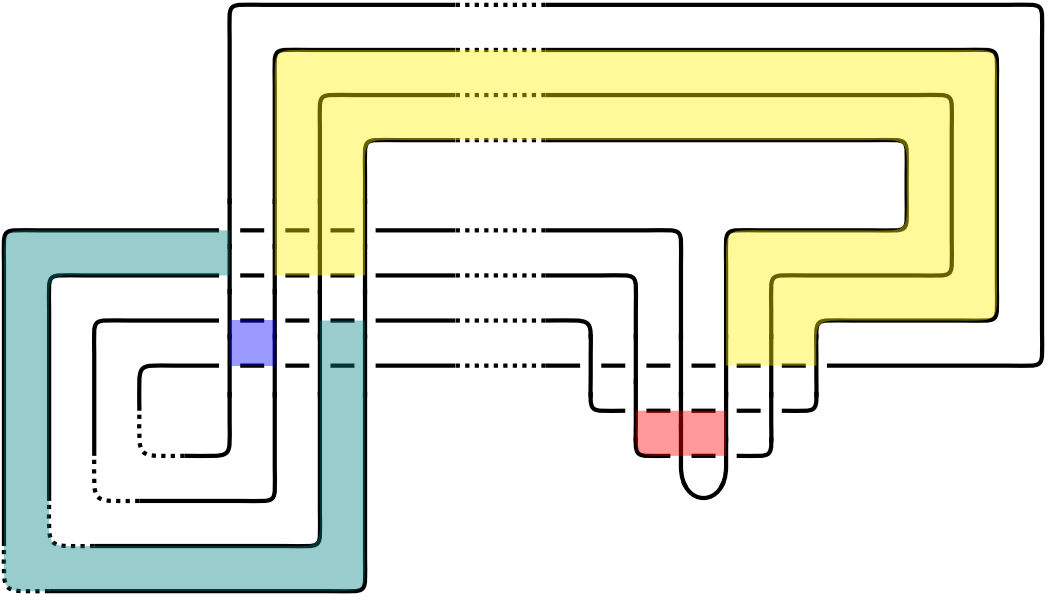
\begin{figure}[t]
\centering
\begin{tikzpicture}[scale = .6] 
   \begin{knot}[
     consider self intersections = true,
     end tolerance = 1pt,
     clip width = 5pt,
     draft mode = off
   ]
   \strand[ultra thick, rounded corners] (-4,0) -- (-2,0) to [out = 0, in = 90, looseness = 2] (-1,-1) -- (-1,-5) to [out =-90,in = -90,looseness = 3]  (0,-5) -- (0,-1) to [out =90,in=180,looseness = 2] (1,0) -- (3,0) to[out =0, in=-90,looseness = 2] (4,1) to[out = 90, in = 0,looseness = 2] (3,2)-- (-4,2);
   \strand[ultra thick, rounded corners,xshift = -3cm](-3,2) -- (-4,2)to[out = 180,in=90,looseness = 2] (-5,1) -- (-5,-7)to [out =-90,in = 0,looseness = 2] (-6,-8) -- (-12,-8);
   \strand[ultra thick, rounded corners,xshift = -3cm,dotted](-12,-8)to[out = 180, in=-90,looseness = 2] (-13,-7);
    \strand[ultra thick, rounded corners,xshift = -3cm](-13,-7)--(-13,-1) to[out = 90, in = 180,looseness = 2](-12,0)--(-3,0);
    \strand[ultra thick, rounded corners] (-4,-1) -- (-3,-1) to [out = 0, in = 90, looseness = 2] (-2,-2) -- (-2,-4) to [out =-90,in = 180,looseness = 2] (-1,-5) -- (0,-5) to [out=0,in=-90,looseness=2] (1,-4) -- (1,-2) to [out =90,in=180,looseness = 2] (2,-1) -- (4,-1)  to[out =0, in=-90,looseness = 2] (5,0) -- (5,2) to[out = 90, in = 0,looseness = 2] (4,3) -- (-4,3);   
    \strand[ultra thick, rounded corners,xshift = -3cm](-3,3)-- (-5,3)to[out = 180,in=90,looseness = 2] (-6,2)--(-6,-6)to [out =-90,in = 0,looseness = 2] (-7,-7) -- (-11,-7);
    \strand[ultra thick, rounded corners,xshift = -3cm,dotted](-11,-7)to[out = 180, in=-90,looseness = 2] (-12, -6);
     \strand[ultra thick, rounded corners,xshift = -3cm](-12, -6)--(-12,-2) to[out = 90, in = 180,looseness = 2](-11,-1)--(-3,-1);
    \strand[ultra thick, rounded corners](-4,-2) to [out = 0,in =90,looseness = 2] (-3,-3) to [out =-90,in=180,looseness = 2] (-2,-4) -- (1,-4) to[out = 0, in=-90,looseness=2] (2,-3) to[out = 90,in=180,looseness=2] (3,-2)--(5,-2)  to[out =0, in=-90,looseness = 2] (6,-1) -- (6,3) to[out = 90, in = 0,looseness = 2] (5,4) -- (3,4)-- (-4,4);  
    \strand[ultra thick, rounded corners,xshift = -3cm](-3,4)--(-6,4)to[out = 180,in=90,looseness = 2] (-7,3) -- (-7,-5)to [out =-90,in = 0,looseness = 2] (-8,-6) --(-10,-6);
    \strand[ultra thick, rounded corners,xshift = -3cm,dotted](-10,-6)to[out = 180, in=-90,looseness = 2] (-11, -5);
     \strand[ultra thick, rounded corners,xshift = -3cm](-11, -5)--(-11,-3) to[out = 90, in = 180,looseness = 2] (-10,-2)--(-3,-2);  
    \strand[ultra thick, rounded corners] (-4,-3) -- (6,-3) to[out =0, in=-90,looseness = 2] (7,-2) -- (7,4) to[out = 90, in = 0,looseness = 2] (6, 5) -- (3,5)--(-4,5);
    \strand[ultra thick, rounded corners,xshift = -3cm](-3,5)-- (-7,5)to[out = 180,in=90,looseness = 2] (-8,4) -- (-8,-4)to[out = -90,in=0,looseness = 2] (-9,-5);
    \strand[ultra thick, rounded corners,xshift = -3cm,dotted](-9,-5)to[out = 180, in=-90,looseness = 2] (-10,-4);
    \strand[ultra thick, rounded corners,xshift = -3cm](-10,-4)to[out = 90,in=180,looseness = 2] (-9,-3)--(-3,-3);
    \strand[ultra thick, dotted] (-4,-3) -- (-6,-3);
    \strand[ultra thick, dotted] (-4,-2) -- (-6,-2);
    \strand[ultra thick, dotted] (-4,-1) -- (-6,-1);
    \strand[ultra thick, dotted] (-4,0) -- (-6,0);
    \strand[ultra thick, dotted] (-4,2) -- (-6,2);
    \strand[ultra thick, dotted] (-4,3) -- (-6,3);
    \strand[ultra thick, dotted] (-4,4) -- (-6,4);
    \strand[ultra thick, dotted] (-4,5) -- (-6,5);
    \flipcrossings{11,12,13,21,22,27};
\end{knot}
\path[xshift = -3cm, fill=teal, fill opacity = .4](-8,0)--(-8,-1)--(-11,-1)to[out=180, in=90,looseness = 2](-12,-2)--(-12,-6) to[out =-90,in=180,looseness =2](-11,-7)--(-7,-7)to[out =0,in=-90,looseness =2] (-6,-6) -- (-6,-2)--(-5,-2) -- (-5,-7)to[out =-90,in=0,looseness =2](-6,-8) -- (-12,-8)to[out =180,in=-90,looseness =2](-13,-7) -- (-13,-1) to[out =90,in=180,looseness =2](-12,0)--cycle;
\path[xshift = -3cm, fill=blue, fill opacity = .4](-7,-2)--(-7,-3)--(-8,-3)--(-8,-2) -- cycle;
\path[fill=red, fill opacity = .4](-2,-4)--(0,-4)--(0,-5)--(-1,-5) to[out = 180,in = -90,looseness=2] cycle;
\path[fill=yellow, fill opacity = .4](-10,-1)--(-8,-1)--(-8,1)to[out =90,in=180,looseness=2](-7,2)--(3,2)to[out=0,in=90,looseness=2](4,1)to[out=-90,in=0,looseness=2](3,0)--(1,0)to[out =180,in=90,looseness=2](0,-1)--(0,-3) -- (2,-3) to[out = 90,in=180,looseness =2] (3,-2) -- (5,-2) to[out=0, in=-90,looseness =2] (6,-1) -- (6,3) to[out=90,in=0,looseness =2] (5,4) -- (-9,4) to[out=180,in=90,looseness =2] (-10,3) -- cycle;
\end{tikzpicture}
\caption{Illustrated here are some examples of thin disks in the $m$-copy of a Legendrian knot. The dotted lines indicate the rest of the $m$-copy, which do not affect the drawn thin disks. The bottom left teal disk contributes to $h_\mrm{cross}$. The top right yellow disk contributes to $h_\mrm{crit-cross}$. The bottom right red disk contributes to $h_\mrm{crit}$. The small blue square contributes to $h_\mrm{triv}$.}
\label{thin disks}
\end{figure}
\subsection{Twisted \texorpdfstring{$m$}{m}-Copy LSFT Differentials}
\label{twisted m-copy LSFT}
In this subsection, we twist the $m$-copy LSFT differentials by curved augmentations. From the twisted $m$-copy differentials, we extract the key ingredients for building the $A_\infty$ Sabloff map and higher homotopies in the next two subsections. 

Let $J^\mrm{mix}\inn \mcalL_m$ be the two-sided ideal generated by mixed generators. Since $d_m$ preserves $e_\pm$, it also preserves $J^\mrm{mix}$. By \Cref{LSFT consistency}, the quotient curved dga $\mcalL_m/J^\mrm{mix}$ is isomorphic to the free product (over $\k^m$) \[\mcalL
_m^{\{1\}}\star\cdots\star \mcalL_m^{\{m\}} = e_1\mcalL e_1\star \cdots\star e_m\mcalL e_m\]with curvature $e_1Fe_1+\cdots+e_mFe_m$, where $F$ is the curvature of $\mcalL$. Note that the underlying algebra of this free product is just $\mcalL^{\times m}$.

For $1\leq \al\leq m$, let $\tep_\al$ be a $c_\al$-curved augmentation $\mcalL$ for some $c_\al\in\k$. Define a map $\tilde{\ep}$ as the composite
\begin{equation}
\label{curved augmentation tuple}
    \begin{tikzcd}
    \mcalL_m\rar &\mcalL_m/J^\mrm{mix}= e_1\mcalL e_1\star \cdots\star e_m\mcalL e_m\arrow{rrr}{(\tep_1,\cdots,\tep_m)} & &&\k^m\bbra{\eta}\,.
\end{tikzcd}
\end{equation} Explicitly, we have 
\[\tep(s) = \begin{cases}
    \tep_\alpha(s)e_\al & e_-^m(s)=e_+^m(s) =\alpha\\
    0 & e_-^m(s)\neq e_+^m(s)
\end{cases}\] for any $s\in\mbfS_m$. By abuse of notation, we also write $\tep = (\tep_1,\cdots,\tep_m)$. Let $\eta_\al\defeq e_\al\eta = \eta e_\al\in \k^m\bbra{\eta}$.

\begin{lemma}
\label{twisted differential lemma}
The map $\tep:\mcalL_m\to\k^m\bbra{\eta}$ defined above satisfies
\[
\begin{cases}
    \tilde{\ep}\circ d_m = 0\\
    \tilde{\ep}(F_m) = c_1\eta_1 + \cdots + c_m\eta_m
\end{cases} \hspace{-1ex}.
\]
\end{lemma}
\begin{proof}
    The first equation follows because $\tep_\al\circ d_m = 0$ for each $\al$. For the second, we have
    \[\tilde{\ep}(F_m) = \tep_1(F)e_1+\cdots +\tep_m(F)e_m = c_1\eta_1 + \cdots + c_m\eta_m\]as claimed.
\end{proof}

 Write $\mcalL_m\agl{\eta}$ for $\mcalL_m$ adjoined a non-commuting variable $\eta$ with $|\eta| = -2$ and $d_m\eta=0$. Let $\widehat{\mcalL_m\agl{\eta}}$ be $\mcalL_m\agl{\eta}$ completed with respect to $\eta$ and $t^1,\cdots,t^m$. (Strictly speaking, what we mean is to complete with respect to $t^\al$'s before imposing the relations $t^\al(t^\al)^{-1} = (t^\al)^{-1}t^\al = e_\al$.) Now define a graded $\k^m$-algebra automorphism $\phi_{\tilde{\ep}}:\widehat{\mcalL_m\agl{\eta}}\to \widehat{\mcalL_m\agl{\eta}}$ by \[\begin{cases}\phi_{\tilde{\ep}}(s) = s+ \tilde{\ep}(s)\\
 \phi_{\tilde{\ep}}(\eta)=\eta\end{cases}	\]
 where $s\in \mbfS_m$. Note that we needed to complete with respect to $t^\al$ so that $\phi_{\tilde{\ep}}((t^\al)^{-1}) = (t^\al+\tilde{\ep}(t^\al))^{-1}$ is well-defined. Define the \tit{twisted $m$-copy LSFT differential} on $\widehat{\mcalL_m\agl{\eta}}$
\[d^{\tilde{\ep}}\defeq \phi_{\tilde{\ep}}\circ d_m\circ \phi_{\tilde{\ep}}^{-1}\,.\]
Let $\mcalL_m^\eta\inn \widehat{\mcalL_m\agl{\eta}}$ denote the closed subalgebra generated by $\mbfR_m\cup \{t^1,\cdots, t^m\}\cup\{\eta\}$.
\begin{lemma}
\label{no constant terms}
    The twisted differential $d^{\tilde{\ep}}$ preserves $\mcalL_m^\eta$, and it has no constant terms (\tit{i.e.}, elements of $\k\bbra{\eta}$).
\end{lemma}
\begin{proof}
    Note that \[\phi_{\tilde{\ep}}((t^\al)^{-1}) = (t^\al+\tilde{\ep}(t^\al))^{-1} = \sum_{k=0}^\infty (-1)^k(t^\al)^k\in\mcalL_m^\eta\]
     so $d^{\tilde{\ep}}$ preserves $\mcalL_m^\eta$. The constant term of $d^{\tilde{\ep}}s = \phi_{\tilde{\ep}}d_ms$ is exactly $\tilde{\ep}\circ d_m (s)$, which is zero by \Cref{twisted differential lemma}.
\end{proof}




Additionally, the algebra $\mcalL_m^\eta$ has two descending filtrations compatible with $d^{\tilde{\ep}}$\,: the filtration $\mcalF$ on $\mcalL_m$ (see \eqref{composable m-copy LSFT filt degrees}) extends to a filtration on $\mcalL_m^\eta$ by setting $\mcalF(\eta) = 1$; another filtration $\mcalF_\eta$ is given by \[\mcalF_\eta(\eta)\defeq 1\,,\ \mcalF_\eta(s)\defeq  0\]
for any $s\in\mbfS_m$.

Let $J_m^\geq \inn \mcalL_m^\eta$ denote the two-sided ideal generated by $\mbfS_m^\geq$, and write $\mcalL_m^+ \defeq \mcalL_m^\eta/J_m^\geq$.

\begin{lemma}
\label{augmented LSFT diff squared}
    The twisted differential $d^{\tilde{\ep}}$ preserves $J_m^\geq$. On $\mcalL_m^+$, we have
    \[(d^{\tilde{\ep}})^2s = [\tilde{\ep}(F_m),s]\] for any $s\in\mbfS_m$.
\end{lemma}
\begin{proof}
    For the first assertion, if $w$ is a word appearing in $d^{\tilde{\ep}}s$ where $s\in \mbfS_m^\geq$, then
    \[e_-^m(w) = e_-^m(s) \geq e_+^m(s) = e_+^m(w)\]since $d^{\tilde{\ep}}$ is $\k^m$-linear. If $w$ contains some $s\in \mbfS_m^\geq$, then we are done. Otherwise, $w$ must be some power of $\eta$, which is impossible by \Cref{no constant terms}.

    For the second one, note that on $\mcalL_m^\eta$, we have \[(d^{\tilde{\ep}})^2s=\phi_{\tilde{\ep}} d^2 \phi_{\tilde{\ep}}^{-1}s = \phi_{\tilde{\ep}} d^2s=\phi_{\tilde{\ep}} [F_m,s] = [\phi_{\tilde{\ep}}(F_m),s+\tilde{\ep}(s)]\,.\] The element $\phi_{\tilde{\ep}}(F_m)$ is given by $\tilde{\ep}(F_m)$ plus non-constant terms. Moreover, each term in $F_m = -\mathord{\bullet}(h_{\mcalL_m})$ is a cyclically composable word, which is either a product of pure generators or has both positive and negative generators as factors. Thus, since $\tilde{\ep}$ is pure, the non-constant terms in $\phi_{\tilde{\ep}}(F_m)$ must lie in $J_m^\geq$. Hence, on $\mcalL_m^+$, we have \[(d^{\tilde{\ep}})^2s = [\tilde{\ep}(F_m),s+\tilde{\ep}(s)] = [\tilde{\ep}(F_m),s]\] where the last equality holds since $k\bbra{\eta}$ is graded commutative.
\end{proof}

Let $\mbfR_m^{+,\eta} = \mbfR_m^+\cup \{\eta_1,\cdots,\eta_m\}$, and set $e_-^m(\eta_\al) = e_+^m(\eta_\al) = \al$ for each $\al$. Then the underlying algebra of $\mcalL_m^+$ is $T_{\k^m}V(\mbfR_m^{+,\eta};e_-^m,e_+^m)$ completed with respect to the variables in $\mbfP_m^+\cup\{\eta_1,\cdots,\eta_m\}$. By \Cref{dualizing derivations} and \Cref{completed dualizing}, dualizing the derivation $d^{\tilde{\ep}}$ on $\mcalL_m^+$ gives us a ``pre-$A_\infty$'' $\k^m$-algebra structure on
\[R_m^{+,\eta}\defeq V(\mbfR_m^{+,\eta};e_-^m,e_+^m)^\vee\,.\]
(Here, by ``pre-$
A_\infty$'' we mean its bar coalgebra is equipped with a coderivation whose square is not necessarily zero.)

Now consider the bar complex $(B^{R_m^{+,\eta}}R_m^{+,\eta}, b^{\tilde{\ep}})$ of the diagonal bimodule $R_m^{+,\eta}$. Since $d_m\eta = 0$, there is a filtered subcomplex
\[E_m^{+,\eta}\defeq B^{R_m^+}Q_m^{+,\eta}\inn B^{R_m^{+,\eta}}R_m^{+,\eta}\] where $R_m^+\inn R_m^{+,\eta}$ is spanned by $\{s^\vee| s\in \mbfR_m^+\}$ and $Q_m^{+,\eta}\inn R_m^{+,\eta}$ is spanned by $\{\eta^\vee\}\cup\{s^\vee| s\in \mbfQ_m^+\}$. 
As an $\k^m$-bimodule, $E_m^{+,\eta}$ decomposes as follows:
\[E_m^{+,\eta} = \bigoplus_{j=0}^\infty R_m^{+,j} = \bigoplus_{j=0}^\infty Q_m^{+,j}\oplus P_m^{+,j}\]
where $R_m^{+,j} = Q_m^{+,j}\oplus P_m^{+,j}$, and $Q_m^{+,j}$ (resp. $P_m^{+,j}$) is spanned by elements of the form $[\cdots |\cdots] \otimes s^\vee\otimes [\cdots |\cdots]$ where $s\in \mbfQ_m^+$ (resp. $s\in \mbfP_m^+$) and the number of times $\eta^\vee$ appears is $j$. In particular,
\[Q_m^{+,0}= B^{Q_m^+}Q_m^+\,,\ P_m^{+,0}= B^{P_m^+}Q_m^+\,.\]
Note that the two filtrations are given by
\begin{align*}
    &\mcalF_jE_m^{+,\eta} = R_m^{+,0}\oplus \cdots\oplus R_m^{+,j-1}\oplus Q_m^{+,j}\\
    &\mcalF_j^\eta E_m^{+,\eta} = R_m^{+,0}\oplus \cdots\oplus R_m^{+,j}\,.
\end{align*}
Since $\mcalF_j^\eta E_m^{+,\eta}\inn \mcalF_{j+1}E_m^{+,\eta}\inn \mcalF^\eta_{j+1}E_m^{+,\eta}\inn \mcalF_{j+2}E_m^{+,\eta}$, we can define a refined filtration $\mcalF'$ on $E_m^{+,\eta}$: \[\mcalF'_{2j}E_m^{+,\eta} \defeq \mcalF_jE_m^{+,\eta},\ \mcalF'_{2j+1}E_m^{+,\eta} \defeq \mcalF^\eta_j E_m^{+,\eta}\] for $j \geq 0$. In other words, the generators have filtration degrees:
\[\mcalF'(\eta) = 2\,,\ \mcalF'(s) = \mcalF(s)\] for any $s\in\mbfS_m$.

By the $\mcalF^\eta$ filtration, $b^{\tilde{\ep}}$ is an upper-triangular matrix:
\[b^{\tilde{\ep}} = 
    \begin{bNiceMatrix}[first-row,first-col]
    & R_m^{+,0} & R_m^{+,1} & R_m^{+,2}  & \cdots\\
    R_m^{+,0} & b^{\tilde{\ep}}_{00} & b^{\tilde{\ep}}_{01} & b^{\tilde{\ep}}_{02} & \cdots \\
    R_m^{+,1} & 0 & b^{\tilde{\ep}}_{11} & b^{\tilde{\ep}}_{12} &  \\ 
    R_m^{+,2}& 0 & 0 & b^{\tilde{\ep}}_{22}&  \\
    \vdots \phantom{A} &\vdots & &&\ddots
    \end{bNiceMatrix}
    \] where $b^{\tilde{\ep}}_{ij}: R_m^{+,j}\to R_m^{+,i}$ (not to be confused with the components $b^{\tilde{\ep}}_{i|j}$ as in the discussion right after \Cref{A infty bimodule def}). Since $R_m^{+,j} = Q_m^{+,j}\oplus P_m^{+,j}$, we can write
\begin{equation}
\label{def of thetas}
    b^{\tilde{\ep}}_{ij} = \begin{bNiceMatrix}[first-row,first-col]
    & Q_m^{+,j} & P_m^{+,j} \\
    Q_m^{+,i} & \tht^{\tilde{\ep}}_{2i,2j} & \tht^{\tilde{\ep}}_{2i,2j+1} \\
    P_m^{+,i} & \tht^{\tilde{\ep}}_{2i+1,2j} & \tht^{\tilde{\ep}}_{2i+1,2j+1}
    \end{bNiceMatrix}\,. \end{equation}
    Now by the refined filtration $\mcalF'$, we have the decomposition
\begin{equation}
\label{upper tri matrix}
   b^{\tilde{\ep}} = 
   \begin{bNiceMatrix}[first-row,first-col, columns-width = 1cm]
   & Q_m^{+,0} & P_m^{+,0} & Q_m^{+,1} & P_m^{+,1} & Q_m^{+,2} & \cdots\\
   Q_m^{+,0} & \tht^{\tilde{\ep}}_{00} & \tht^{\tilde{\ep}}_{01}& \tht_{02}^{\tilde{\ep}}& \tht_{03}^{\tilde{\ep}} & \tht_{04}^{\tilde{\ep}} &\cdots \\
   P_m^{+,0} & 0 & \tht^{\tilde{\ep}}_{11} & \tht_{12}^{\tilde{\ep}} &\tht_{13}^{\tilde{\ep}} & \tht_{14}^{\tilde{\ep}} & \\ 
   Q_m^{+,1}& 0 & 0 & \tht^{\tilde{\ep}}_{22} & \tht^{\tilde{\ep}}_{23} & \tht_{24}^{\tilde{\ep}} & \\
   P_m^{+,1} & 0 & 0 & 0 & \tht^{\tilde{\ep}}_{33} & \tht_{34}^{\tilde{\ep}} & \\
   Q_m^{+,2} & 0 & 0 & 0 & 0 & \tht^{\tilde{\ep}}_{44} & \\
   \vdots \phantom{A} &\vdots &&&&& \ddots
   \end{bNiceMatrix}\,.    
\end{equation}

\begin{lemma}
\label{b squared}
    On $E_m^{+,\eta}$, we have $(b^{\tilde{\ep}})^2_{i|j} = 0$ for $i+j\neq 1$, and the only nonzero values of $(b^{\tilde{\ep}})^2_{0|1}$ and $(b^{\tilde{\ep}})^2_{1|0}$ are
    \begin{equation}
    \label{b tilde ep squared}(b^{\tilde{\ep}})^2_{0|1} [s^\vee| \eta^\vee] = c_{e_-^m(s)}s^\vee,\ (b^{\tilde{\ep}})^2_{1|0} [\eta^\vee |s^\vee] = -c_{e_+^m(s)}s^\vee\end{equation} for any $s\in \mbfR_m^+$. 
\end{lemma}
\begin{proof}
By \Cref{augmented LSFT diff squared} and \Cref{twisted differential lemma}, we have
\[(d^{\tilde{\ep}})^2s = [\tilde{\ep}(F_m),s] = [c_1\eta_1+\cdots+c_m\eta_m,s] = c_{e_-^m(s)}\eta s - c_{e_+^m(s)}s\eta\,. \]
Now the statement follows from $(b^{\tilde{\ep}})^2 = ((d^{\tilde{\ep}})^2)^*$, as in \Cref{dualizing derivations section}.
\end{proof}
Hence, in matrix form, we have
\begin{equation}
     \label{upper tri matrix squared}
    (b^{\tilde{\ep}})^2 = 
    \begin{bNiceMatrix}[first-row,first-col]
    & Q_m^{+,0} & P_m^{+,0} & Q_m^{+,1} & P_m^{+,1} & Q_m^{+,2} & \cdots\\
    Q_m^{+,0} & 0 & 0 & \xi^{\tep}_{02} & 0 & 0 &\cdots \\
    P_m^{+,0} & 0 & 0 & 0 & \xi^{\tep}_{13} & 0 & \\ 
    Q_m^{+,1}& 0 & 0 & 0 & 0 & \xi^{\tep}_{24} & \\
    P_m^{+,1} & 0 & 0 & 0 & 0 & 0 & \\
    Q_m^{+,2} & 0 & 0 & 0 & 0 & 0 & \\
    \vdots\phantom{A}&\vdots &&&&& \ddots
    \end{bNiceMatrix}
\end{equation}
where $\xi^{\tep}_{j,j+2}$ is determined by \eqref{b tilde ep squared}.

Equating the square of the matrix \eqref{upper tri matrix} with \eqref{upper tri matrix squared}, we get a quadratic relation among the $\tht^{\tilde{\ep}}_{ij}$'s for each entry in the upper triangle. We restrict our attention to the first two rows:
   \begin{equation}
   \label{eta homotopy equations}	
   \begin{cases}
   	\tht^{\tilde{\ep}}_{00}\tht^{\tilde{\ep}}_{02}+\tht^{\tilde{\ep}}_{01}\tht^{\tilde{\ep}}_{12}+\tht^{\tilde{\ep}}_{02}\tht^{\tilde{\ep}}_{22} = \xi^{\tep}_{02}\\
 	\tht^{\tilde{\ep}}_{11}\tht^{\tilde{\ep}}_{13}+\tht^{\tilde{\ep}}_{12}\tht^{\tilde{\ep}}_{23}+\tht^{\tilde{\ep}}_{13}\tht^{\tilde{\ep}}_{33} = \xi^{\tep}_{13}\\
	\displaystyle\sum_{i=0}^j \tht^{\tilde{\ep}}_{0i}\tht^{\tilde{\ep}}_{ij}  = 0 & j\neq 2\\
	\displaystyle\sum_{i=1}^j \tht^{\tilde{\ep}}_{1i}\tht^{\tilde{\ep}}_{ij}  =0 &j\neq 3
   \end{cases}
 \end{equation}

 For the purpose of building the $A_\infty$ Sabloff map in the next subsection, we only need the $2\times 2$ submatrix
 \begin{equation}
\label{2 by 2 matrix}
    b^{\tep}_{00}=\begin{bNiceMatrix}[first-row,first-col, columns-width = 1cm]
   & Q_m^{+,0} & P_m^{+,0} \\
   Q_m^{+,0} & \tht^{\tilde{\ep}}_{00} & \tht^{\tilde{\ep}}_{01}\\
   P_m^{+,0} & 0 & \tht^{\tilde{\ep}}_{11}
   \end{bNiceMatrix}   
\end{equation}
of \eqref{upper tri matrix}. Moreover, this submatrix does not depend on the curved lifts $\tep = (\tep_1,\cdots,\tep_m)$:
\begin{lemma}
\label{only depend on ep}
    Let $\ep = (\ep_1,\cdots,\ep_m)$ be the tuple of augmentations induced by $\tilde{\ep}$. The maps $\tht^{\tilde{\ep}}_{00}$, $\tht^{\tilde{\ep}}_{01}$, $\tht^{\tilde{\ep}}_{11}$ only depend on $\ep$, and $\tht^{\tilde{\ep}}_{00} = b^{Q_m^+}_\ep$.
\end{lemma}
\begin{proof}
    Since these maps do not involve $\eta$, setting $\eta = 0$ leaves them unchanged, but that is exactly how $\ep$ is obtained from $\tilde{\ep}$. 
\end{proof}

Hence, to define the maps $\tht^{\tilde{\ep}}_{00}$, $\tht^{\tilde{\ep}}_{01}$, $\tht^{\tilde{\ep}}_{11}$, we can (and will) take $\tilde{\ep}$ to be the trivial lift of $\ep$ to a $0$-curved augmentation, as in \Cref{0-curved augmentations}. Let $b^{P_m^+}_\ep\defeq \tht^{\tilde{\ep}}_{11}$ and $f^\ep\defeq \tht^{\tilde{\ep}}_{01}$.  
\begin{lemma}
\label{m-copy Sabloff relations}
The pair $(P_m^+,b^{P_m^+}_\ep)$ is a $Q_m^+$-bimodule, and $f^\ep: P_m^+\to Q_m^+[1]$ is a morphism of $Q_m^+$-bimodules.
\end{lemma}
\begin{proof}
    By \eqref{eta homotopy equations}, we have $(b^{P_m^+}_\ep)^2 = 0$ and $
    b^{Q_m^+}_\ep f^\ep+f^\ep b^{P_m^+}_\ep = 0$, so $b^{Q_m^+[1]}_\ep f^\ep=f^\ep b^{P_m^+}_\ep$. Because $b^{\tep}$ is a coderivation, $b^{P_m^+}_\ep$ and $f^\ep$ satisfy the appropriate co-Leibniz rules.
\end{proof}

\subsection{\texorpdfstring{$A_\infty$}{A-infinity} Sabloff Duality}
\label{A infty Sabloff}
In this subsection, we prove \Cref{intro A infty Sabloff}. We start by discussing the consistency between the $m$-copy LSFT algebras, which is very similar to the discussion at the end of \Cref{m-copy}. 

Given any $m'$-element subset $I\inn [m] = \{1,\cdots,m\}$, let $J_I\inn \mcalL_m$ be the two-sided ideal generated by \[\{s\in \mbfS_m\mid e_-^m(s)\notin I \text{ or } e_+^m(s)\notin I\}\,.\] Since the LSFT differential $d_m$ preserves $e_\pm$, it also preserves $J_I$, so we get a quotient curved dga $\mcalL_m^I\defeq \mcalL_m/J_I$. The unique order-preserving bijection $\iota:[m']\to I\inn [m]$ defines an injection $h_I:\mbfS_{m'}\to\mbfS_m$ extending \eqref{hI}, by setting \[
    h_I(p^{\alpha\beta}_j) = p^{\iota(\alpha)\iota(\beta)}_j\,,\ h_I(\bar{x}^{\alpha\beta}_k) = \bar{x}^{\iota(\alpha)\iota(\beta)}_k\,,\ h_I(\bar{y}^{\alpha\beta}_j) = \bar{y}^{\iota(\alpha)\iota(\beta)}_k
\] on the extra generators. Hence, we get an inclusion of $\k^{m'}$-algebras $h_I:\mcalL_{m'}\to \mcalL_m$, where $\mcalL_m$ is viewed as an $\k^{m'}$-algebra via the map $\k^{m'}\to \k^m$ induced by $\iota$. This inclusion does not respect the differentials, but we have:

\begin{lemma}
\label{LSFT consistency}
    The composite $\mcalL_{m'}\xrar{h_I} \mcalL_m\to \mcalL_m^I$ is an isomorphism of curved dg $\k^{m'}$-algebras.  
\end{lemma}
\begin{proof}
    Both $d_\DE$ and $\de$ are defined using the Lagrangian projection. Hence, the statement is true simply because deleting any number of copies from $\LA^m_f$ leaves $\LA^{m'}_f$ for some $m' < m$.
\end{proof}

The map $h_I$ extends to $h_I: \mcalL_{m'}^\eta\to \mcalL_m^\eta$ by $h_I(\eta) = \eta$. Let $J_I^\eta\inn \mcalL_m^\eta$ be the two-sided ideal generated by $\{s\in \mbfS_m\mid e_-^m(s)\notin I \text{ or } e_+^m(s)\notin I\}$. Since $d_m$ preserves $J_I$ and $\tilde{\ep}$ is pure, $d^{\tilde{\ep}}$ preserves $J_I^\eta$. Let $(\mcalL_m^\eta)^I = \mcalL_m^\eta/J_I^\eta$ denote the quotient pre-dga. (By ``pre-dga'', we mean a graded algebra equipped with a degree $-1$ derivation whose square is not necessarily zero.) Use the subsequence $\tilde{\ep}_I = \{\tilde{\ep}_\al\mid \al\in I\}$ to define a map $\tilde{\ep}_I: \mcalL_{m'}\to \k^{m'}\bbra{\eta}$, just like \eqref{curved augmentation tuple}. Then we have:
\begin{lemma}
    \label{augmented LSFT consistency}
 The composite $\mcalL_{m'}^\eta\xrar{h_I}\mcalL_m^\eta\to  (\mcalL_m^\eta)^I$ gives an isomorphism of $\k^{m'}$-pre-dgas\[(\mcalL_{m'}^\eta,d^{\tilde{\ep}_{I}})\iso ((\mcalL_m^\eta)^I,d^{\tilde{\ep}})\,.\] Moreover, we can identify $((\mcalL_m^\eta)^I,d^{\tilde{\ep}})  = ((\mcalL_m^I)^\eta,d^{\tilde{\ep}_{I}})$.
\end{lemma}
\begin{proof}
    The first statement follows from \Cref{LSFT consistency} and the definition of $d^{\tilde{\ep}}$. For the second one, note that the $(\mcalL_m^\eta)^I = (\mcalL_m^I)^\eta$ as $\k^{m'}$-algebras, and $d^{\tilde{\ep}} = d^{\tilde{\ep}_I}$ on this algebra because $\tilde{\ep}$ is pure.
\end{proof}
With the above consistency result, we are ready to construct the $A_\infty$ Sabloff map. Let $\mcalP_+$ denote the $\k$-quiver given by \[\mcalP_+(\ep,\ep')\defeq P^+\defeq P^+_2\,.\] Define the codifferential $b^{\mcalP_+}:\mcalP_+\to\mcalP_+$ and the $\Aug_+$-bimodule map\[f:\mcalP_+\to \Aug_+[1]\]as follows. For each augmentation $\ep:\mcalA\to \k$, let $\ep^0:\mcalL\to\k$ denote the canonical $0$-curved augmentation lifting $\ep$, as in \Cref{0-curved augmentations}. Given $i,j\geq 0$ and augmentations $\ep_1,\cdots,\ep_{j+1},\ep_{j+2},\cdots,\ep_{m+1}$,  where $m=i+j+1$, let $\tilde{\ep} = (\ep^0_1,\cdots,\ep^0_{j+1},\ep_{j+2}^0,\cdots,\ep_{m+1}^0)$ and
 \[h_{[m+1]}\defeq h_{m,m+1}\otimes\cdots\otimes h_{j+2,j+3}\otimes h_{j+1,j +2}\otimes h_{j,j+1}\otimes \cdots \otimes h_{12}\,.\]
 Then on 
\[\Aug_+[1](\ep_{j+1},\cdots,\ep_{m+1})\otimes\mcalP_+[1](\ep_{j+1},\ep_{j+2})\otimes\Aug_+[1](\ep_1,\cdots,\ep_{j+1}) \]
define the components of $b^{\mcalP^+}$ and $f$ by
 \[\begin{cases}
  b^{\mcalP^+}_{i|j}\defeq  h_{1,m+1}^{-1}\circ (b^{P_{m+1}^+}_\ep)_{i|j}\circ h_{[m+1]}:Q^+[1]^{\otimes i}\otimes P^+[1]\otimes Q^+[1]^{\otimes j} \to P^+[1]\\
     f_{i|j}\defeq  h_{1,m+1}^{-1}\circ f^\ep_{i|j}\circ h_{[m+1]}:Q^+[1]^{\otimes i}\otimes P^+[1]\otimes Q^+[1]^{\otimes j} \to Q^+[2] 
 \end{cases}\hspace{-1ex}.\]

\begin{lemma}
\label{A infty Sabloff lemma}
 The above maps satisfy the relations	
   \begin{equation}
   \label{A infty Sabloff map relations}	
   \begin{cases}
\displaystyle(b^{\mcalP_+})^2 =0\\
\displaystyle b^{\Aug_+[1]}f - fb^{\mcalP_+} = 0
   \end{cases}\hspace{-1ex}.
 \end{equation}
Hence, $(\mcalP_+,b^{\mcalP_+})$ is an $\Aug_+$-bimodule and $f$ is a morphism of $\Aug_+$-bimodules. 
 \end{lemma}
 \begin{proof}
The proof is entirely analogous to \Cref{Q- def}. Indeed, since the components of $b^{\Aug_+}$, $b^{\mcalP_+}$, and $f$ are built out of $b_\ep^{Q_m^+}$, $b_\ep^{P_m^+}$, and $f^\ep$, respectively, the lemma is reduced to \Cref{m-copy Sabloff relations} by \Cref{augmented LSFT consistency}. 
\end{proof}
The structure maps of the bimodule $\mcalP_+$ count the same set of disks as those of $\mcalQ_-$, defined in \Cref{negative aug bimod}. In fact, we have:
\begin{lemma}
    There is a strict isomorphism of $\Aug_+$-bimodules $\mcalP_+\simeq \mcalQ_-^\dagger[-1]$, given by \[p_k^+\mapsto (-1)^{|p_k|}s^{-1}(q_k^-)^\dagger\,.\]
\end{lemma}
\begin{proof}
Given $i,j\geq 0$, augmentations $\ep_1,\cdots,\ep_{j+1},\ep_{j+2},\cdots,\ep_{m+1}$ where $m=i+j+1$, and $s_1,\cdots,s_j,s_{j+2},\cdots,s_m \in \mbfQ\cup\{x_k,y_k\}_{k=1}^\ell$, we have
\begin{align*}
&\quad\, (\beta^{\mcalP_+}_{i|j})^{\ep_1,\cdots,\ep_{j+1}}_{\ep_{j+2},\cdots,\ep_{m+1}}(s^+_n,\cdots,s^+_{j+2},p_k^+,s^+_j,\cdots , s^+_1)\\ 
&= (-1)^{\sum_{\al_1 < \al_2\leq m}|s_{\al_1}||s_{\al_2}|+\btri}\sum_{p\in \mbfP}\Coeff_{s_1^{12}\cdots s_j^{j,j+1}p_k^{j+2,j+1}s_{j+2}^{j+2,j+3}\cdots s_m^{m,m+1}}(d^{\tep}_{m+1}p^{m+1,1})\cdot p^+
\end{align*} where $s_{j+1} = p_k$, $\tep = (\ep^0_1,\cdots,\ep^0_{j+1},\ep_{j+2}^0,\cdots,\ep_{m+1}^0)$, and $\btri = (m-1) |s_m^+| + \cdots + (j+1)|s_{j+2}^+| + j|p_k^+| + (j-1)|s_j^+| +\cdots +|s_2^+|$. 
On the other hand, by \Cref{linear dual bimodule def} and \eqref{bQ- formula}, we have
\begin{align*}
	&\quad\, (\beta^{\mcalQ_-^\dagger}_{i|j})^{\ep_1,\cdots,\ep_{j+1}}_{\ep_{j+2},\cdots,\ep_{m+1}}(s^+_m,\cdots,s^+_{j+2},(q_k^-)^\dagger,s^+_j,\cdots , s^+_1)\\ 
   &= \sum_{q\in \mbfQ}(-1)^{\sum_{\al_1 < \al_2\neq j+1}|s_{\al_1}||s_{\al_2}|  + \btri'+\bsq}\Coeff_{s_{j+2}^{j+2,j+3}\cdots  s_m^{m,m+1}q^{m+1,1}s_1^{12}\cdots s_j^{j,j+1}}\left(\del^\ep_{m+1}q_k^{j+2,j+1}\right)\cdot (q^-)^\dagger 
\end{align*} where $s_{m+1} = q$, $\ep = (\ep_1,\cdots,\ep_{j+1},\ep_{j+2},\cdots,\ep_{m+1})$, $\btri' = (m-1)|s_j^+|+\cdots +(i+1)|s_1^+| + i|q^-|+(i-1)|s_m^+|+\cdots +|s_{j+3}^+|$, and $\bsq = (|s_{j+2}|+\cdots +|s_m|+i)(|q_k|+|s_1|+\cdots +|s_j|+|q|+j)+(|q_k|+1)(1+i+j)+(i+1)(j+1)$ as in \eqref{linear dual}. 

The terms in the two summations correspond under the bijection $\mbfP\simeq \mbfQ$, $p_k\leftrightarrow q_k$. Moreover, the coefficients count the same set of disks which have a single positive corner at the Reeb chord $a_k^{j+2,j+1}$ and negative corners at the other involved Reeb chords, with the pure chords being augmented by $\tep,\ep$, resepctively. Since $\tep$ lifts $\ep$, one can see from the matching superscripts that augmented generators make equal contributions to both coefficients. It remains to check the signs. By \Cref{SFT bracket}, the coefficients are off by the sign $-(-1)^{(|s_{j+2}|+\cdots +|s_m |+|q|)(|s_1|+\cdots +|s_j|+|p_k|)}$.
The shift in $\mcalQ_-^\dagger[-1]$ brings another sign $(-1)^{1+i
+j+|s_{j+2}^+|+\cdots +|s_m ^+|}$. 

Now it is straightforward to check that the total sign difference is
\[(-1)^{|q_k|(|s_1|+\cdots +|s_j|+|q| + |s_{j+2}|+\cdots +|s_m |) +|p_k|+|p|}=(-1)^{|p_k|+|p|}\]
because $|q_k| + 1\equiv |s_1|+\cdots +|s_j|+|q| + |s_{j+2}|+\cdots +|s_m |$ mod $2$ by \Cref{Hamiltonian}. Hence, the strict pre-morphism $p_k^+\mapsto (-1)^{|p_k|}s^{-1}(q_k^-)^\dagger$ is indeed a morphism $\mcalP_+\to \mcalQ_-^\dagger[-1]$, with an obvious strict inverse.
\end{proof}

By abuse of notation, we write $f: \mcalQ_-^\dagger[-1]\to\Aug_+[1]$, where $\mcalQ_-^\dagger[-1]$ is identified with $\mcalP_+$ via the above lemma. We call $f$ the \tit{$A_\infty$ Sabloff map}, which is justified by the following:

\begin{prop}
\label{extends usual Sabloff}
    The un-suspended component $\varphi_{0|0}$ agrees with the usual Sabloff duality map. In particular, $f$ is a quasi-isomorphism.
\end{prop}
\begin{proof}
    The usual Sabloff map is defined via the \tit{separated $2$-copy} $\td{\LA}^2_f$, obtained by shifting $\LA_1$ very far in the $+z$ direction, so that every mixed Reeb chord goes from $\LA_2$ to $\LA_1$. Since $\LA^2_f$ and $\td{\LA}^2_f$ have the same Lagrangian projections, there are canonical bijections between their LSFT generators:
    \[\iota: \mbfQ^{12}\cup\mbfP^{21}\iso \td{\mbfQ}_2\,,\  \iota:\mbfP^{12}\cup\mbfQ^{21}\iso \td{\mbfP}_2\] where $\mbfQ^{\al\be}\inn\mbfQ_2,\mbfP^{\al\be}\inn\mbfP_2$ are the subsets of generators with superscript ${}^{\al\be}$, and $\td{\mbfQ}_2 = \mbfQ(\td{\LA}_f^2),\td{\mbfP}_2 = \mbfP(\td{\LA}_f^2)$.
    
    The Sabloff map counts disks in $\td{\LA}^2_f$ with a positive corner at some $\td{p}\in \iota(\mbfP^{12})\inn \td{\mbfP}_2$ and a negative corner at some $\td{q}\in \iota(\mbfP^{21})\inn \td{\mbfQ}_2$. (For more details, see for example \cite[Proposition 5.4]{NRSSZ}.) These disks correspond to disks in $\LA^2_f$ with a positive corner at some $p\in\mbfP^{12}$ and a positive corner at some $p'\in \mbfP^{21}$, which define exactly $\varphi_{0|0}$. Hence, $\varphi_{0|0}$ agrees with the Sabloff map, which is a quasi-isomorphism, so $f$ is a quasi-isomorphism. 
\end{proof}

We summarize the results with an explicit formula:
\begin{thm}
    \label{A infty Sabloff thm}
    Let $\LA$ be a Legendrian knot with a single marked point and any Morse function. There is a quasi-isomorphism of $\Aug_+$-bimodules \[f: \mcalQ_-^\dagger[-1]\iso \Aug_+[1]\]extending the linear Sabloff map, specified as follows. Given any $i,j\geq 0$, augmentations $\ep_1,\cdots,\ep_{j+1},\ep_{j+2},\cdots,\ep_{m+1}$ where $m=i+j+1$, and $s_1,\cdots,s_j,s_{j+2},\cdots,s_m \in \mbfQ\cup\{x_k,y_k\}_{k=1}^\ell$, we have
\begin{align*}
&\quad\, s^{-1}(f_{i|j})^{\ep_1,\cdots,\ep_{j+1}}_{\ep_{j+2},\cdots,\ep_{m+1}}[s^+_m|\cdots|s^+_{j+2}|s^{-1}(q_k^-)^\dagger|s^+_j|\cdots |s^+_1]\\ 
&= (-1)^{|p_k|+\sum_{\al < \be}|s_\al||s_\be|}\sum_{q\in\mbfQ\cup\{x_k,y_k\}_{k=1}^\ell}\Coeff_{s_1^{12}\cdots s_j^{j,j+1}p_k^{j+2,j+1}s_{j+2}^{j+2,j+3}\cdots s_m^{m,m+1}}(d^\ep_{n+1}q^{n+1,1})\cdot s(q^+)
\end{align*}
where $s_{j+1} = p_k$, $\ep = (\ep_1,\cdots,\ep_{j+1},\ep_{j+2},\cdots,\ep_{m+1})$, and $s^{\pm1}$ refers to the degree shift.
\end{thm}
\subsection{Higher Homotopies}\label{higher homotopies}
In the previous subsection, we built the $A_\infty$ Sabloff map using the $2\times 2$ submatrix \eqref{2 by 2 matrix}. This subsection deploys the rest of the upper-triangular matrix \eqref{upper tri matrix} to prove \Cref{intro homotopies}, which supplies an explicit homotopy inverse together with all higher homotopies.

For each $j\geq 0$, define a map $\la_j: R_m^{+,0}\to R_m^{+,j}$ of degree $-2j$ as follows:
\begin{align*}
    &\quad\ \la_j([s_{1}|\cdots |s_{k} ]\otimes s\otimes [s_{k+1}|\cdots |s_{k+\ell} ])\\
    &= \sum_{r=0}^j \sum_{\substack{0\leq k_1 < \cdots < k_r\leq k\\ k+1\leq k_{r+1} < \cdots < k_j\leq k+\ell}} \begin{multlined}
 	[s_1|\cdots |s_{k_1}|\eta^\vee| s_{k_1+1}|\cdots |s_{k_r}|\eta^\vee| s_{k_r+1}|\cdots|s_k ]\otimes s\\
    \otimes [s_{k+1}|\cdots |s_{k_{r+1}}|\eta^\vee| s_{k_{r+1}+1}|\cdots |s_{k_j}|\eta^\vee| s_{k_j+1}|\cdots|s_{k+\ell} ]
 \end{multlined}
\end{align*}
where $s_1,\cdots,s_{k+\ell}\in \{s^\vee\mid s\in\mbfQ_m^+\}$ and $s\in \{s^\vee\mid s\in\mbfR_m^+\}$. In short, $\la_j([\cdots |\cdots ]\otimes s\otimes [\cdots |\cdots ])$ is the sum of all possible results of inserting $j$ copies of $\eta^\vee$ into $[\cdots |\cdots ]\otimes s\otimes [\cdots |\cdots ]$. Note that $\la_0$ is just the identity. 

The following lemma tells us that only the first two rows in the matrix \eqref{upper tri matrix} have information, and the remaining rows are just shifted copies of the first two.

\begin{lemma}
\label{inserting eta}
    For any $i,j\geq 0$, we have \[\begin{cases}b^{\tilde{\ep}}_{i,i+j}\la_{i+j} = \la_i b^{\tilde{\ep}}_{0j}\la_j\\\tht^{\tilde{\ep}}_{ij}\la_{\flr{j/2}} = \la_{\flr{i/2}} \tht^{\tilde{\ep}}_{i-2\flr{i/2}, j-2\flr{i/2}}\la_{\flr{j/2}-\flr{i/2}}\end{cases}\hspace{-1ex}.\]
\end{lemma}
\begin{proof}
The second equation follows from the first by definition of the $\tht$'s in \eqref{def of thetas}. For the first equation, the proof below is hard to read, but the idea is simple: both sides insert a total of $i+j$ copies of $\eta^\vee$ and the codifferential $b^{\tilde{\ep}}$ eliminates $j$ copies of them. 

The left-hand side is given by
\begin{align*}
	&\quad\ b^{\tilde{\ep}}_{i,i+j}\la_{i+j}([s_{1}|\cdots |s_{k} ]\otimes s\otimes [s_{k+1}|\cdots |s_{k+\ell} ])\\
	&= b^{\tilde{\ep}}_{i,i+j}\sum_{r=0}^{i+j}\sum_{\substack{0\leq k_1 < \cdots < k_r\leq k\\ k+1\leq k_{r+1} < \cdots < k_{i+j} \leq k+\ell}} \begin{multlined}
 	[s_1|\cdots |s_{k_1}|\eta^\vee| \cdots |s_{k_r}|\eta^\vee| \cdots|s_k ]\otimes s\\
    \otimes [s_{k+1}|\cdots |s_{k_{r+1}}|\eta^\vee| \cdots |s_{k_{i+j}}|\eta^\vee| \cdots|s_{k+\ell} ]
 \end{multlined}\\
 &=\sum_{r=0}^{i+j}\sum_{\substack{0\leq k_1 < \cdots < k_r\leq k\\ k+1\leq k_{r+1} < \cdots < k_{i+j} \leq k+\ell}} \sum_{r'=0}^r\sum_{\substack{k_{r'} < k_0 < k_{r'+1}\\ k_{j+r'} < k_{i+j+1} < k_{j+r'+1}}}\\
 &\quad \begin{multlined}
 	(-1)^{\bsq} [s_1|\cdots |s_{k_{1}}|\eta^\vee|\cdots | s_{k_{r'}}|\eta^\vee|\cdots|s_{k_0-1}]\otimes
 	b^{\tilde{\ep}}_{k-k_0+r-r'+1|k_{i+j+1}-k+j-r}\big([s_{k_0}|\\
 	\cdots|s_{k_{r'+1}}|\eta^\vee| \cdots |s_{k_r}|\eta^\vee|\cdots|s_k ]\otimes s\otimes[s_{k+1}|\cdots|s_{k_{r+1}}|\eta^\vee| \cdots
 	|s_{k_{j+r'}}|\eta^\vee|\\
 	 \cdots|s_{k_{i+j+1}}]\big)\otimes[s_{k_{i+j+1}+1}|\cdots |s_{k_{j+2}}|\eta^\vee|\cdots |s_{k_{j+i-r'+1}}|\eta^\vee|\cdots  |s_{k+\ell} ]
 \end{multlined}
\end{align*} where $\bsq = \|s_1\|+\cdots +\|s_{k_0-1}\|$. The right-hand side is given by
\begingroup
\allowdisplaybreaks
\begin{align*}
    &\quad\ \la_i b^{\tilde{\ep}}_{0j}\la_j([s_{1}|\cdots |s_{k} ]\otimes s\otimes [s_{k+1}|\cdots |s_{k+\ell} ])\\
 &= \la_i\sum_{r=0}^j \sum_{\substack{0\leq k_1 < \cdots < k_r\leq k\\ k+1\leq k_{r+1} < \cdots < k_j\leq k+\ell}}\ \sum_{\substack{k_0 < k_1\\k_j < k_{j+1}}}
 \begin{multlined}
 	(-1)^{\btri} [s_1|\cdots |s_{k_0-1}]\otimes b^{\tilde{\ep}}_{k-k_0+r+1|k_{j+1}-k+j-r}\big([s_{k_0}|\\
 	\cdots|s_{k_1}|\eta^\vee| \cdots |s_{k_r}|\eta^\vee| \cdots|s_k ]\otimes s\otimes[s_{k+1}|\cdots|s_{k_{r+1}}|\eta^\vee|\\
 	\cdots
 	|s_{k_j}|\eta^\vee|\cdots|s_{k_{j+1}}]\big)\otimes[s_{k_{j+1}+1}|\cdots |s_{k+\ell} ]
 \end{multlined}\\
 &= \sum_{r=0}^j \sum_{\substack{0\leq k_0 < k_1 < \cdots < k_r\leq k\\ k+1\leq k_{r+1} < \cdots < k_j < k_{j+1}\leq k+\ell}}\sum_{r'=0}^i\sum_{\substack{0\leq k_{-r'} < \cdots < k_{-1} < k_0\\ k_{j+1} < k_{j+2} < \cdots < k_{j+i-r'+1}\leq k+\ell}}\\
 &\qquad \begin{multlined}
 	(-1)^{\btri} [s_1|\cdots |s_{k_{-r'}}|\eta^\vee|\cdots | s_{k_{-1}}|\eta^\vee|\cdots|s_{k_0-1}]\otimes
 	b^{\tilde{\ep}}_{k-k_0+r+1|k_{j+1}-k+j-r}\big([s_{k_0}|\\
 	\cdots|s_{k_1}|\eta^\vee| \cdots |s_{k_r}|\eta^\vee|\cdots|s_k ]\otimes s\otimes[s_{k+1}|\cdots|s_{k_{r+1}}|\eta^\vee| \cdots
 	|s_{k_j}|\eta^\vee|\\
 	\cdots|s_{k_{j+1}}]\big)\otimes[s_{k_{j+1}+1}|\cdots |s_{k_{j+2}}|\eta^\vee|\cdots |s_{k_{j+i-r'+1}}|\eta^\vee|\cdots  |s_{k+\ell} ]
 \end{multlined}
\end{align*}
\endgroup where $\btri= \|s_1\|+\cdots + \|s_{k_0-1}\|$. We see that the two final expressions are equal by rearranging the sums and reindexing: replace $r$ by $r + r'$; replace $k_1 < \cdots < k_{r'} < k_0 < k_{r'+1} < \cdots < k_{j+r'} < k_{i+j+1} < k_{j+r'+1} < \cdots < k_{i+j}$ by $k_{-r'} < \cdots < k_{-1} < k_0 < k_1 < \cdots < k_r < k_{r+1} < k_{r+2} <\cdots <k_{j+1}$.
\end{proof}

For each $j > 0$, define the following maps: 
\[ f^{\tilde{\ep}}_j\defeq  \begin{cases}
    \displaystyle s\circ\tht^{\tilde{\ep}}_{0j}\circ \la_{j/2}\circ s^{-1}: Q_m^{+,0}[1]\to Q_m^{+,0}[1] & j \text{ even}\\
    \displaystyle -s\circ \tht^{\tilde{\ep}}_{0j}\circ \la_{(j-1)/2}: P_m^{+,0}\to Q_m^{+,0}[1] & j\text{ odd}\\
\end{cases}\]\[g^{\tilde{\ep}}_j\defeq\begin{cases}
 \displaystyle   -\tht^{\tilde{\ep}}_{1,j+1}\circ \la_{j/2}: P_m^{+,0}\to P_m^{+,0} & j \text{ even}\\
    \displaystyle \tht^{\tilde{\ep}}_{1,j+1}\circ \la_{(j+1)/2}\circ s^{-1}:Q_m^{+,0}[1]\to P_m^{+,0} & j\text{ odd}	
 \end{cases}
\]
of degree $1-j$, where $s:Q_m^{+,0}\to Q_m^{+,0}[1]$ is the obvious degree $-1$ isomorphism. Up to shifts, we can write $f_j^{\tilde{\ep}} = (-1)^j\tht_{0j}^{\tilde{\ep}}\la_{\flr{j/2}}$ and $g_j^{\tilde{\ep}} = (-1)^{j+1}\tht^{\tilde{\ep}}_{1,j+1}\la_{\flr{(j+1)/2}}$. Note that $f_1^{\tep} = f^\ep$, by \Cref{only depend on ep}.

\begin{lemma}
The above maps satisfy the following relations:
   \begin{equation}
   \label{m-copy homotopy equations}	
   \begin{cases}
   	\displaystyle b^{Q_m^+[1]}f^{\tilde{\ep}}_2+f^{\tilde{\ep}}_2\,b^{Q_m^+[1]} = -(\xi^{\tep}_{02}\la_1)[1] - f^{\tilde{\ep}}_1g^{\tilde{\ep}}_1\\
 	\displaystyle b^{P_m^+}g^{\tilde{\ep}}_2+ g^{\tilde{\ep}}_2\,b^{P_m^+} = -\xi^{\tep}_{13}\la_1-g^{\tilde{\ep}}_1f^{\tilde{\ep}}_1 \\
	\displaystyle b^{Q_m^+[1]}f^{\tilde{\ep}}_j + f^{\tilde{\ep}}_j\,b^{Q_m^+[1]}= - f^{\tilde{\ep}}_1g^{\tilde{\ep}}_{j-1} + f^{\tilde{\ep}}_2f^{\tilde{\ep}}_{j-2} -f^{\tilde{\ep}}_3g^{\tilde{\ep}}_{j-3} +\cdots -f_{j-1}g_1 & j > 2\text{ even}\\
	\displaystyle b^{Q_m^+[1]}f^{\tilde{\ep}}_j - f^{\tilde{\ep}}_j\,b^{P_m^+}= - f^{\tilde{\ep}}_1g^{\tilde{\ep}}_{j-1} + f^{\tilde{\ep}}_2f^{\tilde{\ep}}_{j-2} -f^{\tilde{\ep}}_3g^{\tilde{\ep}}_{j-3} +\cdots +f_{j-1}f_1 & j \text{ odd}\\
\displaystyle b^{P_m^+}g^{\tilde{\ep}}_j + g^{\tilde{\ep}}_j\,b^{P_m^+} = - g^{\tilde{\ep}}_1f^{\tilde{\ep}}_{j-1} + g^{\tilde{\ep}}_2g^{\tilde{\ep}}_{j-2} - g^{\tilde{\ep}}_3f^{\tilde{\ep}}_{j-3}+\cdots  - g^{\tilde{\ep}}_{j-1}f^{\tilde{\ep}}_1& j > 2\text{ even}\\
\displaystyle b^{P_m^+}g^{\tilde{\ep}}_j - g^{\tilde{\ep}}_j\,b^{Q_m^+[1]} = - g^{\tilde{\ep}}_1f^{\tilde{\ep}}_{j-1} + g^{\tilde{\ep}}_2g^{\tilde{\ep}}_{j-2} - g^{\tilde{\ep}}_3f^{\tilde{\ep}}_{j-3}+\cdots  + g^{\tilde{\ep}}_{j-1}g^{\tilde{\ep}}_1& j\text{ odd}\\
   \end{cases}
 \end{equation}
 \end{lemma}
  \begin{proof}
Up to shifts, the idea is to precompose \eqref{eta homotopy equations} with $\la_{\flr{j/2}}$ and then use \Cref{inserting eta}. We will demonstrate for even $j$'s and leave the odd cases to the reader. 

For $j > 0$ even, we have   \[s\tht^{\tilde{\ep}}_{0i}\tht^{\tilde{\ep}}_{ij}\la_{j/2}s^{-1} = s\tht^{\tilde{\ep}}_{0i}\la_{\flr{i/2}} \tht^{\tilde{\ep}}_{i-2\flr{i/2}, j-2\flr{i/2}}\la_{j/2-\flr{i/2}}s^{-1} = \begin{cases}
 f^{\tilde{\ep}}_i f^{\tilde{\ep}}_{j-i} & i \text{ even}\\
-f^{\tilde{\ep}}_ig^{\tilde{\ep}}_{j-i} & i\text{ odd}\\
 \end{cases}\] by \Cref{inserting eta}. Since $\tht_{00}^{\tep} = b^{Q_m^+}_\ep$, applying $s\circ \cdot\circ \la_{j/2}s^{-1}$ to the first and third equations in \eqref{eta homotopy equations} proves the first and third equations in \eqref{m-copy homotopy equations}.
 
Similarly, for $j$ odd (note the index shift between $g_j^{\tilde{\ep}}$ and $\tht_{1j}^{\tilde{\ep}}$), we have
 \[\tht^{\tilde{\ep}}_{1i}\tht^{\tilde{\ep}}_{ij}\la_{j/2} = \tht^{\tilde{\ep}}_{1i}\la_{\flr{i/2}} \tht^{\tilde{\ep}}_{i-2\flr{i/2}, j-2\flr{i/2}}\la_{j/2-\flr{i/2}} = \begin{cases}
- g^{\tilde{\ep}}_{i-1} f^{\tilde{\ep}}_{j-i} & i \text{ even}\\
 g^{\tilde{\ep}}_{i-1}g^{\tilde{\ep}}_{j-i} & 1<i<j\text{ odd}\\
 \end{cases}
\]
 by \Cref{inserting eta}. Since $b^{P_m^+} = \tht_{11}^{\tep}$, applying $\cdot\circ \la_{j/2}$ to the second and fourth equations in \eqref{eta homotopy equations} proves the second and fifth equations in \eqref{m-copy homotopy equations}.
\end{proof}

We know $f_1^{\tilde{\ep}} = f^\ep: P_m^+\to Q_m^+[1]$ is a morphism of $Q_m^+$-bimodules, by \Cref{m-copy Sabloff relations}. The following lemma shows that all the other $f_j^{\tep}$'s and $g_j^{\tep}$'s are also pre-morphisms of $Q_m^+$-bimodules.
\begin{lemma}
\label{f and g are bicomodule maps}
    For any $j > 0$, the maps $f^{\tilde{\ep}}_j$ and $g^{\tilde{\ep}}_j$ are maps of $BQ_m^+$-bicomodules. 
\end{lemma}
\begin{proof}
 We will show $f^{\tilde{\ep}}_{2j}$ commutes with the left coactions: \[(1\otimes f^{\tilde{\ep}}_{2j})\DE_L = \DE_L\circ f^{\tilde{\ep}}_{2j}\] and the remaining cases are entirely analogous. The left-hand side is given by
\begingroup
\allowdisplaybreaks
\begin{align*}
    &\quad\ (1\otimes f^{\tilde{\ep}}_{2j})\DE_L([s_1|\cdots |s_k]\otimes s\otimes [s_{k+1}|\cdots |s_{k+\ell}])\\
    &= \sum_{i=0}^k [s_1|\cdots |s_i]\otimes \tht^{\tilde{\ep}}_{0j}\circ \la_{j}( [s_{i+1}|\cdots |s_k]\otimes s\otimes [s_{k+1}|\cdots | s_{k+\ell}])\\
    &= \sum_{i=0}^k\sum_{\substack{i+1\leq k_0 < k_1 < \cdots < k_r\leq k\\ k+1\leq k_{r+1} < \cdots < k_j < k_{j+1}\leq k+\ell}}\begin{multlined}
    [s_1|\cdots |s_i]\otimes [s_{i+1}|\cdots|s_{k_0-1}]\otimes b^{\tilde{\ep}}_{k-k_0+r+1|k_{j+1}-k+j-r}([s_{k_0}|\\
    \cdots|s_{k_1}|\eta^\vee|\cdots|s_{k_r}|\eta^\vee|\cdots |s_k]\otimes s\otimes [s_{k+1}| \cdots|s_{k_{r+1}}|\eta^\vee|\\
    \cdots|s_{k_j}|\eta^\vee|\cdots |s_{k_{j+1}}])\otimes [s_{k_{j+1}+1}|\cdots | s_{k+\ell}]
 \end{multlined}\,.
\end{align*}
\endgroup
    The right-hand side is given by
     \begin{align*}
        &\quad\ \DE_L\circ f^{\tilde{\ep}}_{2j}([s_1|\cdots |s_k]\otimes s\otimes [s_{k+1}|\cdots |s_{k+\ell}])\\
 		&= \sum_{\substack{0\leq k_0 < k_1 < \cdots < k_r\leq k\\ k+1\leq k_{r+1} < \cdots < k_j < k_{j+1}\leq k+\ell}}\sum_{i=0}^{k_0-1}\begin{multlined}
 				 [s_1|\cdots |s_i]\otimes [s_{i+1} |s_{k_0-1}]\otimes b^{\tilde{\ep}}_{k-k_0+r+1|k_{j+1}-k+j-r}([s_{k_0}|\\
 				 \cdots|s_{k_1}|\eta^\vee|\cdots|s_{k_r}|\eta^\vee|\cdots |s_k]\otimes s\otimes [s_{k+1}| \cdots|s_{k_{r+1}}|\eta^\vee|\\
 				 \cdots|s_{k_j}|\eta^\vee|\cdots |s_{k_{j+1}}])\otimes [s_{k_{j+1}+1}|\cdots | s_{k+\ell}]
 			 \end{multlined}
    \end{align*}
    so the two sides are equal by swapping the sums.
\end{proof}
By the above two lemmas, $f_1^{\tilde{\ep}} = f^\ep: P_m^+\to Q_m^+[1]$ and $g_1^{\tilde{\ep}}: Q_m^+[1]\to P_m^+$ are morphisms of $Q_m^+$-bimodules, and the higher $f^{\tep}_j$'s and $g^{\tep}_j$'s are pre-morphisms of $Q_m^+$-bimodules. 

For the remainder of this subsection, fix $c_1,c_2\in\k$ such that $c_1+c_2 = 1$, and assume every augmentation $\ep:\mcalA\to \k$ lifts to a ($-c_1$)-curved augmentation $\ep^-:\mcalL\to\k$ and a $c_2$-curved augmentation $\ep^+:\mcalL\to\k$. 

Define the following pre-morphisms of $\Aug_+$-bimodules: \[\begin{cases}
     f_{2k}: \Aug_+[1]\to \Aug_+[1] & i > 0\\
     f_{2k+1}: \mcalP_+\to \Aug_+[1] & i>0\\
     g_{2k}: \mcalP_+\to \mcalP_+ & i > 0\\
     g_{2k+1}:\Aug_+[1]\to\mcalP_+ & i\geq 0
 \end{cases}\]as follows. Given $i,j\geq 0$ and augmentations $\ep_1,\cdots,\ep_{m+1}$,  where $m=i+j+1$, let \begin{equation}\label{1/2 curved augmentations}\tilde{\ep} = \tilde{\ep}_{i,j} = (\ep^-_1,\cdots,\ep^-_{j+1},\ep^+_{j+2},\cdots,\ep^+_{m+1})\end{equation} and $h_{[m+1]}= h_{m,m+1}\otimes\cdots\otimes h_{j+2,j+3}\otimes h_{j+1, j +2}\otimes h_{j, j+1}\otimes \cdots \otimes h_{12}$.
 Then on \[\Aug_+[1](\ep_{j+2},\cdots,\ep_{m+1})\otimes\Aug_+[1](\ep_{j+1},\ep_{j+2})\otimes\Aug_+[1](\ep_1,\cdots,\ep_{j+1}) \] define the components of $f_{2k}$ and $g_{2k+1}$ as the compositions
 \[\begin{cases}
     (f_{2k})_{i|j}\defeq  h_{1,m+1}^{-1}\circ (f_{2k}^{\tilde{\ep}})_{i|j}\circ h_{[m+1]}:Q^+[1]^{\otimes i}\otimes Q^+[2]\otimes Q^+[1]^{\otimes j} \to Q^+[2]  \\
     (g_{2k+1})_{i|j}\defeq  h_{1,m+1}^{-1}\circ (g_{2k+1}^{\tilde{\ep}})_{i|j}\circ h_{[m+1]}:Q^+[1]^{\otimes i}\otimes Q^+[2]\otimes Q^+[1]^{\otimes j} \to P^+[1]
 \end{cases}\hspace{-1ex}.\]
 Similarly, on 
\[\Aug_+[1](\ep_{j+2},\cdots,\ep_{m+1})\otimes\mcalP_+[1](\ep_{j+1},\ep_{j+2})\otimes\Aug_+[1](\ep_1,\cdots,\ep_{j+1}) \]
define the components of $f_{2k+1}$ and $g_{2k}$ by
 \[\begin{cases}
     (f_{2k+1})_{i|j}\defeq  h_{1,m+1}^{-1}\circ (f_{2k+1}^{\tilde{\ep}})_{i|j}\circ h_{[m+1]}:Q^+[1]^{\otimes i}\otimes P^+[1]\otimes Q^+[1]^{\otimes j} \to Q^+[2]  \\
     (g_{2k})_{i|j}\defeq  h_{1,m+1}^{-1}\circ (g_{2k}^{\tilde{\ep}})_{i|j}\circ h_{[m+1]}:Q^+[1]^{\otimes i}\otimes P^+[1]\otimes Q^+[1]^{\otimes j} \to P^+[1]
 \end{cases}\hspace{-1ex}.\]
Note that we can talk about components of the $f^{\tep}_j$'s and $g^{\tep}_j$'s because of \Cref{f and g are bicomodule maps}. Since $f_1^{\tep} = f^\ep$, we have $f_1 = f$.
\begin{prop}
\label{homotopy equations lemma}
 The above maps satisfy the relations	
   \begin{equation}
   \label{homotopy equations}	
   \begin{cases}
   	\displaystyle b^{\Aug_+[1]}f_2+f_2b^{\Aug_+[1]} = 1_{\Aug_+[1]} - f_1g_1\\
 	\displaystyle b^{\mcalP_+}g_2+ g_2b^{\mcalP_+} = 1_{\mcalP_+}-g_1f_1 \\
	\displaystyle b^{\Aug_+[1]}f_j + f_jb^{\Aug_+[1]}= - f_1g_{j-1} + f_2f_{j-2} -f_3g_{j-3} +\cdots -f_{j-1}g_1 & j > 2\text{ even}\\
	\displaystyle b^{\Aug_+[1]}f_j - f_jb^{\mcalP_+}= - f_1g_{j-1} + f_2f_{j-2} -f_3g_{j-3} +\cdots +f_{j-1}f_1 & j \text{ odd}\\
\displaystyle b^{\mcalP_+}g_j + g_jb^{\mcalP_+} = - g_1f_{j-1} + g_2g_{j-2} - g_3f_{j-3}+\cdots  - g_{j-1}f_1& j > 2\text{ even}\\
\displaystyle b^{\mcalP_+}g_j - g_jb^{\Aug_+[1]} = - g_1f_{j-1} + g_2g_{j-2} - g_3f_{j-3}+\cdots  + g_{j-1}g_1& j\text{ odd}
\end{cases}
\end{equation}
 \end{prop}
 \begin{proof}
The proof is entirely analogous to \Cref{Q- def}. Indeed, since the components of $f_j$, $g_j$ are given by $f_j^{\tilde{\ep}}$, $g_j^{\tilde{\ep}}$, the lemma is reduced to \eqref{m-copy homotopy equations} by \Cref{augmented LSFT consistency}. In order to use \Cref{augmented LSFT consistency}, it is crucial that deleting any subset of parallel copies changes $\tep_{i,j}$ to $\tep_{i',j'}$ for some $i'\leq i$, $j'\leq j$ (see \eqref{1/2 curved augmentations}).

The first two equations require additional justifications. The above argument and \eqref{m-copy homotopy equations} tell us that on \[\Aug_+[1](\ep_{j+2},\cdots,\ep_{m+1})\otimes\Aug_+[1](\ep_{j+1},\ep_{j+2})\otimes\Aug_+[1](\ep_1,\cdots,\ep_{j+1}) \]            the map $b^{\Aug_+[1]}f_2+f_2b^{\Aug_+[1]} + f_1g_1$ is equal to $-h_{1,m+1}^{-1}\circ (\xi^{\tep}_{02}\la_1)_{i|j}\circ h_{[m+1]}$. By \Cref{b squared}, the pre-morphism $\xi^{\tep}_{02}\la_1: Q_{n+1}^+\to Q_{n+1}^+$ is in fact strict, given by
 \[(\xi^{\tep}_{02}\la_1)_{0|0}(s) = (c_{e_-^m(s)}-c_{e_+^m(s)})s\]
for any $s\in \mbfQ_{n+1}$.       
Hence, $b^{\Aug_+[1]}f_2+f_2b^{\Aug_+[1]} + f_1g_1$ is also strict. On $\Aug_+[1](\ep_1,\ep_2)$, it is given by $-h_{12}^{-1}\circ (\xi^{\tep}_{02}\la_1)_{0|0}\circ h_{12}$, where $\tep = (\ep_1^-,\ep_2^+)$. If $e_-^m(s) = 1$ and $e_+^m(s) = 2$, then
 \[(\xi^{\tep}_{02}\la_1)_{0|0}(s) = (-c_1-c_2)s=-s\] since $c_1+c_2 = 1$ by choice. Therefore, we get $b^{\Aug_+[1]}f_2+f_2b^{\Aug_+[1]} + f_1g_1 = 1_{\Aug_+[1]}$, as claimed. The proof of the second equation is very similar.
\end{proof}
\Cref{homotopy equations lemma} says $f_2$ and $g_2$ are homotopies witnessing $f_1$ and $g_1$ are homotopy inverses, and the higher $f_j$'s and $g_j$'s serve as higher homotopies. When $\rot(\LA) = 0$, so that everything is $\Z$-graded, we can make this statement more precise using the language of $\infty$-categories, as follows.

 Let $J$ be the nerve of the groupoid $0\leftrightarrows 1$ consisting of two objects and a single isomorphism between them. Just like an arrow in a category is an isomorphism iff it extends to a map from the groupoid $0\leftrightarrows 1$ to the category, an arrow $\DE^1\to X$ in an $\infty$-category $X$ is an equivalence iff it extends along the inclusion $\DE^1\to J$. In fact, each equivalence in $X$ has a contractible space of lifts to a map $J\to X$ \cite{Joy02}\cite[Lemma 2.4.5]{Lan21}. 
 
 In our situation, $X$ is the dg nerve $N^\mrm{dg}_\bullet(\bimod{\Aug_+})$, where equivalences are exactly quasi-isomorphisms of bimodules, and the $A_\infty$ Sabloff map $f=f_1: \mcalQ_-^\dagger[-1]\to \Aug_+[1]$ is the quasi-isomorphism in question. It turns out \Cref{homotopy equations lemma} furnishes exactly a lift of $f$ to a map $J\to N^\mrm{dg}_\bullet(\bimod{\Aug})$, which is what we will prove in \Cref{quasi-inverse}. For the definition of the dg nerve of a dg category, we refer the reader to \cite[Section 1.3]{Lur}. What is important to us is how to prescribe a map into a dg nerve:

Let $K_\bullet$ be a simplicial set and $\mcalA$ a dg category. A map of simplicial sets $\psi:K_\bullet\to N^\mrm{dg}_\bullet(\mcalA)$ is equivalent to the following data:
\begin{itemize}
    \item a set map $\psi: K_0\to \ob \mcalA$;
    \item for each $k$-simplex $\sigma: \DE^k\to K_\bullet$ with $k > 0$, an element $\psi(\sigma)\in \mcalA(\psi(\sigma(0)),\psi(\sigma(k)))$ of degree $k-1$,
\end{itemize}
satisfying the following conditions:
\begin{enumerate}[label = (\arabic*)]
    \item $\psi$ sends degenerate edges to identity morphisms;
    \item $\psi$ sends degenerate simplices of dimension $\geq 2$ to zero;
    \item for each $k$-simplex $\sigma:\DE^k\to K_\bullet$ with $k > 0$, define composites $\sigma_{\leq j}:\DE^j\mono \DE^k\xrar{\sigma} K_\bullet$ and $\sigma_{\geq j}:\DE^{k-j}\mono \DE^k\xrar{\sigma} K_\bullet$ for each $0 < j < k$, where $\DE^{k-j}\mono \DE^k$ is given on vertices by $i\mapsto i+j$, then we have 
    \begin{equation}
	\label{map into dg nerve}
 	\del \psi(\sigma) = \sum_{j = 1}^{k-1} (-1)^{k-j}(\psi(\sigma_{\geq j})\circ \psi(\sigma_{\leq j})) - \psi(d_j^k\sigma))
 \end{equation}
     where $d_j^k$ are the face maps of $K_\bullet$.
\end{enumerate}

We also need a more explicit description of $J$: for each $k > 0$, the $\infty$-groupoid $J$ has exactly two non-degenerate $k$-simplices, namely $\sigma_k : \cdots\to 0\to 1\to 0$ and $\tau_k:\cdots\to 1\to 0\to 1$. The face maps are given by \begin{equation}
\label{J face maps}
d_j^k(\sigma_k) = \begin{cases}\sigma_{k-1} & j = 0\\
 \tau_{k-1} & j = k	\\
 \text{degenerate } & \text{otherwise}
 \end{cases}\,,\quad d_j^k(\tau_k) = \begin{cases}\tau_{k-1} & j = 0\\
 \sigma_{k-1} & j =k\\
 \text{degenerate } & \text{otherwise}	
 \end{cases}\,.
\end{equation}
We are now ready to prove:
\begin{thm}
\label{quasi-inverse}
The assignment \begin{equation}
\label{quasi inverse element}\begin{cases}
0\mapsto \Aug_+[1]\\
 1\mapsto \mcalQ_-^\dagger[-1]\\
 \sigma_k\mapsto f_k & k\geq 1\\
 \tau_k\mapsto g_k & k\geq 1
 \end{cases}\end{equation}
defines a map of $\infty$-categories
$\psi:J\to N^\mrm{dg}_\bullet(\bimod{\Aug_+})$.
\end{thm}

\begin{proof} 
Conditions (1) and (2) determine $\psi$ on degenerate simplicies. 
Condition (3) is easy to check for degenerate simplicies, and we will leave it to the reader. It remains to check \eqref{map into dg nerve} for $\sigma_k$'s and $\tau_k$'s. Note that $(\sigma_k)_{\geq j} = \sigma_{k-j}$, $(\tau_k)_{\geq j} = \tau_{k-j}$, \[(\sigma_k)_{\leq j} = \begin{cases}\sigma_j & k-j\text{ even}\\
 \tau_j & k-j\text{ odd}	
 \end{cases}\,,\quad (\tau_k)_{\leq j} = \begin{cases}\tau_j & k-j\text{ even}\\
 \sigma_j & k-j\text{ odd}	
 \end{cases}\,.
\] Since $d_j^k\sigma_k$ and $d_j^k\tau_k$ are degnerate for each $1\leq j\leq k-1$, $\psi$ sends them to identities if $ k = 2$, and zero otherwise. Hence, for $k = 2$, we have
\[-(\psi((\sigma_2)_{\geq 1})\circ \psi((\sigma_2)_{\leq 1})) - \psi(d_1^2\sigma_2)) = -f_1g_1 + 1_{\Aug_+[1]}= \del\psi(\sigma_2) \] by \eqref{homotopy equations}.
For  $k > 2$, we have
\begin{align*}&\quad\,\sum_{j = 1}^{k-1} (-1)^{k-j}(\psi((\sigma_k)_{\geq j})\circ \psi((\sigma_k)_{\leq j})) - \psi(d_j^k\sigma_k))\\& = -f_1g_{k-1} + f_2f_{k-2}-\cdots = \del\psi(\sigma_k)
\end{align*}	
 by \eqref{J face maps} and \eqref{homotopy equations}. The argument is similar for $\tau_k$'s.
\end{proof}
Since each augmentation canonically lifts to a $0$-curved augmentation by \Cref{0-curved augmentations}, picking $c_1 = 0$ and $c_2 = 1$ gives us:
\begin{thm}
\label{higher homotopy thm}
Let $\LA$ be a Legendrian knot of rotation number zero, equipped with a single marked point and any Morse function. Assume every augmentation of $\LA$ lifts to a $1$-curved augmentation. Then there is an explicit map of $\infty$-categories \[\psi:J\to N^\mrm{dg}_\bullet(\bimod{\Aug_+})\]specified by \eqref{quasi inverse element}, which sends the edge $1\to 0$ to the $A_\infty$ Sabloff map $f: \mcalQ^\dagger_-[-1]\to \Aug_+[1]$.	
\end{thm}
\subsection{Weak Calabi-Yau Structures}
\label{CY structures}
This subsection centers around \Cref{CY conj} below. We begin by
defining weak (relative) proper Calabi-Yau structures of $A_\infty$ categories. Although the definitions usually involve (relative) Hochschild complexes \cite{BD19}\cite{Gan23}, we have chosen the equivalent phrasing in terms of bimodules, since it works more directly with our constructions.

\begin{defn}
An $A_\infty$ category is \tit{proper} if all of its hom-spaces are perfect $\k$-complexes.    
\end{defn}

Let $\mcalA,\mcalC$ be strictly unital proper $A_\infty$ categories and $f:\mcalA\to\mcalC$ an $A_\infty$ functor.

\begin{defn}
	A \tit{weak proper $n$-dimensional Calabi-Yau ($n$-CY) structure} on $\mcalA$ is a quasi-isomorphism of $\mcalA$-bimodules \[\mcalA[n]\iso\mcalA^\dagger\,.\]
\end{defn}
\begin{defn}
\label{weak relative CY def}
A \tit{weak relative proper $n$-CY structure} on $f$ consists of
\begin{itemize}
	\item a weak proper ($n-1$)-CY structure $v$ on $\mcalC$
	\item a degree $-1$ pre-morphism $u:\mcalA[n-1]\to \mcalA^\dagger$
\end{itemize}
such that 
\begin{enumerate}[label = (\arabic*)]
	\item $u$ is a null-homotopy for the composition
\[
\begin{tikzcd}
\mcalA[n-1]\rar{f} & f^*\mcalC[n-1]\rar{v} & \mcalC^\dagger \rar{f^\dagger} & \mcalA^\dagger 	
\end{tikzcd}
\]    	
\item the induced map of exact triangles in the derived category $D(\bimod{\mcalA})$ \[\begin{tikzcd}        \mcalA[n-1]\rar{f}\dar{} & f^*\mcalC[n-1]\rar{}\dar{v} &\cone(f)\rar\dar{} & {}\\
        \cone(f^\dagger)[-1]\rar{} &(f^*\mcalC)^\dagger\rar{f^\dagger} & \mcalA^\dagger\rar & {}
    \end{tikzcd}\]
    is an isomorphism.
\end{enumerate}
\end{defn}
\begin{conj}
\label{CY conj}
    For any Legendrian link $\LA$ with marked points and any Morse function, the strict $A_\infty$ functor $\pi:\Aug_+\to \mcalC$ admits a weak relative proper $2$-CY structure.
\end{conj}
\begin{rmk}
    The upcoming work of Ma and Sabloff \cite{MS25} will establish this conjecture over $\Z/2$, although the weak CY structure they produce is non-explicit. The author believes $\pi$ actually admits a (strong) CY structure, in the sense of Brav and Dyckerhoff \cite{BD19}.

    Legout \cite{Leg23} showed that the Chekanov-Eliashberg dga of some Legendrian spheres admits a weak smooth CY structure. In her joint work with Dimitroglou Rizell \cite{DRL23}\cite{DRL24}, this result will be extended to a weak relative smooth CY structure on the inclusion of coefficients into the Chekanov-Eliashberg dga with loop space coefficients. Kuo and Li \cite{KL24} showed the existence of a strong relative smooth CY structure on the microlocalization functor of constructible sheaves with Legendrian singular support. It would be interesting to see whether these various CY structures are compatible with each other.
\end{rmk}
 Note that $\Aug_+$ and $\mcalC$ are proper because their hom-spaces are bounded complexes of finite free $\k$-modules. Below we present partial results towards the above conjecture, starting with a weak proper $1$-CY structure for $\mcalC$. 

For the rest of this subsection, let $\LA$ be a Legendrian knot of rotation number zero, equipped with a single marked point and a Morse function with exactly two critical points, and assume they are positioned as in \Cref{matrix formula prop}. Using \eqref{C for knot}, it is easy to verify the following:
\begin{thm}
\label{CY on circle}
     There is a strict morphism of $\mcalC$-bimodules $v:\mcalC[1]\to\mcalC^\dagger$ given by
\[\nu_{0|0}(sx^+) = (y^+)^\dagger,\ \nu_{0|0}(sy^+) = (x^+)^\dagger\]
which has an obvious strict inverse.
\end{thm}
Consider the following diagram of $\Aug_+$-bimodules:
\begin{equation}
\label{CY square}
    \begin{tikzcd}[column sep = 50]
        \Aug_+[1]\rar{\pi} & \pi^*\mcalC[1]\dar{v} \\
        \mcalQ_-^\dagger[-1]\uar{f}\rar{a^\dagger[-1]} &(\pi^*\mcalC)^\dagger 
    \end{tikzcd}
\end{equation}
and let $c = v\pi f -  a^\dagger[-1]$ denote the difference. It is straightforward to calculate $c$ explicitly, using \Cref{LSFT matrix formula prop}. For cusp-pointed Legendrian knots, we get the following: 
\begin{lemma}
\label{left square difference}
Assume $\LA$ is cusp-pointed and the boundary of the disk at the resolved cusp containing the marked point agrees with the orientation of $\LA$. Then the only nonzero un-suspended components of $c$ are given by
\begin{equation}
	\label{wrong map difference}
	\begin{cases}
        \gamma_{1|0}(x^+,(q_\bullet^-)^\dagger) = (x^+)^\dagger\\
		\gamma_{1|0}(x^+,(q_k^-)^\dagger) =  (\ep_1(q_k) - \ep_2(q_k)) (y^+)^\dagger\\
		\gamma_{1|1}(x^+,(q_k^-)^\dagger,q_k^+) = -\gamma_{2|0}(x^+,q_k^+,(q_k^-)^\dagger) = (y^+)^\dagger\\
	\end{cases}\hspace{-1ex}.
\end{equation}
In particular, the linear term $\gamma_{0|0}$ is zero, so \eqref{CY square} commutes on cohomology.
\end{lemma}

 Although $c_{0|0} = 0$, we could not show $c$ itself is null-homotopic. We suspect the reason is that our choice of the absolute CY structure $v$ of $\mcalC$, although simple, is too artificial. A more natural choice should come from the continuation map for a homotopy connecting the perturbing Morse function to its negative. If one can find such a $v$ so that \eqref{CY square} commutes up to homotopy, then we have the following recipe for an explicit weak relative CY structure on $\pi:\Aug_+\to \mcalC$. Write $\Aug_+^\dagger = \mcalQ_-^\dagger \oplus \mcalC^\dagger$ and let $\iota: \mcalQ_-^\dagger\to \Aug_+^\dagger$ be the inclusion. 

\begin{prop}
\label{CY prop}
    Assume every augmentation of $\LA$ lifts to a $1$-curved augmentation, and assume $c$ admits a null-homotopy $r$. Let $g\defeq g_1$ and $f_2$ be the maps defined in \Cref{higher homotopies}, and define the degree $-1$ pre-morphism $u: \Aug_+[1]\to\Aug_+^\dagger$ by 
\[u \defeq \iota g + \pi^\dagger rg + \pi^\dagger v\pi f_2 \,.\]
    Then the pair $(v,u)$ defines a weak relative proper $2$-CY structure on $\pi:\Aug_+\to \mcalC$.
\end{prop}
\begin{proof}
We need to check conditions (1) and (2) in \Cref{weak relative CY def}.

    By \Cref{Q+- triangle}, we have $b^{\Aug_+^\dagger}\iota+\iota b^{\Aug_+[1]} =\pi^\dagger a^\dagger[-1]$. By \Cref{homotopy equations lemma}, we have $b^{\Aug_+^\dagger}f_2 + f_2b^{\Aug_+[1]} = 1-fg$, where $f$ is the $A_\infty$ Sabloff map. Also $b^{\Aug_+^\dagger}r+rb^{\Aug_+[1]} = c$ by assumption. Therefore, we get 
    \begin{align*}
        &\quad\, b^{\Aug_+^\dagger}u+u b^{\Aug_+[1]}\\ &= (b^{\Aug_+^\dagger}\iota + \iota b^{\Aug_+[1]})g + \pi^\dagger(b^{\Aug_+^\dagger}r+rb^{\Aug_+[1]})g +\pi^\dagger v\pi(b^{\Aug_+^\dagger}f_2 + f_2b^{\Aug_+[1]})  \\
        &= \pi^\dagger a^\dagger[-1] g+\pi^\dagger(v\pi f - a^\dagger[-1])g + \pi^\dagger v\pi (1-fg) = \pi^\dagger v\pi
    \end{align*} so condition (1) is satisfied.
    
    For condition (2), we need to check that the induced map of exact triangles \[\begin{tikzcd}
        \Aug_+[1]\rar{\pi}\dar{} & \pi^*\mcalC[1]\rar{}\dar{v} &\cone(\pi)\rar\dar{} & {}\\
        \cone(\pi^\dagger)[-1]\rar{} &(\pi^*\mcalC)^\dagger\rar{\pi^\dagger} & \Aug_+^\dagger\rar & {}
    \end{tikzcd}\] is an isomorphism in $D(\bimod{\Aug_+})$. Since $v$ is a strict isomorphism, it suffices to show the first vertical map $\Aug_+[1]\to \cone (\pi^\dagger)[-1] = (\pi^*\mcalC)^\dagger\oplus\Aug_+^\dagger[-1]$, given by $(v\pi, u)$, is an isomorphism in $D(\bimod{\Aug_+})$. Recall from \Cref{negative aug bimod} that $\Aug_+$ is the mapping cone of a morphism $a: \pi^*\mcalC[-1]\to\mcalQ_-$, so $\Aug_+^\dagger$ is the mapping cone of $a^\dagger[-1]: \mcalQ_-^\dagger[-1]\to (\pi^*\mcalC)^\dagger$. Thus, there is a canonical quasi-isomporphism $\cone (\pi^\dagger)\to \mcalQ_-^\dagger$, given by quotienting out $(\pi^*\mcalC)^\dagger$ and the image of $\pi^\dagger:(\pi^*\mcalC)^\dagger\to\Aug_+^\dagger$. By definition of $u$, it is easy to see that the composition \[\begin{tikzcd}
        \Aug_+[1]\arrow[rr, "{(v\pi, u)}"] & & \cone (\pi^\dagger)[-1]\rar &  \mcalQ_-^\dagger[-1] 
    \end{tikzcd}\]is just $g: \Aug_+[1]\to\mcalQ_-^\dagger[-1]$, which is a homotopy equivalence. Hence, $\Aug_+[1]\to \cone(\pi^\dagger)[-1]$ is an isomorphism in $D(\bimod{\Aug_+})$.
\end{proof}

\newpage
\bibliographystyle{amsalpha}
\bibliography{references.bib}

\end{document}